\documentclass[a4, 12pt, leqno]{article}
\usepackage{amsmath}
\usepackage{amssymb}
\newtheorem{df}{Definition}[section]
\newtheorem{pr}[df]{Proposition}
\newtheorem{Th}[df]{Theorem}

\newtheorem{lm}[df]{Lemma}

\begin{document}

	\title{\bf Motion of a Vortex Filament \\ on a Slanted Plane}
\author{Masashi A{\sc iki}}
\date{}
\maketitle
\vspace*{-0.5cm}

\begin{abstract}
We consider a nonlinear model equation, known as the 
Localized Induction Equation, describing the motion of a vortex filament
immersed in an incompressible and inviscid fluid. We prove the unique solvability
of an initial-boundary value problem describing the motion of a vortex filament on a slanted plane.
\end{abstract}

\section{Introduction and Problem Setting}

A vortex filament is a space curve on which the vorticity of the fluid is concentrated. 
Vortex filaments are used to model very thin vortex structures such as vortices that 
trail off airplane wings or propellers. 
In this paper, we prove the solvability of the following initial-boundary value problem
which describes the motion of a vortex filament moving on a slanted plane.
\begin{eqnarray}
\left\{
\begin{array}{ll}
\mbox{\mathversion{bold}$x$}_{t} =\mbox{\mathversion{bold}$x$}_{s}\times 
\mbox{\mathversion{bold}$x$}_{ss}, & s\in I, \ t>0, \\[3mm]
\mbox{\mathversion{bold}$x$}(s,0)=\mbox{\mathversion{bold}$x$}_{0}(s), & s\in I, \ t>0, \\[3mm]
\mbox{\mathversion{bold}$x$}_{s}(0,t)=\mbox{\mathversion{bold}$a$}, \ 
\mbox{\mathversion{bold}$x$}_{s}(1,t)=\mbox{\mathversion{bold}$e$}_{3}, & t>0,
\end{array}\right.
\label{slant}
\end{eqnarray}
where \( \mbox{\mathversion{bold}$x$}(s, t)= {}^{t}(x_{1}(s ,t), x_{2}(s ,t), x_{3}(s ,t))\) is the 
position vector of the vortex filament parametrized by 
its arc length \( s \) at time \( t\), 
\( \times \) is the exterior product in the three dimensional Euclidean space,
\( I=(0,1)\subset \mathbf{R}\) is an open interval, 
\( \mbox{\mathversion{bold}$a$}\in \mathbf{R}^{3}\) is an arbitrary vector satisfying 
\( |\mbox{\mathversion{bold}$a$}|=1 \), \( \mbox{\mathversion{bold}$e$}_{3}={}^{t}(0,0,1) \),
and subscripts \( s\) and \(t \) are differentiations with the respective variables.
Problem (\ref{slant}) describes the motion of a segment of a vortex filament moving on a slanted plane.
We can see that, by taking the trace \( s=0\) in the equation of problem
(\ref{slant}), a filament moving according to problem (\ref{slant}) satisfies
\begin{align*}
\mbox{\mathversion{bold}$x$}_{t}(0,t)&=\mbox{\mathversion{bold}$a$}\times \mbox{\mathversion{bold}$x$}_{ss}(0,t),
\end{align*}
hence the end-point \( \mbox{\mathversion{bold}$x$}(0,t)\)
of the filament moves along the plane perpendicular to \( \mbox{\mathversion{bold}$a$}\).
The reason we also impose a boundary condition at \( s=1 \) is for the following reason. 
A more intuitive problem setting for a vortex filament moving on a plane would be
\begin{eqnarray}
\left\{
\begin{array}{ll}
\mbox{\mathversion{bold}$x$}_{t} =\mbox{\mathversion{bold}$x$}_{s}\times 
\mbox{\mathversion{bold}$x$}_{ss}, & s>0, \ t>0, \\[3mm]
\mbox{\mathversion{bold}$x$}(s,0)=\mbox{\mathversion{bold}$x$}_{0}(s), & s>0, \ t>0, \\[3mm]
\mbox{\mathversion{bold}$x$}_{s}(0,t)=\mbox{\mathversion{bold}$a$}, & t>0,
\end{array}\right.
\label{hp}
\end{eqnarray}
which is a problem describing an infinitely long filament with one end moving along the plane perpendicular to 
\( \mbox{\mathversion{bold}$a$}\). The solvability of problem (\ref{hp}) is a direct consequence of 
a previous work by the author and Iguchi \cite{11}, which proved the solvability of problem 
(\ref{hp}) with \( \mbox{\mathversion{bold}$a$}=\mbox{\mathversion{bold}$e$}_{3} \), because the solution of
problem (\ref{hp}) can be obtained by rotating the solution obtained in \cite{11} in a way that 
\( \mbox{\mathversion{bold}$a$} \) is trasformed to \( \mbox{\mathversion{bold}$e$}_{3} \). Hence,
problem (\ref{hp}) for general \( \mbox{\mathversion{bold}$a$} \) is essentially the same as the case 
\( \mbox{\mathversion{bold}$a$}=\mbox{\mathversion{bold}$e$}_{3} \). 
So to describe the motion of a vortex filament on a slanted plane, we imposed a boundary condition at 
\( s=1 \) to set a reference plane which allows us to express the slanted-ness of the plane that the 
filament is moving on.
The motivation for considering problem (\ref{slant}) comes from the following 
problem.
\begin{align}
\left\{
\begin{array}{ll}
\label{gb}
\mbox{\mathversion{bold}$x$}_{t}=\mbox{\mathversion{bold}$x$}_{s}\times 
\mbox{\mathversion{bold}$x$}_{ss}, & s>0, \ t>0, \\[3mm]
\mbox{\mathversion{bold}$x$}(s,0)= \mbox{\mathversion{bold}$x$}_{0}(s), & s>0, \\[3mm]
\displaystyle 
\mbox{\mathversion{bold}$x$}_{s}(0,t)=
\frac{\nabla B(\mbox{\mathversion{bold}$x$}(0,t))}{|\nabla B(\mbox{\mathversion{bold}$x$}(0,t))|}, & t>0,
\end{array}\right.
\end{align}
where \( \nabla = {}^{t}(\partial_{x_{1}},\partial_{x_{2}},
\partial_{x_{3}})\) and 
for \( {}^{t}(x_{1},x_{2},x_{3})\in \mathbf{R}^{3} \), 
\( B:\mathbf{R}^{3}\to \mathbf{R}\) is a given function of the form
\begin{align*}
B(x_{1},x_{2},x_{3}) = x_{3}-b(x_{2},x_{3})
\end{align*}
for a given scalar function \( b:\mathbf{R}^{2}\to \mathbf{R}\). 
Problem (\ref{gb}) describes an infinitely long vortex filament moving on a 
surface given as the graph of \( b\) in the three-dimensional Euclidean space.
Problem (\ref{gb}) is a generalization of the problem setting in \cite{11} and 
can be seen as a simplified model for the motion of a tornado, where the 
ground is given by the graph of \( b\), but the solvability for a general \( b\) 
seems hard, and as a first step, we chose the special case where the ground is a 
slanted plane.

The equation in problem (\ref{slant}) is called the 
Localized Induction Equation (LIE) which is derived by applying the so-called localized induction approximation to the 
Biot--Savart integral.
The LIE was first derived by 
Da Rios in 1906 and was re-derived twice independently by Murakami et al. in 1937 and by 
Arms and Hama in 1965. Many researches have been done on the LIE and many results have been obtained.
Nishiyama and Tani \cite{5,13} proved the unique solvability of the initial value problem in Sobolev 
spaces. Koiso considered a geometrically generalized setting in which he proved rigorously the 
equivalence of the LIE and a nonlinear Schr\"odinger equation. This equivalence was first shown by 
Hasimoto \cite{14} in which he studied the formation of solitons on a vortex filament. He defined a transformation of variable known as the Hasimoto transformation to transform the LIE into a
nonlinear Schr\"odinger equation. 
The Hasimoto transformation was proposed by Hasimoto \cite{14}
and is a change of variable given by 
\[ \psi = \kappa \exp \left( {\rm i} \int ^{s}_{0} \tau \,{\rm d}s \right), \]
where \( \kappa \) is the curvature and \( \tau \) is the torsion of the filament.
Defined as such, it is well known that \( \psi \) satisfies the 
nonlinear Schr\"odinger equation given by 
\begin{eqnarray}
{\rm i}\frac{\partial \psi}{\partial t}
= \frac{ \partial ^{2}\psi }{\partial s^{2}} + \frac{1}{2} \left| \psi \right| ^{2}\psi .
\label{NLS}
\end{eqnarray}
The original transformation proposed by Hasimoto uses the torsion of the filament in its definition,
which means that the transformation is undefined at points where the curvature of the filament is zero. 
Koiso \cite{12} constructed a transformation, sometimes referred to as the 
generalized Hasimoto transformation, and gave a mathematically rigorous 
proof of the equivalence of the LIE and (\ref{NLS}). More recently, 
Banica and Vega \cite{15,16} and Guti\'errez, Rivas, and Vega \cite{17} constructed and analyzed
a family of self-similar solutions of the LIE which forms a corner in finite time. 
The authors \cite{11} proved the unique solvability of an initial-boundary value problem for the LIE
in which the filament moved in the three-dimensional half space. Nishiyama and Tani \cite{5} also
considered initial-boundary value problems with different boundary conditions. 

Setting \( \mbox{\mathversion{bold}$v$}:= \mbox{\mathversion{bold}$x$}_{s}\) and taking the 
\( s\)-derivative of the equation in problem (\ref{slant}), we see that problem (\ref{slant}) is trasformed to
\begin{eqnarray}
\left\{
\begin{array}{ll}
\mbox{\mathversion{bold}$v$}_{t} =\mbox{\mathversion{bold}$v$}\times 
\mbox{\mathversion{bold}$v$}_{ss}, & s\in I, \ t>0, \\[3mm]
\mbox{\mathversion{bold}$v$}(s,0)=\mbox{\mathversion{bold}$v$}_{0}(s), & s\in I, \ t>0, \\[3mm]
\mbox{\mathversion{bold}$v$}(0,t)=\mbox{\mathversion{bold}$a$}, \ 
\mbox{\mathversion{bold}$v$}(1,t)=\mbox{\mathversion{bold}$e$}_{3}, & t>0,
\end{array}\right.
\label{vslant}
\end{eqnarray}
where \( \mbox{\mathversion{bold}$v$}_{0}:= \mbox{\mathversion{bold}$x$}_{0s} \).
If we can prove the solvability of problem (\ref{vslant}), the solution \( \mbox{\mathversion{bold}$x$} \)
of problem (\ref{slant}) can be constructed from \( \mbox{\mathversion{bold}$v$} \) by 
\begin{eqnarray*}
\mbox{\mathversion{bold}$x$}(s,t) = \mbox{\mathversion{bold}$x$}_{0}(s) + \int^{t}_{0} 
\mbox{\mathversion{bold}$v$}(s,\tau )\times \mbox{\mathversion{bold}$v$}_{ss}(s, \tau ) {\rm d}\tau ,
\end{eqnarray*}
hence the solvability of problem (\ref{slant}) and (\ref{vslant}) are equivalent.
To this end, we consider the solvability of problem (\ref{vslant}) from here on.

The contents of the rest of the paper are as follows. In Section 2, we define notations used 
in this paper and state our main theorem. In Section 3, we introduce a regularized nonlinear 
problem and construct a corrected initial datum associated to the 
regularized problem. The correction is necessary to insure that the 
compatibility conditions for the regularized problem are satisfied.
In Section 4, we give a brief description of the method used to prove the solvability of the 
linear problem associated to the regularized nonlinear problem given in Section 3 and 
state the existence theorem for the regularized nonlinear problem.
Finally, in Section 5 we construct the solution of (\ref{vslant}), and
derive estimates of the solution to prove the time-global 
solvability of (\ref{vslant}).


\section{Function Spaces, Notations, and Main Theorem}

We introduce some function spaces that will be used throughout this paper, and notations associated with the spaces.
For a non-negative integer \( m\), and \( 1\leq p \leq \infty \), \( W^{m,p}(I)\) is the Sobolev space 
containing all real-valued functions that have derivatives in the sense of distribution up to order \( m\) 
belonging to \( L^{p}(I)\).
We set \( H^{m}(I) := W^{m,2}(I) \) as the Sobolev space equipped with the usual inner product.
The norm in \( H^{m}(I) \) is denoted by \( \| \cdot \|_{m} \) and we simply write \( \| \cdot \| \) 
for \( \|\cdot \|_{0} \). Otherwise, for a Banach space \( X\), the norm in \( X\) is written as \( \| \cdot \| _{X}\).
The inner product in \( L^{2}(I)\) is denoted by \( (\cdot ,\cdot )\).

For \( 0<T \leq \infty \) and a Banach space \( X\), 
\( C^{m}([0,T];X) \)
( \( C^{m}\big( [0,\infty);X\big)\) when \( T= \infty \)),
denotes the space of functions that are \( m\) times continuously differentiable 
in \( t\) with respect to the norm of \( X\).

For any function space described above, we say that a vector valued function belongs to the function space 
if each of its components does. 

\vspace*{5mm}
We further define notations to express the 
compatibility conditions for problem (\ref{vslant}). 
First we set 
\( \mbox{\mathversion{bold}$P$}_{0}(\mbox{\mathversion{bold}$v$}):=
\mbox{\mathversion{bold}$v$}\) and
\( \mbox{\mathversion{bold}$P$}_{1}(\mbox{\mathversion{bold}$v$}):=
\mbox{\mathversion{bold}$v$}\times \mbox{\mathversion{bold}$v$}_{ss}\).
\( \mbox{\mathversion{bold}$P$}_{1}(\mbox{\mathversion{bold}$v$})\) is the 
right-hand side of the equation in (\ref{vslant}). We also use the notation
\( \mbox{\mathversion{bold}$P$}_{1}(s,t):= \mbox{\mathversion{bold}$P$}_{1}
(\mbox{\mathversion{bold}$v$})\) and sometimes omit \( (s,t) \) for simplicity. 
We successively define
\begin{eqnarray*}
\mbox{\mathversion{bold}$P$}_{m}:=
\sum ^{m-1}_{j=0}
\left(
\begin{array}{c}
m-1\\
j
\end{array}
\right)
\mbox{\mathversion{bold}$P$}_{j}\times \partial ^{2}_{s}\mbox{\mathversion{bold}$P$}_{m-1-j},
\end{eqnarray*}
where 
\( 
\left(
\begin{array}{c}
m-1\\
j
\end{array}
\right)
\)
is the binomial coeffecient.
\( \mbox{\mathversion{bold}$P$}_{m}\) gives the expression of 
\( \partial ^{m}_{t}\mbox{\mathversion{bold}$v$} \) with only \( s\)-derivatives of 
\( \mbox{\mathversion{bold}$v$}\). From the boundary condition for 
problem (\ref{vslant}),  we arrive at the following definitions for the 
compatibility conditions.
\begin{df}(Compatibility conditions for (\ref{vslant})).
For \( m\in \mathbf{N}\cup \{ 0\}\), we say that \( \mbox{\mathversion{bold}$v$}_{0}\)
satisfies the \( m\)-th order compatibility condition for {\rm (\ref{vslant})} if
\( \mbox{\mathversion{bold}$v$}_{0}\in H^{2m+1}(I)\) and 
\begin{eqnarray*}
\mbox{\mathversion{bold}$v$}_{0}(0) = \mbox{\mathversion{bold}$a$}, \quad 
\mbox{\mathversion{bold}$v$}_{0}(1) = \mbox{\mathversion{bold}$e$}_{3},
\end{eqnarray*}
when \( m=0\), and 
\begin{eqnarray*}
\mbox{\mathversion{bold}$P$}_{m}(\mbox{\mathversion{bold}$v$}_{0}(0))
=
\mbox{\mathversion{bold}$P$}_{m}(\mbox{\mathversion{bold}$v$}_{0}(1)) 
=
\mbox{\mathversion{bold}$0$}
\end{eqnarray*}
when \( m \geq 1 \). We also say that \( \mbox{\mathversion{bold}$v$}_{0}\) satisfies 
the compatibility conditions for {\rm (\ref{vslant})} up to order \( m\) if 
\( \mbox{\mathversion{bold}$v$}_{0}\) satisfies the \( k\)-th order 
compatibility condition for all \( k\) with \( 0\leq k \leq m \).
\label{cco}
\end{df}

\vspace*{1cm}

Now we state our main theorem regarding the solvability of (\ref{vslant}). 
\begin{Th}
For an integer \( l\geq 0\), if \( \mbox{\mathversion{bold}$v$}_{0}
\in H^{l+3}(I)\), 
\( |\mbox{\mathversion{bold}$v$}_{0}|\equiv 1\), and 
\( \mbox{\mathversion{bold}$v$}_{0}\) satisfy the compatibility conditions 
up to order \( [\frac{l+3}{2}]\), 
then there exists a unique solution \( \mbox{\mathversion{bold}$v$}\)
satisfying \( | \mbox{\mathversion{bold}$v$}|\equiv 1 \) and 
\begin{eqnarray*}
\mbox{\mathversion{bold}$v$}\in \bigcap ^{[\frac{l+3}{2}]}_{j=0}
C^{j}\big( [0,\infty);H^{l+3-2j}(I) \big).
\end{eqnarray*}
Here, \( [\frac{l+3}{2}]\) is the largest integer not exceeding 
\( \frac{l+3}{2}\).
\label{TH}
\end{Th}
The above theorem gives the time-global unique solvability of problem (\ref{vslant}) and thus, 
for (\ref{slant}). Note that if the initial datum \( \mbox{\mathversion{bold}$v$}_{0}\)
satisfies \( | \mbox{\mathversion{bold}$v$}_{0}|\equiv 1 \), then a solution
\( \mbox{\mathversion{bold}$v$}\) of 
(\ref{vslant}) also satisfies \( |\mbox{\mathversion{bold}$v$}|\equiv 1 \) automatically.
This is because
\begin{eqnarray*}
\frac{{\rm d}}{{\rm d}t}|\mbox{\mathversion{bold}$v$}|^{2}
= 2\mbox{\mathversion{bold}$v$}\cdot \mbox{\mathversion{bold}$v$}_{t}
= 2\mbox{\mathversion{bold}$v$}\cdot ( 
\mbox{\mathversion{bold}$v$}\times \mbox{\mathversion{bold}$v$}_{ss})
=0,
\end{eqnarray*}
and the arc length parameter is preserved throughout the motion.
Here, \( \cdot \) is the inner product in the three-dimensional Euclidean space.
This property will play a crucial role in the upcoming analysis.

%
%
%
%
%
%
%
\section{Regularized nonlinear problem and its compatibility conditions}
\setcounter{equation}{0}
We construct the solution of problem (\ref{vslant}) by taking the limit 
\( \varepsilon \to +0 \) in the following regularized problem.
\begin{eqnarray}
\left\{
\begin{array}{ll}
\mbox{\mathversion{bold}$v$}^{\varepsilon}_{t}=
\mbox{\mathversion{bold}$v$}^{\varepsilon}\times \mbox{\mathversion{bold}$v$}^{\varepsilon}_{ss}
+ \varepsilon \mbox{\mathversion{bold}$v$}^{\varepsilon}_{ss} 
+ \varepsilon |\mbox{\mathversion{bold}$v$}^{\varepsilon}_{s}
|^{2}\mbox{\mathversion{bold}$v$}^{\varepsilon},
& s\in I, \ t>0, \\[3mm]
\mbox{\mathversion{bold}$v$}^{\varepsilon}(s,0)= \mbox{\mathversion{bold}$v$}^{\varepsilon}_{0}
(s), & s\in I, \ t>0, \\[3mm]
\mbox{\mathversion{bold}$v$}^{\varepsilon}(0,t)= \mbox{\mathversion{bold}$a$}, \ 
\mbox{\mathversion{bold}$v$}^{\varepsilon}(1,t)=\mbox{\mathversion{bold}$e$}_{3},
& t>0.
\end{array}\right.
\label{rnl}
\end{eqnarray}
From here on, it is assumed that
\( |\mbox{\mathversion{bold}$v$}_{0}|\equiv 1 \) holds, i.e. the initial datum 
for the original problem (\ref{vslant}) is 
paramtrized by its arc length. 
The regularization shown above was chosen for two main reasons.
Firstly, \( \varepsilon \mbox{\mathversion{bold}$v$}^{\varepsilon}_{ss}\) was added 
so that the associated linear equation becomes a parabolic equation. Secondly, 
the term \( \varepsilon |\mbox{\mathversion{bold}$v$}^{\varepsilon}_{s}|^{2}
\mbox{\mathversion{bold}$v$}^{\varepsilon} \) was added so that if the initial datum 
\( \mbox{\mathversion{bold}$v$}^{\varepsilon}_{0} \) satisfies 
\( | \mbox{\mathversion{bold}$v$}^{\varepsilon}_{0}|\equiv 1 \), then a solution
of (\ref{rnl}) also satisfies \( | \mbox{\mathversion{bold}$v$}^{\varepsilon}| \equiv 1 \).

Since we modified the equation, we must make corrections to the initial datum to 
insure that the compatibility conditions for problem 
(\ref{rnl}) are satisfied.
\subsection{Compatibility conditions for (\ref{rnl})}
We first derive the compatibility conditions for problem (\ref{rnl}).
We set \( \mbox{\mathversion{bold}$Q$}_{0}(\mbox{\mathversion{bold}$v$}):= 
\mbox{\mathversion{bold}$v$}\),
\( \mbox{\mathversion{bold}$Q$}_{1}(\mbox{\mathversion{bold}$v$}):=
\mbox{\mathversion{bold}$v$}\times \mbox{\mathversion{bold}$v$}_{ss} + 
\varepsilon \mbox{\mathversion{bold}$v$}_{ss} + 
\varepsilon |\mbox{\mathversion{bold}$v$}_{s}|^{2}\mbox{\mathversion{bold}$v$}\), and
\begin{align*}
\mbox{\mathversion{bold}$Q$}_{m}(\mbox{\mathversion{bold}$v$})&:=
\sum ^{m-1}_{j=0}
\left(
\begin{array}{c}
m-1\\
j
\end{array}
\right)
\mbox{\mathversion{bold}$Q$}_{j}\times 
\partial ^{2}_{s}\mbox{\mathversion{bold}$Q$}^{\varepsilon}_{m-1-j}
+
\varepsilon \partial^{2}_{s}\mbox{\mathversion{bold}$Q$}_{m-1}
(\mbox{\mathversion{bold}$v$}) \\[3mm]
& \qquad + \sum ^{m-1}_{j=0}
\sum^{m-1-j}_{k=0}
\left(
\begin{array}{c}
m-1\\
j
\end{array}
\right)
\left(
\begin{array}{c}
m-1-j\\
k
\end{array}
\right)
(\partial_{s}\mbox{\mathversion{bold}$Q$}_{j}\cdot 
\partial_{s}\mbox{\mathversion{bold}$Q$}_{k})
\mbox{\mathversion{bold}$Q$}_{m-1-j-k}
\end{align*}
for \( m\geq 2 \). We arrive at the following compatibility conditions.
\begin{df}(Compatibility conditions for (\ref{rnl})).
For \( m\in \mathbf{N}\cup \{0\}\), we say that 
\( \mbox{\mathversion{bold}$v$}^{\varepsilon}_{0}\) satisfies the \( m\)-th order 
compatibility condition for (\ref{rnl}) if 
\( \mbox{\mathversion{bold}$v$}^{\varepsilon}_{0}\in H^{2m+1}(I) \) and
\begin{eqnarray*}
\mbox{\mathversion{bold}$v$}^{\varepsilon}_{0}(0)=
\mbox{\mathversion{bold}$a$}, \quad 
\mbox{\mathversion{bold}$v$}^{\varepsilon}_{0}(1)=
\mbox{\mathversion{bold}$e$}_{3},
\end{eqnarray*}
when \( m=0 \), and 
\begin{eqnarray*}
\mbox{\mathversion{bold}$Q$}_{m}(
\mbox{\mathversion{bold}$v$}^{\varepsilon}_{0}(0)) =
\mbox{\mathversion{bold}$Q$}_{m}(
\mbox{\mathversion{bold}$v$}^{\varepsilon}_{0}(1)) = 0,
\end{eqnarray*}
when \( m\geq 1 \).
We also say that \( \mbox{\mathversion{bold}$v$}^{\varepsilon}_{0}\) satisfies the 
compatibility conditions for (\ref{rnl}) up to order \( m \) if 
\( \mbox{\mathversion{bold}$v$}^{\varepsilon}_{0} \) satisfies the 
\( k\)-th order compatibility condition for all \( k\) with 
\( 0\leq k \leq m \).
\end{df}
\subsection{Corrections to the initial datum}
Now we construct a corrected initial datum \( \mbox{\mathversion{bold}$v$}^{\varepsilon}_{0}\)
such that given an initial datum \( \mbox{\mathversion{bold}$v$}_{0} \) that satisfies 
the compatibility conditions for (\ref{vslant}), 
\( \mbox{\mathversion{bold}$v$}^{\varepsilon}_{0}\) satisfies the compatibility conditions
for (\ref{rnl}) and \( \mbox{\mathversion{bold}$v$}^{\varepsilon}_{0} \to 
\mbox{\mathversion{bold}$v$}_{0}\) in the appropriate function space.
We must also correct the initial datum so that \( |\mbox{\mathversion{bold}$v$}^{\varepsilon}_{0}|
\equiv 1 \) is satisfied.

Suppose that we have an initial datum 
\( \mbox{\mathversion{bold}$v$}_{0}\in H^{2m+1}(I)\) satisfying the compatibility conditions 
for (\ref{vslant}) up to order \( m\) (the case when 
the initial datum belongs to a Sobolev space with even index will be remarked on at 
the end). We construct
\( \mbox{\mathversion{bold}$v$}^{\varepsilon}_{0}\) in the form 
\begin{eqnarray}
\mbox{\mathversion{bold}$v$}^{\varepsilon}_{0}=
\frac{\mbox{\mathversion{bold}$v$}_{0} + 
\mbox{\mathversion{bold}$h$}^{\varepsilon}}{|
\mbox{\mathversion{bold}$v$}_{0} + 
\mbox{\mathversion{bold}$h$}^{\varepsilon }|},
\label{cor}
\end{eqnarray}
where \( \mbox{\mathversion{bold}$h$}^{\varepsilon}\) is constructed so that
\( \mbox{\mathversion{bold}$h$}^{\varepsilon} \to \mbox{\mathversion{bold}$0$}\)
as \( \varepsilon \to +0 \) in \( H^{2m+1}(I) \). We do this by 
determining the differential coeffecients of 
\( \mbox{\mathversion{bold}$h$}^{\varepsilon}\) at \( s=0,1\) and extend it to \( I\) so that 
\( \mbox{\mathversion{bold}$h$}^{\varepsilon}\) belongs to \( H^{2m+1}(I) \).
We introduce some notations. We set
\begin{align*}
\mbox{\mathversion{bold}$g$}^{\varepsilon}_{0}
( \mbox{\mathversion{bold}$V$})&:= \mbox{\mathversion{bold}$V$},\\[3mm]
\mbox{\mathversion{bold}$g$}^{\varepsilon}_{1}
( \mbox{\mathversion{bold}$V$})&:= 
\mbox{\mathversion{bold}$V$}\times \mbox{\mathversion{bold}$V$}_{ss}
+ \varepsilon \mbox{\mathversion{bold}$V$}_{ss} 
+ \varepsilon | \mbox{\mathversion{bold}$V$}_{s}|^{2}\mbox{\mathversion{bold}$V$}, \\[3mm]
\mbox{\mathversion{bold}$g$}^{\varepsilon}_{m+1}
( \mbox{\mathversion{bold}$V$})&:= D\mbox{\mathversion{bold}$g$}^{\varepsilon}_{m}
(\mbox{\mathversion{bold}$V$})
\big[ \mbox{\mathversion{bold}$g$}^{\varepsilon}_{1}
( \mbox{\mathversion{bold}$V$})\big],
\end{align*}
where \( m \geq 1\) and \( D\) is the derivative with respect to 
\( \mbox{\mathversion{bold}$V$}\), i.e. 
\( D\mbox{\mathversion{bold}$g$}^{\varepsilon}_{m}
(\mbox{\mathversion{bold}$V$})
\big[ \mbox{\mathversion{bold}$W$} \big] = \frac{{\rm d}}{{\rm d}r}
\mbox{\mathversion{bold}$g$}^{\varepsilon}_{m}(
\mbox{\mathversion{bold}$V$}+r\mbox{\mathversion{bold}$W$})|_{r=0} \).
Here, \( |_{r=0} \) is the trace at \( r=0 \).
Under these notations, the \( m\)-th order compatibility condition for 
(\ref{rnl}) can be expressed as 
\( \mbox{\mathversion{bold}$g$}^{\varepsilon}_{m}(
\mbox{\mathversion{bold}$v$}^{\varepsilon}_{0}(0))=
\mbox{\mathversion{bold}$g$}^{\varepsilon}_{m}(
\mbox{\mathversion{bold}$v$}^{\varepsilon}_{0}(1))=
\mbox{\mathversion{bold}$0$}\) because
\( \mbox{\mathversion{bold}$g$}^{\varepsilon}_{m}(\mbox{\mathversion{bold}$V$})
=
\mbox{\mathversion{bold}$Q$}_{m}(\mbox{\mathversion{bold}$V$}) \).
We gave a different notation because it is more convenient for the upcoming calculations.
We first prove 
\begin{lm}
If \( |\mbox{\mathversion{bold}$V$}|\equiv 1 \), then for 
any \( m\geq 1 \),
\begin{eqnarray}
\sum^{m}_{k=0}
\left(
\begin{array}{c}
m\\
k
\end{array}\right)
\mbox{\mathversion{bold}$g$}^{\varepsilon}_{k}(\mbox{\mathversion{bold}$V$})
\cdot
\mbox{\mathversion{bold}$g$}^{\varepsilon}_{m-k}(\mbox{\mathversion{bold}$V$})
\equiv 0.
\label{zero}
\end{eqnarray}
\label{lmzero}
\end{lm}
{\it Proof.}
We show this by induction.
By direct calculation, we see that
\begin{eqnarray*}
\mbox{\mathversion{bold}$V$}\cdot \mbox{\mathversion{bold}$g$}^{\varepsilon}_{1}
( \mbox{\mathversion{bold}$V$})=
\frac{\varepsilon }{2} \left(|\mbox{\mathversion{bold}$V$}|^{2}\right)_{ss}
+ \varepsilon |\mbox{\mathversion{bold}$V$}_{s}|^{2}
\left( |\mbox{\mathversion{bold}$V$}|^{2}-1\right)
\equiv 0,
\end{eqnarray*}
which proves (\ref{zero}) with \( m=1 \).
Suppose (\ref{zero}) holds up to some \( m\) with \( m\geq 1 \). From the assumption of 
induction, we have for any vector \( \mbox{\mathversion{bold}$W$}\) and \( r\in \mathbf{R}\),
\begin{eqnarray*}
\sum^{m}_{k=0}
\left(
\begin{array}{c}
m\\
k
\end{array}\right)
\mbox{\mathversion{bold}$g$}^{\varepsilon}_{k}\left(
\frac{\mbox{\mathversion{bold}$V$}+r\mbox{\mathversion{bold}$W$}}
{|\mbox{\mathversion{bold}$V$}+r\mbox{\mathversion{bold}$W$}|}\right)
\cdot
\mbox{\mathversion{bold}$g$}^{\varepsilon}_{m-k}
\left(\frac{\mbox{\mathversion{bold}$V$}+r\mbox{\mathversion{bold}$W$}}
{|\mbox{\mathversion{bold}$V$}+r\mbox{\mathversion{bold}$W$}|}\right)
\equiv 0.
\end{eqnarray*}
Differentiating with respect to \( r\) and setting \( r=0\) yields
\begin{align*}
\sum^{m}_{k=0}
\left(
\begin{array}{c}
m\\
k
\end{array}\right)
&
\{
D\mbox{\mathversion{bold}$g$}^{\varepsilon}_{k}
(\mbox{\mathversion{bold}$V$})[\mbox{\mathversion{bold}$W$}
-(\mbox{\mathversion{bold}$V$}\cdot \mbox{\mathversion{bold}$W$})
\mbox{\mathversion{bold}$V$}]
\cdot
\mbox{\mathversion{bold}$g$}^{\varepsilon}_{m-k}
(\mbox{\mathversion{bold}$V$}) \\[3mm]
&+
\mbox{\mathversion{bold}$g$}^{\varepsilon}_{k}
(\mbox{\mathversion{bold}$V$})
\cdot
D\mbox{\mathversion{bold}$g$}^{\varepsilon}_{m-k}
(\mbox{\mathversion{bold}$V$})[\mbox{\mathversion{bold}$W$}
-(\mbox{\mathversion{bold}$V$}\cdot \mbox{\mathversion{bold}$W$})
\mbox{\mathversion{bold}$V$}]\}
\equiv 0.
\end{align*}
By choosing \( \mbox{\mathversion{bold}$W$}=
\mbox{\mathversion{bold}$g$}^{\varepsilon}_{1}(\mbox{\mathversion{bold}$V$})\), we have 
\begin{align*}
0 &\equiv 
\sum^{m}_{k=0}
\left(
\begin{array}{c}
m\\
k
\end{array}\right)
\{
\mbox{\mathversion{bold}$g$}^{\varepsilon}_{k+1}
(\mbox{\mathversion{bold}$V$})
\cdot
\mbox{\mathversion{bold}$g$}^{\varepsilon}_{m-k}
(\mbox{\mathversion{bold}$V$})
+
\mbox{\mathversion{bold}$g$}^{\varepsilon}_{k}
(\mbox{\mathversion{bold}$V$})
\cdot
\mbox{\mathversion{bold}$g$}^{\varepsilon}_{m+1-k}
(\mbox{\mathversion{bold}$V$}) \\[3mm]
&=
\sum^{m+1}_{k=0}
\left(
\begin{array}{c}
m+1\\
k
\end{array}\right)
\mbox{\mathversion{bold}$g$}^{\varepsilon}_{k}
(\mbox{\mathversion{bold}$V$})
\cdot
\mbox{\mathversion{bold}$g$}^{\varepsilon}_{m+1-k}
(\mbox{\mathversion{bold}$V$}),
\end{align*}
which proves (\ref{zero}) for the case \( m+1 \), and this finishes the proof. \hfill \( \Box. \)

Next, we make the following notations.
\begin{align*}
\mbox{\mathversion{bold}$f$}_{0}(\mbox{\mathversion{bold}$V$})
&:=\mbox{\mathversion{bold}$V$},\\[3mm]
\mbox{\mathversion{bold}$f$}_{1}(\mbox{\mathversion{bold}$V$})
&:=\mbox{\mathversion{bold}$V$}\times \mbox{\mathversion{bold}$V$}_{ss}, \\[3mm]
\mbox{\mathversion{bold}$f$}_{m+1}(\mbox{\mathversion{bold}$V$})
&:= D\mbox{\mathversion{bold}$f$}_{m}(\mbox{\mathversion{bold}$V$})
[\mbox{\mathversion{bold}$f$}_{1}(\mbox{\mathversion{bold}$V$})],
\end{align*}
which is equivalent to taking \( \varepsilon =0 \) in 
\( \mbox{\mathversion{bold}$g$}^{\varepsilon}_{m}\). Hence, 
\begin{eqnarray*}
\sum^{m}_{k=0}
\left(
\begin{array}{c}
m\\
k
\end{array}\right)
\mbox{\mathversion{bold}$f$}_{k}(\mbox{\mathversion{bold}$V$})
\cdot
\mbox{\mathversion{bold}$f$}_{m-k}(\mbox{\mathversion{bold}$V$})
\equiv 0
\end{eqnarray*}
for any vector \( \mbox{\mathversion{bold}$V$}\) satisfying
\( | \mbox{\mathversion{bold}$V$}|\equiv 1 \). Also, the 
\( m\)-th order compatibility condition for (\ref{vslant}) can be expressed as 
\( \mbox{\mathversion{bold}$f$}_{m}(\mbox{\mathversion{bold}$v$}_{0}(0))
=\mbox{\mathversion{bold}$f$}_{m}(\mbox{\mathversion{bold}$v$}_{0}(1))
=\mbox{\mathversion{bold}$0$}\) for \( m \geq 1 \).

\medskip
From here we look into the structure of 
\( \mbox{\mathversion{bold}$f$}_{m} \) and \( \mbox{\mathversion{bold}$g$}^{\varepsilon}_{m}\)
in more detail. 
This will allow us to determine the differential coefficients of 
the correction term \( \mbox{\mathversion{bold}$h$}^{\varepsilon} \).
We prove 
\begin{lm}
For \( m \geq 1 \),
\begin{eqnarray}
\mbox{\mathversion{bold}$g$}^{\varepsilon}_{m}(\mbox{\mathversion{bold}$V$})
=
\mbox{\mathversion{bold}$f$}_{m}(\mbox{\mathversion{bold}$V$})+
\varepsilon \mbox{\mathversion{bold}$r$}^{\varepsilon}_{m}(\mbox{\mathversion{bold}$V$}),
\label{differ}
\end{eqnarray}
where
\( \mbox{\mathversion{bold}$r$}^{\varepsilon}_{1}(\mbox{\mathversion{bold}$V$}):=
\mbox{\mathversion{bold}$V$}_{ss} + |\mbox{\mathversion{bold}$V$}_{s}|^{2}
\mbox{\mathversion{bold}$V$} \) and 
\begin{eqnarray*}
\mbox{\mathversion{bold}$r$}^{\varepsilon}_{m}(\mbox{\mathversion{bold}$V$}):=
D\mbox{\mathversion{bold}$r$}^{\varepsilon}_{m-1}(\mbox{\mathversion{bold}$V$})
[\mbox{\mathversion{bold}$g$}^{\varepsilon}_{1}(\mbox{\mathversion{bold}$V$})]
+
D\mbox{\mathversion{bold}$f$}_{m-1}(\mbox{\mathversion{bold}$V$})
[\mbox{\mathversion{bold}$r$}^{\varepsilon}_{1}(\mbox{\mathversion{bold}$V$})]
\end{eqnarray*}
for \( m\geq 2\). Formula {\rm (\ref{differ})} tells us that the difference 
\( \mbox{\mathversion{bold}$g$}^{\varepsilon}_{m}-\mbox{\mathversion{bold}$f$}_{m}\) is of order \( \varepsilon \).
\end{lm}
{\it Proof.}
We prove (\ref{differ}) by induction. It is obvious that 
(\ref{differ}) holds for \( m=1 \) from the definition of 
\( \mbox{\mathversion{bold}$g$}^{\varepsilon}_{1}\), \( \mbox{\mathversion{bold}$f$}_{1}\),
and \( \mbox{\mathversion{bold}$r$}^{\varepsilon}_{1}\).
Suppose it holds up to \( m-1 \) for some \( m \geq 2\). From the assumption of 
induction, we have for any vector \( \mbox{\mathversion{bold}$W$}\) and 
\( r\in \mathbf{R}\),
\begin{eqnarray*}
\mbox{\mathversion{bold}$g$}^{\varepsilon}_{m-1}(\mbox{\mathversion{bold}$V$}
+r\mbox{\mathversion{bold}$W$})
=
\mbox{\mathversion{bold}$f$}_{m-1}(\mbox{\mathversion{bold}$V$}
+r\mbox{\mathversion{bold}$W$})+
\varepsilon \mbox{\mathversion{bold}$r$}^{\varepsilon}_{m-1}(\mbox{\mathversion{bold}$V$}
+r\mbox{\mathversion{bold}$W$})
\end{eqnarray*}
holds. Differentiating with respect to \( r\) and setting \( r=0 \) yields
\begin{eqnarray*}
D\mbox{\mathversion{bold}$g$}^{\varepsilon}_{m-1}(\mbox{\mathversion{bold}$V$})
[\mbox{\mathversion{bold}$W$}]
=
D\mbox{\mathversion{bold}$f$}_{m-1}(\mbox{\mathversion{bold}$V$})
[\mbox{\mathversion{bold}$W$}]+
\varepsilon D\mbox{\mathversion{bold}$r$}^{\varepsilon}_{m-1}
(\mbox{\mathversion{bold}$V$})[\mbox{\mathversion{bold}$W$}].
\end{eqnarray*}
Choosing \( \mbox{\mathversion{bold}$W$}=
\mbox{\mathversion{bold}$g$}^{\varepsilon}_{1}(\mbox{\mathversion{bold}$V$}) \) yields
\begin{align*}
\mbox{\mathversion{bold}$g$}^{\varepsilon}_{m}(\mbox{\mathversion{bold}$V$})
&=
D\mbox{\mathversion{bold}$f$}_{m-1}(\mbox{\mathversion{bold}$V$})
[\mbox{\mathversion{bold}$g$}^{\varepsilon}_{1}(\mbox{\mathversion{bold}$V$})]+
\varepsilon D\mbox{\mathversion{bold}$r$}^{\varepsilon}_{m-1}
(\mbox{\mathversion{bold}$V$})[\mbox{\mathversion{bold}$g$}^{\varepsilon}_{1}(\mbox{\mathversion{bold}$V$}]\\[3mm]
&=D\mbox{\mathversion{bold}$f$}_{m-1}(\mbox{\mathversion{bold}$V$})
[\mbox{\mathversion{bold}$f$}_{1}(\mbox{\mathversion{bold}$V$})
+ \varepsilon \mbox{\mathversion{bold}$r$}^{\varepsilon }_{1}
(\mbox{\mathversion{bold}$V$})]+
\varepsilon D\mbox{\mathversion{bold}$r$}^{\varepsilon}_{m-1}
(\mbox{\mathversion{bold}$V$})[\mbox{\mathversion{bold}$g$}^{\varepsilon}_{1}(\mbox{\mathversion{bold}$V$})]\\[3mm]
&=\mbox{\mathversion{bold}$f$}_{m}(\mbox{\mathversion{bold}$V$})
+
\varepsilon \big\{
D\mbox{\mathversion{bold}$f$}_{m-1}(\mbox{\mathversion{bold}$V$})
[\mbox{\mathversion{bold}$r$}^{\varepsilon}_{1}
(\mbox{\mathversion{bold}$V$})]
+
D\mbox{\mathversion{bold}$r$}^{\varepsilon}_{m-1}(\mbox{\mathversion{bold}$V$})
[\mbox{\mathversion{bold}$g$}^{\varepsilon}_{1}(\mbox{\mathversion{bold}$V$})]
\big\}\\[3mm]
&=
\mbox{\mathversion{bold}$f$}_{m}(\mbox{\mathversion{bold}$V$})+
\varepsilon \mbox{\mathversion{bold}$r$}^{\varepsilon}_{m}(\mbox{\mathversion{bold}$V$}),
\end{align*}
which finishes the proof. \hfill \( \Box\).
Next we prove 
\begin{lm}
For \( m\geq 1\), 
\begin{eqnarray}
\mbox{\mathversion{bold}$f$}_{m}(\mbox{\mathversion{bold}$V$})
= A_{m}(\mbox{\mathversion{bold}$V$})\partial ^{2m}_{s}\mbox{\mathversion{bold}$V$}
+
\mbox{\mathversion{bold}$\phi $}_{m}(\mbox{\mathversion{bold}$V$})
\label{fhigh}
\end{eqnarray}
holds, where the vector \( \mbox{\mathversion{bold}$\phi$}_{m}(\mbox{\mathversion{bold}$V$})\)
 and the operator \( A_{m}(\mbox{\mathversion{bold}$V$})\) are defined as follows.
\begin{align*}
A_{0}(\mbox{\mathversion{bold}$V$})\mbox{\mathversion{bold}$W$}=
\mbox{\mathversion{bold}$W$}, \ 
A_{1}(\mbox{\mathversion{bold}$V$})\mbox{\mathversion{bold}$W$}&=
\mbox{\mathversion{bold}$V$}\times \mbox{\mathversion{bold}$W$}\\[3mm]
\mbox{\mathversion{bold}$\phi $}_{1}(\mbox{\mathversion{bold}$V$}) &= \mbox{\mathversion{bold}$0$}.
\end{align*}
and for \( m\geq 2 \),
\begin{align*}
A_{m}(\mbox{\mathversion{bold}$V$})\mbox{\mathversion{bold}$W$} &= \mbox{\mathversion{bold}$V$}\times 
\big( A_{m-1}(\mbox{\mathversion{bold}$V$})\mbox{\mathversion{bold}$W$}\big)\\[3mm]
\mbox{\mathversion{bold}$\phi $}_{m}(\mbox{\mathversion{bold}$V$})&=
D\mbox{\mathversion{bold}$\phi $} _{m-1}(\mbox{\mathversion{bold}$V$})[\mbox{\mathversion{bold}$f$}_{1}
(\mbox{\mathversion{bold}$V$})]+
A_{m-1}(\mbox{\mathversion{bold}$V$})
\sum ^{2m-3}_{j=0}
\left(
\begin{array}{c}
2(m-1)\\
j
\end{array}
\right)
\partial^{2(m-1)-j}_{s}\mbox{\mathversion{bold}$V$}\times
\partial^{j+2}_{s}\mbox{\mathversion{bold}$V$}\\[3mm]
& \quad \quad +
 \left( DA_{m}(\mbox{\mathversion{bold}$V$})[\mbox{\mathversion{bold}$f$}_{1}
(\mbox{\mathversion{bold}$V$})]\right) \partial^{2(m-1)}_{s}\mbox{\mathversion{bold}$V$}.
\end{align*}
Although we do not need \( A_{0} \) for this lemma, we defined it because 
we will use it later.
We also note that from the definition of \( \mbox{\mathversion{bold}$\phi $} _{m}\), 
\( \mbox{\mathversion{bold}$\phi $}_{m}(\mbox{\mathversion{bold}$V$}) \) satisfies
\begin{eqnarray*}
|\mbox{\mathversion{bold}$\phi $}_{m}(\mbox{\mathversion{bold}$V$})|\leq 
C\big( |\mbox{\mathversion{bold}$V$}|+|\mbox{\mathversion{bold}$V$}_{s}| + 
\cdots +|\partial ^{2m-1}_{s}\mbox{\mathversion{bold}$V$}|\big)
\end{eqnarray*}
if
\( |\mbox{\mathversion{bold}$V$}|+|\mbox{\mathversion{bold}$V$}_{s}| + 
\cdots +|\partial ^{2m-1}_{s}\mbox{\mathversion{bold}$V$}| \leq M \) and \( C>0\) depends on 
\( M\).
Hence, formula {\rm (\ref{fhigh})} gives the explicit form of the term with the highest order of 
\( s\)-derivatives in \( \mbox{\mathversion{bold}$f$}_{m}(\mbox{\mathversion{bold}$V$}) \).
\label{fhighlm}
\end{lm}
{\it Proof.}
We see that (\ref{fhigh}) holds for \( m=1 \) from the definition of 
\( \mbox{\mathversion{bold}$f$}_{1}(\mbox{\mathversion{bold}$V$})\). Suppose it holds up to 
\( m-1 \) for some \( m \geq 2 \). Then for any vector \( \mbox{\mathversion{bold}$W$}\) and
\( r\in \mathbf{R}\) we have
\begin{eqnarray*}
\mbox{\mathversion{bold}$f$}_{m-1}(\mbox{\mathversion{bold}$V$}+r\mbox{\mathversion{bold}$W$})
=A_{m-1}(\mbox{\mathversion{bold}$V$}+r\mbox{\mathversion{bold}$W$})\partial^{2(m-1)}_{s}
(\mbox{\mathversion{bold}$V$}+r\mbox{\mathversion{bold}$W$}) + 
\mbox{\mathversion{bold}$\phi $}_{m-1}(\mbox{\mathversion{bold}$V$}+r\mbox{\mathversion{bold}$W$}).
\end{eqnarray*}
Differentiating with respect to \( r\), setting \( r=0\), and choosing 
\( \mbox{\mathversion{bold}$W$}= \mbox{\mathversion{bold}$f$}_{1}(\mbox{\mathversion{bold}$V$}) \) yields
\begin{align*}
\mbox{\mathversion{bold}$f$}_{m}(\mbox{\mathversion{bold}$V$}) 
&= 
A_{m-1}(\mbox{\mathversion{bold}$V$})\partial^{2(m-1)}_{s}(\mbox{\mathversion{bold}$f$}_{1}(\mbox{\mathversion{bold}$V$}))\\[3mm]
& \qquad +
\big( DA_{m-1}(\mbox{\mathversion{bold}$V$})[\mbox{\mathversion{bold}$f$}_{1}(\mbox{\mathversion{bold}$V$})]\big)
\partial^{2(m-1)}_{s}\mbox{\mathversion{bold}$V$}
+
D\mbox{\mathversion{bold}$\phi $}_{m-1}(\mbox{\mathversion{bold}$V$})
[\mbox{\mathversion{bold}$f$}_{1}(\mbox{\mathversion{bold}$V$})]
\end{align*}
Since \( \mbox{\mathversion{bold}$f$}_{1}(\mbox{\mathversion{bold}$V$})=
\mbox{\mathversion{bold}$V$}\times \mbox{\mathversion{bold}$V$}_{ss} \), we have
\begin{align*}
A_{m-1}(\mbox{\mathversion{bold}$V$})\partial^{2(m-1)}_{s}\big( \mbox{\mathversion{bold}$f$}_{1}(\mbox{\mathversion{bold}$V$})\big)
&=
A_{m-1}(\mbox{\mathversion{bold}$V$})(\mbox{\mathversion{bold}$V$}\times \partial^{2m}_{s}\mbox{\mathversion{bold}$V$})\\[3mm]
& \qquad +
A_{m-1}(\mbox{\mathversion{bold}$V$})\left\{
\sum^{2m-3}_{j=0}
\left(
\begin{array}{c}
2(m-1)\\
j
\end{array}
\right)
\partial^{2(m-1)-j}_{s}\mbox{\mathversion{bold}$V$}\times
\partial^{j+2}_{s}\mbox{\mathversion{bold}$V$}
\right\}\\[3mm]
&=A_{m}(\mbox{\mathversion{bold}$V$})\partial^{2m}_{s}\mbox{\mathversion{bold}$V$}\\[3mm]
& \qquad +
A_{m-1}(\mbox{\mathversion{bold}$V$})\left\{
\sum^{2m-3}_{j=0}
\left(
\begin{array}{c}
2(m-1)\\
j
\end{array}
\right)
\partial^{2(m-1)-j}_{s}\mbox{\mathversion{bold}$V$}\times
\partial^{j+2}_{s}\mbox{\mathversion{bold}$V$}
\right\}.
\end{align*}
Hence we have
\begin{align*}
\mbox{\mathversion{bold}$f$}_{m}(\mbox{\mathversion{bold}$V$})
&=
A_{m}(\mbox{\mathversion{bold}$V$})\partial^{2m}_{s}\mbox{\mathversion{bold}$V$}\\[3mm]
&\qquad +
A_{m-1}(\mbox{\mathversion{bold}$V$})
\sum^{2m-3}_{j=0}
\left(
\begin{array}{c}
2(m-1)\\
j
\end{array}
\right)
\partial^{2(m-1)-j}_{s}\mbox{\mathversion{bold}$V$}\times
\partial^{j+2}_{s}\mbox{\mathversion{bold}$V$}\\[3mm]
& \qquad +
\big( DA_{m-1}(\mbox{\mathversion{bold}$V$})[\mbox{\mathversion{bold}$f$}_{1}(\mbox{\mathversion{bold}$V$})]\big)
\partial^{2(m-1)}_{s}\mbox{\mathversion{bold}$V$}
+
D\mbox{\mathversion{bold}$\phi $}_{m-1}(\mbox{\mathversion{bold}$V$})
[\mbox{\mathversion{bold}$f$}_{1}(\mbox{\mathversion{bold}$V$})]\\[3mm]
&=A_{m}(\mbox{\mathversion{bold}$V$})\partial^{2m}_{s}\mbox{\mathversion{bold}$V$}
+
\mbox{\mathversion{bold}$\phi $}_{m}(\mbox{\mathversion{bold}$V$}),
\end{align*}
which finishes the proof. \hfill \( \Box \).
Next we prove 
\begin{lm}
For \( m \geq 1 \),
\begin{eqnarray}
\mbox{\mathversion{bold}$r$}^{\varepsilon}_{m}(\mbox{\mathversion{bold}$V$})
=
\sum^{m-1}_{j=0}c_{m,j}\varepsilon ^{j}A_{m-1-j}(\mbox{\mathversion{bold}$V$})\partial^{2m}_{s}\mbox{\mathversion{bold}$V$}
+
\mbox{\mathversion{bold}$\varphi $}^{\varepsilon}_{m}(\mbox{\mathversion{bold}$V$})
\label{remainder}
\end{eqnarray}
holds. Here, \( \mbox{\mathversion{bold}$\varphi $}^{\varepsilon}_{1}(\mbox{\mathversion{bold}$V$})
= |\mbox{\mathversion{bold}$V$}_{s}|^{2}\mbox{\mathversion{bold}$V$} \) and for \( m\geq 2\),
\begin{align*}
\mbox{\mathversion{bold}$\varphi $}^{\varepsilon}_{m}(\mbox{\mathversion{bold}$V$})
&=
D\mbox{\mathversion{bold}$\varphi $}^{\varepsilon}_{m-1}(\mbox{\mathversion{bold}$V$})
[\mbox{\mathversion{bold}$g$}^{\varepsilon}_{1}(\mbox{\mathversion{bold}$V$})]
+
D\mbox{\mathversion{bold}$\phi $}_{m-1}(\mbox{\mathversion{bold}$V$})
[\mbox{\mathversion{bold}$r$}^{\varepsilon}_{1}(\mbox{\mathversion{bold}$V$})]\\[3mm]
& \qquad +
\sum ^{m-1}_{j=0}d_{m-1,j}\varepsilon^{j}A_{m-1-j}(\mbox{\mathversion{bold}$V$})
\left\{ \partial^{2(m-1)}_{s}\big( |\mbox{\mathversion{bold}$V$}_{s}|^{2}\mbox{\mathversion{bold}$V$}\big)\right\}\\[3mm]
& \qquad +
\sum ^{m-2}_{j=0}c_{m-1,j}\varepsilon^{j}A_{m-2-j}(\mbox{\mathversion{bold}$V$})
\left\{
\sum^{2m-3}_{k=0}
\left(
\begin{array}{c}
2m-2\\
k
\end{array}\right)
\partial^{2m-2-k}_{s}\mbox{\mathversion{bold}$V$}\times \partial^{j+2}\mbox{\mathversion{bold}$V$}
\right\}\\[3mm]
& \qquad +
\sum^{m-1}_{j=0}d_{m-1,j}\varepsilon^{j}\big(
DA_{m-1-j}(\mbox{\mathversion{bold}$V$})[\mbox{\mathversion{bold}$r$}^{\varepsilon}_{1}(\mbox{\mathversion{bold}$V$})]\big)
\partial^{2m-2}_{s}\mbox{\mathversion{bold}$V$}\\[3mm]
& \qquad +
\sum^{m-2}_{j=0}c_{m-1,j}\varepsilon^{j}\big( 
DA_{m-2-j}(\mbox{\mathversion{bold}$V$})[\mbox{\mathversion{bold}$f$}_{1}(\mbox{\mathversion{bold}$V$})]\big)
\partial^{2m-2}_{s}\mbox{\mathversion{bold}$V$}.
\end{align*}
\( c_{m,j}\) and \( d_{m,j} \) are constants defined as follows. \( c_{1,0}=1 \) and
\begin{align*}
c_{m,j}=
\left\{
\begin{array}{ll}
m, & j=0, \\[3mm]
c_{m-1,j-1}+c_{m-1,j}, & 1\leq j \leq m-2, \\[3mm]
c_{m-1,m-2}, & j=m-1,
\end{array}\right.
\end{align*}
for \( m\geq 2 \). \( d_{m,0}=1 \) and \( d_{m,j}=c_{m,j-1} \) for \( 1\leq j \leq m \).
We also note that from the definition of \( \mbox{\mathversion{bold}$\varphi $}^{\varepsilon} _{m}\), 
\( \mbox{\mathversion{bold}$\varphi $}^{\varepsilon}_{m}(\mbox{\mathversion{bold}$V$}) \) satisfies
\begin{eqnarray*}
|\mbox{\mathversion{bold}$\varphi $}^{\varepsilon}_{m}(\mbox{\mathversion{bold}$V$})|\leq 
C\big( |\mbox{\mathversion{bold}$V$}|+|\mbox{\mathversion{bold}$V$}_{s}| + 
\cdots +|\partial ^{2m-1}_{s}\mbox{\mathversion{bold}$V$}|\big)
\end{eqnarray*}
if
\( |\mbox{\mathversion{bold}$V$}|+|\mbox{\mathversion{bold}$V$}_{s}| + 
\cdots +|\partial ^{2m-1}_{s}\mbox{\mathversion{bold}$V$}| \leq M \) and \( C>0\) depends on 
\( M\).
\label{rmhigh}
\end{lm}
{\it Proof.}
Since \( \mbox{\mathversion{bold}$r$}^{\varepsilon}_{1}(\mbox{\mathversion{bold}$V$}) = 
\mbox{\mathversion{bold}$V$}_{ss}+|\mbox{\mathversion{bold}$V$}_{s}|^{2}\mbox{\mathversion{bold}$V$}\), 
(\ref{remainder}) holds with \( m=1 \). Suppose it holds up to \( m-1\) for some \( m \geq 2\).
Then we have
\begin{eqnarray*}
\mbox{\mathversion{bold}$r$}^{\varepsilon}_{m-1}(\mbox{\mathversion{bold}$V$}+r\mbox{\mathversion{bold}$W$})
=
\sum ^{m-2}_{j=0}c_{m-1,j}\varepsilon^{j}A_{m-2-j}(\mbox{\mathversion{bold}$V$}+r\mbox{\mathversion{bold}$W$})
\partial^{2(m-1)}_{s}(\mbox{\mathversion{bold}$V$}+r\mbox{\mathversion{bold}$W$})
+
\mbox{\mathversion{bold}$\varphi $}^{\varepsilon}_{m-1}(\mbox{\mathversion{bold}$V$}+r\mbox{\mathversion{bold}$W$})
\end{eqnarray*}
for any vector \( \mbox{\mathversion{bold}$W$}\) and \( r\in \mathbf{R}\). Differentiating with respect to \( r\), 
setting \( r=0 \), and choosing \( \mbox{\mathversion{bold}$W$}=
\mbox{\mathversion{bold}$g$}^{\varepsilon}_{1}(\mbox{\mathversion{bold}$V$})\) yields
\begin{align*}
D\mbox{\mathversion{bold}$r$}^{\varepsilon}_{m-1}(\mbox{\mathversion{bold}$V$})
[\mbox{\mathversion{bold}$g$}^{\varepsilon}_{1}(\mbox{\mathversion{bold}$V$})]
&=
\sum^{m-2}_{j=0}c_{m-1,j}\varepsilon^{j}A_{m-2-j}(\mbox{\mathversion{bold}$V$})
\big( \partial^{2(m-1)}_{s}\mbox{\mathversion{bold}$g$}^{\varepsilon}_{1}(\mbox{\mathversion{bold}$V$})\big)\\[3mm]
& \qquad +
\sum^{m-2}_{j=0}c_{m-1,j}\varepsilon^{j}
\big( DA_{m-2-j}(\mbox{\mathversion{bold}$V$})[\mbox{\mathversion{bold}$g$}^{\varepsilon}_{1}(\mbox{\mathversion{bold}$V$})]\big)
\partial^{2(m-1)}_{s}\mbox{\mathversion{bold}$V$}\\[3mm]
& \qquad +
D\mbox{\mathversion{bold}$\varphi $}^{\varepsilon}_{m-1}(\mbox{\mathversion{bold}$V$})
[\mbox{\mathversion{bold}$g$}^{\varepsilon}_{1}(\mbox{\mathversion{bold}$V$})].
\end{align*}
Furthermore, we have proved that 
\( \mbox{\mathversion{bold}$f$}_{m-1}(\mbox{\mathversion{bold}$V$}) = 
A_{m-1}(\mbox{\mathversion{bold}$V$})\partial^{2(m-1)}_{s}\mbox{\mathversion{bold}$V$}+
\mbox{\mathversion{bold}$\phi $}_{m-1}(\mbox{\mathversion{bold}$V$}) \), so that
\begin{align*}
D\mbox{\mathversion{bold}$f$}_{m-1}(\mbox{\mathversion{bold}$V$})
[\mbox{\mathversion{bold}$r$}^{\varepsilon}_{1}(\mbox{\mathversion{bold}$V$})]
&=
A_{m-1}(\mbox{\mathversion{bold}$V$})
\big( \partial^{2(m-1)}_{s}\mbox{\mathversion{bold}$r$}^{\varepsilon}_{1}(\mbox{\mathversion{bold}$V$})\big)
+
\big(DA_{m-1}(\mbox{\mathversion{bold}$V$})[\mbox{\mathversion{bold}$r$}^{\varepsilon}_{1}(\mbox{\mathversion{bold}$V$})]\big)
\partial^{2(m-1)}_{s}\mbox{\mathversion{bold}$V$}\\[3mm]
& \qquad +
D\mbox{\mathversion{bold}$\phi $}_{m-1}(\mbox{\mathversion{bold}$V$})
[\mbox{\mathversion{bold}$r$}^{\varepsilon}_{1}(\mbox{\mathversion{bold}$V$})]
\end{align*}
holds.
Substituting the definition of \( \mbox{\mathversion{bold}$r$}^{\varepsilon}_{1}(\mbox{\mathversion{bold}$V$})\)
into the first term on the right-hand side yields
\begin{align*}
D\mbox{\mathversion{bold}$f$}_{m-1}(\mbox{\mathversion{bold}$V$})
[\mbox{\mathversion{bold}$r$}^{\varepsilon}_{1}(\mbox{\mathversion{bold}$V$})]
&=
A_{m-1}(\mbox{\mathversion{bold}$V$})\partial^{2m}_{s}\mbox{\mathversion{bold}$V$}
+
A_{m-1}(\mbox{\mathversion{bold}$V$})\big( \partial^{2(m-1)}_{s}
(|\mbox{\mathversion{bold}$V$}_{s}|^{2}\mbox{\mathversion{bold}$V$})\big)\\[3mm]
& \qquad +
\big( DA_{m-1}(\mbox{\mathversion{bold}$V$})
[\mbox{\mathversion{bold}$r$}^{\varepsilon}_{1}(\mbox{\mathversion{bold}$V$})]\big)
\partial^{2(m-1)}_{s}\mbox{\mathversion{bold}$V$}
+
D\mbox{\mathversion{bold}$\phi $}_{m-1}(\mbox{\mathversion{bold}$V$})
[\mbox{\mathversion{bold}$r$}^{\varepsilon}_{1}(\mbox{\mathversion{bold}$V$})].
\end{align*}
From the definition of \( \mbox{\mathversion{bold}$r$}^{\varepsilon}_{m}(\mbox{\mathversion{bold}$V$})\), we have
\begin{align*}
\mbox{\mathversion{bold}$r$}^{\varepsilon}_{m}(\mbox{\mathversion{bold}$V$})
&=
D\mbox{\mathversion{bold}$r$}^{\varepsilon}_{m-1}(\mbox{\mathversion{bold}$V$})
[\mbox{\mathversion{bold}$g$}^{\varepsilon}_{1}(\mbox{\mathversion{bold}$V$})]
+
D\mbox{\mathversion{bold}$f$}_{m-1}(\mbox{\mathversion{bold}$V$})
[\mbox{\mathversion{bold}$r$}^{\varepsilon}_{1}(\mbox{\mathversion{bold}$V$})]\\[3mm]
&=\sum^{m-2}_{j=0}c_{m-1,j}\varepsilon^{j}A_{m-2-j}(\mbox{\mathversion{bold}$V$})
\big( \partial^{2(m-1)}_{s}\mbox{\mathversion{bold}$g$}^{\varepsilon}_{1}(\mbox{\mathversion{bold}$V$})\big)\\[3mm]
& \qquad +
\sum^{m-2}_{j=0}c_{m-1,j}\varepsilon^{j}
\big( DA_{m-2-j}(\mbox{\mathversion{bold}$V$})[\mbox{\mathversion{bold}$g$}^{\varepsilon}_{1}(\mbox{\mathversion{bold}$V$})]\big)
\partial^{2(m-1)}_{s}\mbox{\mathversion{bold}$V$}\\[3mm]
& \qquad +
D\mbox{\mathversion{bold}$\varphi $}^{\varepsilon}_{m-1}(\mbox{\mathversion{bold}$V$})
[\mbox{\mathversion{bold}$g$}^{\varepsilon}_{1}(\mbox{\mathversion{bold}$V$})]
+
A_{m-1}(\mbox{\mathversion{bold}$V$})\partial^{2m}_{s}\mbox{\mathversion{bold}$V$}
+
A_{m-1}(\mbox{\mathversion{bold}$V$})\big( \partial^{2(m-1)}_{s}
(|\mbox{\mathversion{bold}$V$}_{s}|^{2}\mbox{\mathversion{bold}$V$})\big)\\[3mm]
& \qquad +
\big( DA_{m-1}(\mbox{\mathversion{bold}$V$})
[\mbox{\mathversion{bold}$r$}^{\varepsilon}_{1}(\mbox{\mathversion{bold}$V$})]\big)
\partial^{2(m-1)}_{s}\mbox{\mathversion{bold}$V$}
+
D\mbox{\mathversion{bold}$\phi $}_{m-1}(\mbox{\mathversion{bold}$V$})
[\mbox{\mathversion{bold}$r$}^{\varepsilon}_{1}(\mbox{\mathversion{bold}$V$})]
\end{align*}
Since \( \mbox{\mathversion{bold}$g$}^{\varepsilon}_{1}(\mbox{\mathversion{bold}$V$})
= \mbox{\mathversion{bold}$V$}\times \mbox{\mathversion{bold}$V$}_{ss}+
\varepsilon ( \mbox{\mathversion{bold}$V$}_{ss}+|\mbox{\mathversion{bold}$V$}_{s}|^{2}\mbox{\mathversion{bold}$V$})\),
we have
\begin{align*}
\partial^{2(m-1)}_{s}\mbox{\mathversion{bold}$g$}^{\varepsilon}_{1}(\mbox{\mathversion{bold}$V$})
&=
\mbox{\mathversion{bold}$V$}\times \partial^{2m}_{s}\mbox{\mathversion{bold}$V$}
+
\varepsilon \partial^{2m}_{s}\mbox{\mathversion{bold}$V$}
+
\sum^{2m-3}_{j=0}
\left(
\begin{array}{c}
2m-2\\
j
\end{array}\right)
\partial^{2m-2-j}_{s}\mbox{\mathversion{bold}$V$}\times 
\partial^{j+2}_{s}\mbox{\mathversion{bold}$V$}\\[3mm]
& \qquad +
\varepsilon \big( \partial^{2(m-1)}_{s}(|\mbox{\mathversion{bold}$V$}_{s}|^{2}\mbox{\mathversion{bold}$V$})\big).
\end{align*}
Substituting the above yields
\begin{align*}
\mbox{\mathversion{bold}$r$}^{\varepsilon}_{m}(\mbox{\mathversion{bold}$V$})
&=
\sum^{m-2}_{j=0}c_{m-1,j}\varepsilon^{j}A_{m-2-j}(\mbox{\mathversion{bold}$V$})
\big( \mbox{\mathversion{bold}$V$}\times \partial^{2m}_{s}\mbox{\mathversion{bold}$V$} 
+ \varepsilon\partial^{2m}_{s}\mbox{\mathversion{bold}$V$}\big)
+ A_{m-1}(\mbox{\mathversion{bold}$V$})\partial^{2m}_{s}\mbox{\mathversion{bold}$V$}\\[3mm]
& \qquad +
\sum^{m-2}_{j=0}c_{m-1,j}\varepsilon^{j}A_{m-2-j}(\mbox{\mathversion{bold}$V$})
\left\{
\sum^{2m-3}_{k=0}
\left(
\begin{array}{c}
2m-2\\
k
\end{array}\right)
\partial^{2(m-1)-k}_{s}\mbox{\mathversion{bold}$V$}\times \partial ^{k+2}_{s}\mbox{\mathversion{bold}$V$}\right.
\\[3mm]
& \hspace*{7cm} \left.
\phantom{
\sum^{m-2}_{j=0}
}
+\varepsilon \big( \partial^{2(m-1)}_{s}(|\mbox{\mathversion{bold}$V$}_{s}|^{2}\mbox{\mathversion{bold}$V$})\big) \right\}
\\[3mm]
& \qquad +
A_{m-1}(\mbox{\mathversion{bold}$V$})\big( \partial^{2(m-1)}_{s}(|\mbox{\mathversion{bold}$V$}_{s}|^{2}\mbox{\mathversion{bold}$V$})\big)
+
D\mbox{\mathversion{bold}$\varphi $}^{\varepsilon}_{m-1}(\mbox{\mathversion{bold}$V$})
[\mbox{\mathversion{bold}$g$}^{\varepsilon}_{1}(\mbox{\mathversion{bold}$V$})]
+
D\mbox{\mathversion{bold}$\phi $}_{m-1}(\mbox{\mathversion{bold}$V$})
[\mbox{\mathversion{bold}$r$}^{\varepsilon}_{1}(\mbox{\mathversion{bold}$V$})]\\[3mm]
& \qquad +\sum^{m-2}_{j=0}c_{m-1,j}\varepsilon^{j}\big( 
DA_{m-2-j}(\mbox{\mathversion{bold}$V$})
[\mbox{\mathversion{bold}$g$}^{\varepsilon}_{1}(\mbox{\mathversion{bold}$V$})]\big)
\partial^{2(m-1)}_{s}\mbox{\mathversion{bold}$V$}\\[3mm]
& \qquad +
\big(D A_{m-1}(\mbox{\mathversion{bold}$V$})[\mbox{\mathversion{bold}$r$}^{\varepsilon}_{1}(\mbox{\mathversion{bold}$V$})]\big)
\partial^{2(m-1)}_{s}\mbox{\mathversion{bold}$V$}.
\end{align*}
The terms with the highest order of derivatives, i.e. terms
containing \( \partial^{2m}_{s}\mbox{\mathversion{bold}$V$}\), can be further calculated as follows.
\begin{align*}
\sum^{m-2}_{j=0}
& c_{m-1,j}\varepsilon^{j}A_{m-2-j}(\mbox{\mathversion{bold}$V$})
\big( \mbox{\mathversion{bold}$V$}\times \partial^{2m}_{s}\mbox{\mathversion{bold}$V$} 
+ \varepsilon\partial^{2m}_{s}\mbox{\mathversion{bold}$V$}\big)
+ A_{m-1}(\mbox{\mathversion{bold}$V$})\partial^{2m}_{s}\mbox{\mathversion{bold}$V$}\\[3mm]
&= \sum^{m-2}_{j=0}c_{m-1,j}\varepsilon^{j}A_{m-1-j}(\mbox{\mathversion{bold}$V$})\partial^{2m}_{s}\mbox{\mathversion{bold}$V$}
+
\sum^{m-2}_{j=0}c_{m-1,j}\varepsilon^{j+1}A_{m-2-j}(\mbox{\mathversion{bold}$V$})\partial^{2m}_{s}\mbox{\mathversion{bold}$V$}
+
A_{m-1}(\mbox{\mathversion{bold}$V$})\partial^{2m}_{s}\mbox{\mathversion{bold}$V$}\\[3mm]
&=
\sum^{m-2}_{j=0}c_{m-1,j}\varepsilon^{j}A_{m-1-j}(\mbox{\mathversion{bold}$V$})\partial^{2m}_{s}\mbox{\mathversion{bold}$V$}
+
\sum^{m-1}_{j=1}c_{m-1,j-1}\varepsilon^{j}A_{m-1-j}(\mbox{\mathversion{bold}$V$})\partial^{2m}_{s}\mbox{\mathversion{bold}$V$}
+
A_{m-1}(\mbox{\mathversion{bold}$V$})\partial^{2m}_{s}\mbox{\mathversion{bold}$V$}\\[3mm]
&= \sum^{m-1}_{j=0}c_{m,j}\varepsilon^{j}A_{m-1-j}(\mbox{\mathversion{bold}$V$})\partial^{2m}_{s}\mbox{\mathversion{bold}$V$},
\end{align*}
where the definition of \( c_{m,j}\) was used in the last equality.
The other terms can also be calculated as follows.
\begin{align*}
\sum^{m-2}_{j=0}&c_{m-1,j}\varepsilon^{j}A_{m-2-j}(\mbox{\mathversion{bold}$V$})
\left\{
\varepsilon \big( \partial^{2(m-1)}_{s}(|\mbox{\mathversion{bold}$V$}_{s}|^{2}\mbox{\mathversion{bold}$V$})\big)
\right\}
+A_{m-1}(\mbox{\mathversion{bold}$V$})\big( \partial^{2(m-1)}_{s}(|\mbox{\mathversion{bold}$V$}_{s}|^{2}\mbox{\mathversion{bold}$V$})\big)
\\[3mm]
& \qquad = \sum ^{m-1}_{j=0}d_{m-1,j}\varepsilon^{j}A_{m-1-j}(\mbox{\mathversion{bold}$V$})
\left\{ \partial^{2(m-1)}_{s}\big( |\mbox{\mathversion{bold}$V$}_{s}|^{2}\mbox{\mathversion{bold}$V$}\big)\right\},\\[5mm]
\sum^{m-2}_{j=0}&c_{m-1,j}\varepsilon^{j}\big( 
DA_{m-2-j}(\mbox{\mathversion{bold}$V$})
[\mbox{\mathversion{bold}$g$}^{\varepsilon}_{1}(\mbox{\mathversion{bold}$V$})]\big)
\partial^{2(m-1)}_{s}\mbox{\mathversion{bold}$V$}
+
\big(D A_{m-1}(\mbox{\mathversion{bold}$V$})[\mbox{\mathversion{bold}$r$}^{\varepsilon}_{1}(\mbox{\mathversion{bold}$V$})]\big)
\partial^{2(m-1)}_{s}\mbox{\mathversion{bold}$V$}\\[3mm]
&\qquad = \sum^{m-2}_{j=0}c_{m-1,j}\varepsilon^{j}\big( 
DA_{m-2-j}(\mbox{\mathversion{bold}$V$})
[\mbox{\mathversion{bold}$r$}^{\varepsilon}_{1}(\mbox{\mathversion{bold}$V$})]\big)
\partial^{2(m-1)}_{s}\mbox{\mathversion{bold}$V$}
+
\big(D A_{m-1}(\mbox{\mathversion{bold}$V$})[\mbox{\mathversion{bold}$r$}^{\varepsilon}_{1}(\mbox{\mathversion{bold}$V$})]\big)
\partial^{2(m-1)}_{s}\mbox{\mathversion{bold}$V$}\\[3mm]
&\hspace*{6cm} 
+
\sum^{m-2}_{j=0}c_{m-1,j}\varepsilon^{j}\big( 
DA_{m-2-j}(\mbox{\mathversion{bold}$V$})
[\mbox{\mathversion{bold}$f$}_{1}(\mbox{\mathversion{bold}$V$})]\big)
\partial^{2(m-1)}_{s}\mbox{\mathversion{bold}$V$}\\[3mm]
& \qquad 
=
\sum^{m-1}_{j=0}d_{m-1,j}\varepsilon^{j}\big( 
DA_{m-1-j}(\mbox{\mathversion{bold}$V$})
[\mbox{\mathversion{bold}$r$}^{\varepsilon}_{1}(\mbox{\mathversion{bold}$V$})]\big)
\partial^{2(m-1)}_{s}\mbox{\mathversion{bold}$V$}\\[3mm]
& \hspace*{6cm}
+
\sum^{m-2}_{j=0}c_{m-1,j}\varepsilon^{j}\big( 
DA_{m-2-j}(\mbox{\mathversion{bold}$V$})
[\mbox{\mathversion{bold}$f$}_{1}(\mbox{\mathversion{bold}$V$})]\big)
\partial^{2(m-1)}_{s}\mbox{\mathversion{bold}$V$},
\end{align*}
where \( \mbox{\mathversion{bold}$g$}^{\varepsilon}_{1}(\mbox{\mathversion{bold}$V$})=
\mbox{\mathversion{bold}$f$}_{1}(\mbox{\mathversion{bold}$V$})+\varepsilon
\mbox{\mathversion{bold}$r$}^{\varepsilon}_{1}(\mbox{\mathversion{bold}$V$}) \) was used.
This finishes the proof. \hfill \( \Box \).

\bigskip
Now we are ready to construct \( \mbox{\mathversion{bold}$h$}^{\varepsilon}\).
We utilize the preceding lemmas to determine the trace of \( \mbox{\mathversion{bold}$h$}^{\varepsilon }\)
at \( s=0,1 \). 
First, we choose \( \mbox{\mathversion{bold}$h$}^{\varepsilon}(0)=\mbox{\mathversion{bold}$h$}^{\varepsilon}(1)
= \mbox{\mathversion{bold}$0$}\). This ensures that \( \mbox{\mathversion{bold}$v$}^{\varepsilon}_{0}\) 
satisfies the \( 0 \)-th order compatibility condition.

From the explicit form of \( \mbox{\mathversion{bold}$v$}^{\varepsilon}_{0}\) we see that
for \( n\in \mathbf{N}\),
\begin{eqnarray}
\partial^{n}_{s}\mbox{\mathversion{bold}$v$}^{\varepsilon}_{0}|_{s=0,1}=
\partial^{n}_{s}\mbox{\mathversion{bold}$v$}_{0}+\partial^{n}_{s}\mbox{\mathversion{bold}$h$}^{\varepsilon}
-(\mbox{\mathversion{bold}$v$}_{0}\cdot 
\partial^{n}_{s}\mbox{\mathversion{bold}$h$}^{\varepsilon})\mbox{\mathversion{bold}$v$}_{0}
+
\mbox{\mathversion{bold}$q$}_{n}(\mbox{\mathversion{bold}$v$}_{0},\mbox{\mathversion{bold}$h$}^{\varepsilon})|_{s=0,1}.
\label{hn}
\end{eqnarray}
Here, \( \mbox{\mathversion{bold}$q$}_{n}(\mbox{\mathversion{bold}$v$}_{0},\mbox{\mathversion{bold}$h$}^{\varepsilon})\)
satisfies
\begin{eqnarray*}
|\mbox{\mathversion{bold}$q$}_{n}(\mbox{\mathversion{bold}$v$}_{0},\mbox{\mathversion{bold}$h$}^{\varepsilon})|
\leq
C(|\mbox{\mathversion{bold}$h$}^{\varepsilon}|+|\mbox{\mathversion{bold}$h$}^{\varepsilon}_{s}|+
\cdots + |\partial^{n-1}_{s}\mbox{\mathversion{bold}$h$}^{\varepsilon}|)
\end{eqnarray*}
if \( |\mbox{\mathversion{bold}$h$}^{\varepsilon}|+|\mbox{\mathversion{bold}$h$}^{\varepsilon}_{s}|+
\cdots + |\partial^{n-1}_{s}\mbox{\mathversion{bold}$h$}^{\varepsilon}|\leq M \), and \( C>0\) depends on 
\( \mbox{\mathversion{bold}$v$}_{0}\) and \( M\).
From (\ref{hn}) and Lemma \ref{fhighlm} we have
\begin{eqnarray*}
\mbox{\mathversion{bold}$f$}_{m}(\mbox{\mathversion{bold}$v$}^{\varepsilon}_{0})\big|_{s=0,1} =
\mbox{\mathversion{bold}$f$}_{m}(\mbox{\mathversion{bold}$v$}_{0})+
A_{m}(\mbox{\mathversion{bold}$v$}^{\varepsilon}_{0})
\big\{
\partial^{2m}_{s}\mbox{\mathversion{bold}$h$}^{\varepsilon}-(\mbox{\mathversion{bold}$v$}_{0}\cdot 
\partial^{2m}_{s}\mbox{\mathversion{bold}$h$}^{\varepsilon})\mbox{\mathversion{bold}$v$}_{0}
\big\}
+
\mbox{\mathversion{bold}$F$}_{m}(\mbox{\mathversion{bold}$v$}_{0},\mbox{\mathversion{bold}$h$}^{\varepsilon})\big|_{s=0,1},
\end{eqnarray*}
for \( m\geq 1\). Here, \( \mbox{\mathversion{bold}$F$}_{m}(\mbox{\mathversion{bold}$v$}_{0},\mbox{\mathversion{bold}$h$}^{\varepsilon})\)
satisfies 
\begin{eqnarray*}
|\mbox{\mathversion{bold}$F$}_{m}(\mbox{\mathversion{bold}$v$}_{0},\mbox{\mathversion{bold}$h$}^{\varepsilon})|
\leq
C(|\mbox{\mathversion{bold}$h$}^{\varepsilon}|+|\mbox{\mathversion{bold}$h$}^{\varepsilon}_{s}|+
\cdots + |\partial^{2m-1}_{s}\mbox{\mathversion{bold}$h$}^{\varepsilon}|)
\end{eqnarray*}
if \( |\mbox{\mathversion{bold}$h$}^{\varepsilon}|+|\mbox{\mathversion{bold}$h$}^{\varepsilon}_{s}|+
\cdots + |\partial^{2m-1}_{s}\mbox{\mathversion{bold}$h$}^{\varepsilon}|\leq M \), and \( C>0\) depends on 
\( \mbox{\mathversion{bold}$v$}_{0}\) and \( M\).
Since we chose \( \mbox{\mathversion{bold}$h$}^{\varepsilon}|_{s=0,1} = \mbox{\mathversion{bold}$0$}\), it follows that
\( \mbox{\mathversion{bold}$v$}^{\varepsilon}_{0}|_{s=0,1} = 
\mbox{\mathversion{bold}$v$}_{0}|_{s=0,1} \), and thus we have
\begin{eqnarray}
\mbox{\mathversion{bold}$f$}_{m}(\mbox{\mathversion{bold}$v$}^{\varepsilon}_{0})|_{s=0,1} = 
\mbox{\mathversion{bold}$f$}_{m}(\mbox{\mathversion{bold}$v$}_{0}) + A_{m}(\mbox{\mathversion{bold}$v$}_{0})
\partial^{2m}_{s}\mbox{\mathversion{bold}$h$}^{\varepsilon} + 
\mbox{\mathversion{bold}$F$}_{m}(\mbox{\mathversion{bold}$v$}_{0},\mbox{\mathversion{bold}$h$}^{\varepsilon})\big|_{s=0,1}.
\label{ccind}
\end{eqnarray}
Now we prove by induction that the differential coefficients of \( \mbox{\mathversion{bold}$h$}^{\varepsilon} \)
at \( s=0,1\) can be chosen so that
\( \mbox{\mathversion{bold}$v$}^{\varepsilon}_{0}\) satisfies the compatibility conditions for (\ref{rnl}) up to order \( m\), and
the differential coefficients of \( \mbox{\mathversion{bold}$h$}^{\varepsilon} \) at \( s=0,1 \) are \( O(\varepsilon ) \).
We have already chosen \(  \mbox{\mathversion{bold}$h$}^{\varepsilon}|_{s=0,1} = 
\mbox{\mathversion{bold}$0$}\), which insures that \(  \mbox{\mathversion{bold}$v$}^{\varepsilon}_{0}\) satisfies the 
\( 0\)-th order compatibility condition. 
Suppose that for a \( m\geq 1 \), the differential coefficients 
\( \partial^{j}_{s}\mbox{\mathversion{bold}$h$}^{\varepsilon}|_{s=0,1}\) for \( 0\leq j\leq 2(m-1) \) have been chosen in a way 
such that they are \( O(\varepsilon ) \) and \( \mbox{\mathversion{bold}$v$}^{\varepsilon}_{0}\) satisfies the 
compatibility conditions up to order \( m-1\).
First, we choose \(  \partial^{2m-1}_{s}\mbox{\mathversion{bold}$h$}^{\varepsilon}|_{s=0,1}= \mbox{\mathversion{bold}$0$}\).
Then from Lemma \ref{differ}, Lemma \ref{rmhigh}, and (\ref{ccind}) we have
\begin{align*}
\mbox{\mathversion{bold}$g$}^{\varepsilon}_{m}( \mbox{\mathversion{bold}$v$}^{\varepsilon}_{0})|_{s=0,1}
&=
\mbox{\mathversion{bold}$f$}_{m}(\mbox{\mathversion{bold}$v$}^{\varepsilon}_{0})
+\varepsilon \mbox{\mathversion{bold}$r$}^{\varepsilon}_{m}(\mbox{\mathversion{bold}$v$}^{\varepsilon}_{0})|_{s=0,1} \\[3mm]
&=
\mbox{\mathversion{bold}$f$}_{m}(\mbox{\mathversion{bold}$v$}_{0})+
A_{m}(\mbox{\mathversion{bold}$v$}_{0})\partial^{2m}_{s}\mbox{\mathversion{bold}$h$}^{\varepsilon}
+\mbox{\mathversion{bold}$F$}_{m}(\mbox{\mathversion{bold}$v$}_{0},\mbox{\mathversion{bold}$h$}^{\varepsilon})\\[3mm]
&\qquad +
\varepsilon 
\bigg\{ \sum^{m-1}_{j=0}c_{m,j}\varepsilon^{j}A_{m-1-j}(\mbox{\mathversion{bold}$v$}_{0})
\partial^{2m}_{s}\mbox{\mathversion{bold}$v$}^{\varepsilon}_{0}
+\mbox{\mathversion{bold}$\varphi$}_{m}^{\varepsilon}
(\mbox{\mathversion{bold}$v$}^{\varepsilon}_{0})\bigg\}
\bigg|_{s=0,1}\\[3mm]
&=A_{m}(\mbox{\mathversion{bold}$v$}_{0})\partial^{2m}_{s}\mbox{\mathversion{bold}$h$}^{\varepsilon}
+\mbox{\mathversion{bold}$F$}_{m}(\mbox{\mathversion{bold}$v$}_{0},\mbox{\mathversion{bold}$h$}^{\varepsilon})\\[3mm]
& \qquad +
\varepsilon 
\bigg\{ \sum^{m-1}_{j=0}c_{m,j}\varepsilon^{j}A_{m-1-j}(\mbox{\mathversion{bold}$v$}_{0})
\partial^{2m}_{s}\mbox{\mathversion{bold}$v$}^{\varepsilon}_{0}
+\mbox{\mathversion{bold}$\varphi$}_{m}^{\varepsilon}
(\mbox{\mathversion{bold}$v$}^{\varepsilon}_{0})\bigg\}\bigg|_{s=0,1},
\end{align*}
where we have used the fact that \( \mbox{\mathversion{bold}$v$}_{0} \) satisfies the \( m\)-th order compatibility condition 
for (\ref{vslant}), i.e. \( \mbox{\mathversion{bold}$f$}_{m}(\mbox{\mathversion{bold}$v$}_{0})|_{s=0,1}=
\mbox{\mathversion{bold}$0$}\). 
From (\ref{hn}), we have
\begin{align*}
\sum^{m-1}_{j=0}c_{m,j}\varepsilon^{j}A_{m-1-j}(\mbox{\mathversion{bold}$v$}_{0})
\partial^{2m}_{s}\mbox{\mathversion{bold}$v$}^{\varepsilon}_{0}
&=
\sum^{m-1}_{j=0}c_{m,j}\varepsilon^{j}A_{m-1-j}(\mbox{\mathversion{bold}$v$}_{0})
\partial^{2m}_{s}\mbox{\mathversion{bold}$v$}_{0}
+
\sum^{m-1}_{j=0}c_{m,j}\varepsilon^{j}A_{m-1-j}(\mbox{\mathversion{bold}$v$}_{0})
\partial^{2m}_{s}\mbox{\mathversion{bold}$h$}^{\varepsilon}_{0}\\[3mm]
& \qquad \qquad -
\sum^{m-1}_{j=0}c_{m,j}\varepsilon^{j}A_{m-1-j}(\mbox{\mathversion{bold}$v$}_{0})
\{(\mbox{\mathversion{bold}$v$}_{0}
\cdot \partial^{2m}_{s}\mbox{\mathversion{bold}$h$}^{\varepsilon})
\mbox{\mathversion{bold}$v$}_{0}\}\\[3mm]
& \qquad \qquad  +
\sum^{m-1}_{j=0}c_{m,j}\varepsilon^{j}A_{m-1-j}(\mbox{\mathversion{bold}$v$}_{0})
\mbox{\mathversion{bold}$q$}_{2m}(\mbox{\mathversion{bold}$v$}_{0},
\mbox{\mathversion{bold}$h$}^{\varepsilon}_{0}),
\end{align*}
which yields
\begin{align}
\label{vperp}
\mbox{\mathversion{bold}$g$}^{\varepsilon}_{m}
( \mbox{\mathversion{bold}$v$}^{\varepsilon}_{0})|_{s=0,1}
&=
A_{m}(\mbox{\mathversion{bold}$v$}_{0})\partial^{2m}_{s}\mbox{\mathversion{bold}$h$}^{\varepsilon}
+
\varepsilon \sum^{m-1}_{j=0}c_{m,j}\varepsilon^{j}A_{m-1-j}(\mbox{\mathversion{bold}$v$}_{0})
\partial^{2m}_{s}\mbox{\mathversion{bold}$h$}^{\varepsilon}
+
\mbox{\mathversion{bold}$F$}_{m}(\mbox{\mathversion{bold}$v$}_{0},\mbox{\mathversion{bold}$h$}^{\varepsilon}) \\[3mm]
& \qquad \qquad -
c_{m,m-1}\varepsilon^{m}
\{(\mbox{\mathversion{bold}$v$}_{0}
\cdot \partial^{2m}_{s}\mbox{\mathversion{bold}$h$}^{\varepsilon})
\mbox{\mathversion{bold}$v$}_{0}\} \nonumber \\[3mm]
& \qquad \qquad +
\varepsilon \sum^{m-1}_{j=0}c_{m,j}\varepsilon^{j}A_{m-1-j}(\mbox{\mathversion{bold}$v$}_{0})
\partial^{2m}_{s}\mbox{\mathversion{bold}$v$}_{0}\nonumber\\[3mm]
& \qquad \qquad +
\varepsilon \sum^{m-1}_{j=0}c_{m,j}\varepsilon^{j}A_{m-1-j}(\mbox{\mathversion{bold}$v$}_{0})
\mbox{\mathversion{bold}$q$}_{2m}(\mbox{\mathversion{bold}$v$}_{0},
\mbox{\mathversion{bold}$h$}^{\varepsilon}_{0})
\bigg|_{s=0,1}\nonumber\\[3mm]
&=:A_{m}(\mbox{\mathversion{bold}$v$}_{0})\partial^{2m}_{s}\mbox{\mathversion{bold}$h$}^{\varepsilon}
+
\varepsilon \sum^{m-1}_{j=0}c_{m,j}\varepsilon^{j}A_{m-1-j}(\mbox{\mathversion{bold}$v$}_{0})
\partial^{2m}_{s}\mbox{\mathversion{bold}$h$}^{\varepsilon}\nonumber\\[3mm]
& \qquad \qquad -
c_{m,m-1}\varepsilon^{m}
\{(\mbox{\mathversion{bold}$v$}_{0}
\cdot \partial^{2m}_{s}\mbox{\mathversion{bold}$h$}^{\varepsilon})
\mbox{\mathversion{bold}$v$}_{0}\}
+
\mbox{\mathversion{bold}$F$}_{m}(\mbox{\mathversion{bold}$v$}_{0},\mbox{\mathversion{bold}$h$}^{\varepsilon})
+
\varepsilon \mbox{\mathversion{bold}$G$}_{m}(\mbox{\mathversion{bold}$v$}_{0},\mbox{\mathversion{bold}$h$}^{\varepsilon})|_{s=0,1}.\nonumber
\end{align}
We note that 
\( \mbox{\mathversion{bold}$F$}_{m}(\mbox{\mathversion{bold}$v$}_{0},\mbox{\mathversion{bold}$h$}^{\varepsilon})\)
and
\( \mbox{\mathversion{bold}$G$}_{m}(\mbox{\mathversion{bold}$v$}_{0},\mbox{\mathversion{bold}$h$}^{\varepsilon}) \)
only contain terms with \( s\)-derivatives of 
\( \mbox{\mathversion{bold}$v$}_{0} \) and 
\( \mbox{\mathversion{bold}$h$}^{\varepsilon}\) less than or equal to 
\( 2m-1\).
From Lemma \ref{lmzero}, we see that
\begin{eqnarray*}
\sum^{m}_{k=0}
\left(
\begin{array}{c}
m\\
k
\end{array}\right)
\mbox{\mathversion{bold}$g$}^{\varepsilon}_{k}
(\mbox{\mathversion{bold}$v$}^{\varepsilon}_{0})\cdot
\mbox{\mathversion{bold}$g$}^{\varepsilon}_{m-k}
(\mbox{\mathversion{bold}$v$}^{\varepsilon}_{0}) \equiv 0,
\end{eqnarray*}
and the assumption of induction implies that
\( \mbox{\mathversion{bold}$g$}^{\varepsilon}_{k}
(\mbox{\mathversion{bold}$v$}^{\varepsilon}_{0})|_{s=0,1} = 
\mbox{\mathversion{bold}$0$}\) for \( 1\leq k \leq m-1 \).
Hence we have
\begin{eqnarray*}
0 \equiv \mbox{\mathversion{bold}$g$}^{\varepsilon}_{0}(\mbox{\mathversion{bold}$v$}^{\varepsilon}_{0})\cdot 
\mbox{\mathversion{bold}$g$}^{\varepsilon}_{m}(\mbox{\mathversion{bold}$v$}^{\varepsilon}_{0})|_{s=0,1}
=
\mbox{\mathversion{bold}$v$}_{0}\cdot
\mbox{\mathversion{bold}$g$}^{\varepsilon}_{m}(\mbox{\mathversion{bold}$v$}^{\varepsilon}_{0})|_{s=0,1}.
\end{eqnarray*}
Substituting (\ref{vperp}) into \( 
\mbox{\mathversion{bold}$g$}^{\varepsilon}_{m}(\mbox{\mathversion{bold}$v$}^{\varepsilon}_{0}) \)
yields
\begin{align}
\label{vperp2}
0 &\equiv 
\mbox{\mathversion{bold}$v$}_{0}\cdot \bigg\{
A_{m}(\mbox{\mathversion{bold}$v$}_{0})\partial^{2m}_{s}\mbox{\mathversion{bold}$h$}^{\varepsilon}
+
\varepsilon \sum^{m-1}_{j=0}c_{m,j}\varepsilon^{j}A_{m-1-j}(\mbox{\mathversion{bold}$v$}_{0})
\partial^{2m}_{s}\mbox{\mathversion{bold}$h$}^{\varepsilon}\\[3mm]
& \qquad \qquad -
c_{m,m-1}\varepsilon^{m}
\{(\mbox{\mathversion{bold}$v$}_{0}
\cdot \partial^{2m}_{s}\mbox{\mathversion{bold}$h$}^{\varepsilon})
\mbox{\mathversion{bold}$v$}_{0}\}
+
\mbox{\mathversion{bold}$F$}_{m}(\mbox{\mathversion{bold}$v$}_{0},\mbox{\mathversion{bold}$h$}^{\varepsilon})
+
\varepsilon \mbox{\mathversion{bold}$G$}_{m}(\mbox{\mathversion{bold}$v$}_{0},\mbox{\mathversion{bold}$h$}^{\varepsilon})\bigg\} \bigg|_{s=0,1}\nonumber \\[3mm]
&=
\mbox{\mathversion{bold}$v$}_{0}\cdot
\bigg\{ \mbox{\mathversion{bold}$F$}_{m}(\mbox{\mathversion{bold}$v$}_{0},\mbox{\mathversion{bold}$h$}^{\varepsilon})
+
\varepsilon \mbox{\mathversion{bold}$G$}_{m}(\mbox{\mathversion{bold}$v$}_{0},\mbox{\mathversion{bold}$h$}^{\varepsilon})\bigg\}\bigg|_{s=0,1}. \nonumber
\end{align}
Furthermore, since
\begin{eqnarray*}
\mbox{\mathversion{bold}$v$}_{0}\times(\mbox{\mathversion{bold}$v$}_{0}\times 
\partial^{2m}_{s}\mbox{\mathversion{bold}$h$}^{\varepsilon})
=
(\mbox{\mathversion{bold}$v$}_{0}\cdot 
\partial^{2m}_{s}\mbox{\mathversion{bold}$h$}^{\varepsilon})\mbox{\mathversion{bold}$v$}_{0}
-
\partial^{2m}_{s}\mbox{\mathversion{bold}$h$}^{\varepsilon},
\end{eqnarray*}
we see that
\begin{eqnarray*}
A_{k}(\mbox{\mathversion{bold}$v$}_{0})
\partial^{2m}_{s}\mbox{\mathversion{bold}$h$}^{\varepsilon}
=
\left\{
\begin{array}{ll}
(-1)^{l}\mbox{\mathversion{bold}$v$}_{0}\times
\partial^{2m}_{s}\mbox{\mathversion{bold}$h$}^{\varepsilon}, & \mbox{when} \ k=2l+1, \\[3mm]
(-1)^{l+1}\big( ( \mbox{\mathversion{bold}$v$}_{0}\cdot
\partial^{2m}_{s}\mbox{\mathversion{bold}$h$}^{\varepsilon})\mbox{\mathversion{bold}$v$}_{0}
-
\partial^{2m}_{s}\mbox{\mathversion{bold}$h$}^{\varepsilon}\big),
& \mbox{when} \ k=2l,
\end{array}\right.
\end{eqnarray*}
where \( l\in \mathbf{N}\cup \{0\}\).
From (\ref{vperp}), we see that when \( m=2l+1 \),
\begin{align*}
\mbox{\mathversion{bold}$g$}^{\varepsilon}_{m}
( \mbox{\mathversion{bold}$v$}^{\varepsilon}_{0})\bigg|_{s=0,1}
&=
(-1)^{l}\mbox{\mathversion{bold}$v$}_{0}\times 
\partial^{2m}_{s}\mbox{\mathversion{bold}$h$}^{\varepsilon}
+
\varepsilon \left\{
\sum^{l}_{j=0}(-1)^{j+1}c_{m,2j}\varepsilon^{2j}
\big( ( \mbox{\mathversion{bold}$v$}_{0}\cdot
\partial^{2m}_{s}\mbox{\mathversion{bold}$h$}^{\varepsilon})\mbox{\mathversion{bold}$v$}_{0}
-
\partial^{2m}_{s}\mbox{\mathversion{bold}$h$}^{\varepsilon}\big)\right. \\[3mm]
& \left. \qquad \qquad +
\sum^{l-1}_{j=0}(-1)^{j}c_{m,2j+1}\varepsilon^{2j+1}
\mbox{\mathversion{bold}$v$}_{0}\times \partial^{2m}_{s}
\mbox{\mathversion{bold}$h$}^{\varepsilon}\right\} \\[3mm]
& \qquad -
c_{m,m-1}\varepsilon^{m}
\{(\mbox{\mathversion{bold}$v$}_{0}
\cdot \partial^{2m}_{s}\mbox{\mathversion{bold}$h$}^{\varepsilon})
\mbox{\mathversion{bold}$v$}_{0}\}
+
\mbox{\mathversion{bold}$F$}_{m}(\mbox{\mathversion{bold}$v$}_{0},\mbox{\mathversion{bold}$h$}^{\varepsilon})
+
\varepsilon \mbox{\mathversion{bold}$G$}_{m}(\mbox{\mathversion{bold}$v$}_{0},\mbox{\mathversion{bold}$h$}^{\varepsilon})\bigg|_{s=0,1} \\[3mm]
&=
\left(
(-1)^{l}+ \varepsilon \sum^{l-1}_{j=0}(-1)^{j}c_{m,2j+1}\varepsilon^{2j+1}\right)
\mbox{\mathversion{bold}$v$}_{0}\times 
\partial^{2m}_{s}\mbox{\mathversion{bold}$h$}^{\varepsilon}\\[3mm]
& \qquad +
\left( 
\varepsilon \sum^{l}_{j=0}(-1)^{j}c_{m,2j}\varepsilon^{2j}\right)
\partial^{2m}_{s}\mbox{\mathversion{bold}$h$}^{\varepsilon} \\[3mm]
& \qquad +
\varepsilon 
\sum^{l}_{j=0}(-1)^{j+1}c_{m,2j}\varepsilon^{2j}
( \mbox{\mathversion{bold}$v$}_{0}\cdot
\partial^{2m}_{s}\mbox{\mathversion{bold}$h$}^{\varepsilon})\mbox{\mathversion{bold}$v$}_{0}
-
c_{m,m-1}\varepsilon^{m}
\{(\mbox{\mathversion{bold}$v$}_{0}
\cdot \partial^{2m}_{s}\mbox{\mathversion{bold}$h$}^{\varepsilon})
\mbox{\mathversion{bold}$v$}_{0}\}\\[3mm]
& \qquad +
\mbox{\mathversion{bold}$F$}_{m}(\mbox{\mathversion{bold}$v$}_{0},\mbox{\mathversion{bold}$h$}^{\varepsilon})
+
\varepsilon \mbox{\mathversion{bold}$G$}_{m}(\mbox{\mathversion{bold}$v$}_{0},\mbox{\mathversion{bold}$h$}^{\varepsilon})\bigg|_{s=0,1}
\end{align*}
holds.
At this point, we choose
\begin{align*}
\partial^{2m}_{s}\mbox{\mathversion{bold}$h$}^{\varepsilon}|_{s=0,1}
=
\mbox{\mathversion{bold}$v$}_{0}\times 
\{ B(\mbox{\mathversion{bold}$v$}_{0})^{-1}\big( 
\mbox{\mathversion{bold}$F$}_{m}(\mbox{\mathversion{bold}$v$}_{0},\mbox{\mathversion{bold}$h$}^{\varepsilon})
+
\varepsilon \mbox{\mathversion{bold}$G$}_{m}(\mbox{\mathversion{bold}$v$}_{0},\mbox{\mathversion{bold}$h$}^{\varepsilon})
\big)\}|_{s=0,1},
\end{align*}
where \( B(\mbox{\mathversion{bold}$v$}_{0}) \) is the \( 3\times 3 \) matrix given by
\begin{align*}
\left(
(-1)^{l}+ \varepsilon \sum^{l-1}_{j=0}(-1)^{j}c_{m,2j+1}\varepsilon^{2j+1}\right)
I_{3}
-
\left( 
\varepsilon \sum^{l}_{j=0}(-1)^{j}c_{m,2j}\varepsilon^{2j}\right)
A(\mbox{\mathversion{bold}$v$}_{0}),
\end{align*}
where \( I_{3}\) is the \( 3\times 3 \) identity matrix, and 
\( A(\mbox{\mathversion{bold}$v$}_{0})\) 
is the representation matrix of \( A_{1}(\mbox{\mathversion{bold}$v$}_{0})\).
The inverse of \( B(\mbox{\mathversion{bold}$v$}_{0})\) 
exists for sufficiently small \( \varepsilon \), 
where the smallness depends only on \( \mbox{\mathversion{bold}$v$}_{0}\).
More precisely, there exists \( \varepsilon _{\ast , 1}>0\) such that 
for any \( \varepsilon \in (0,\varepsilon _{\ast,1}]\), the inverse matrix can be expressed as
\begin{eqnarray*}
B(\mbox{\mathversion{bold}$v$}_{0})^{-1}
=
\left(
(-1)^{l}+ \varepsilon M_{1}\right)
(I_{3}
+
D(\mbox{\mathversion{bold}$v$}_{0})),
\end{eqnarray*}
where
\begin{align*}
M_{1}&=\sum^{l-1}_{j=0}(-1)^{j}c_{m,2j+1}\varepsilon^{2j+1}, \\[3mm]
M_{2}&= \bigg( \varepsilon \sum^{l}_{j=0}(-1)^{j}c_{m,2j}\varepsilon^{2j}\bigg)
\bigg(
(-1)^{l}+ \varepsilon M_{1}\bigg)^{-1}, \\[3mm]
D(\mbox{\mathversion{bold}$v$}_{0})
&=
\sum^{\infty}_{k=1}M_{2}^{k}A(\mbox{\mathversion{bold}$v$}_{0})^{k}.
\end{align*}
Since \( M_{2}=O(\varepsilon) \) for small \( \varepsilon \), 
the above infinite sum absolutely converges.
Noting that by this choice of 
\( \partial^{2m}_{s}\mbox{\mathversion{bold}$h$}^{\varepsilon}|_{s=0,1}\), we have
\( \mbox{\mathversion{bold}$v$}_{0}\cdot \partial^{2m}_{s}
\mbox{\mathversion{bold}$h$}^{\varepsilon}|_{s=0,1}=0\), and we see that
\newpage
\begin{align*}
\mbox{\mathversion{bold}$g$}^{\varepsilon}_{m}
( \mbox{\mathversion{bold}$v$}^{\varepsilon}_{0})\bigg|_{s=0,1}
&=
\left(
(-1)^{l}+ M_{1}\right)
\mbox{\mathversion{bold}$v$}_{0}\times
\big\{
\mbox{\mathversion{bold}$v$}_{0}\times 
\big( B(\mbox{\mathversion{bold}$v$}_{0})^{-1}( 
\mbox{\mathversion{bold}$F$}_{m}(\mbox{\mathversion{bold}$v$}_{0},\mbox{\mathversion{bold}$h$}^{\varepsilon})
+
\varepsilon \mbox{\mathversion{bold}$G$}_{m}(\mbox{\mathversion{bold}$v$}_{0},\mbox{\mathversion{bold}$h$}^{\varepsilon})
)\big) \big\}\\[3mm]
& \qquad +
\big( (-1)^{l}+\varepsilon M_{1}\big)M_{2}
\mbox{\mathversion{bold}$v$}_{0}\times 
\{ B(\mbox{\mathversion{bold}$v$}_{0})^{-1}\big( 
\mbox{\mathversion{bold}$F$}_{m}(\mbox{\mathversion{bold}$v$}_{0},\mbox{\mathversion{bold}$h$}^{\varepsilon})
+
\varepsilon \mbox{\mathversion{bold}$G$}_{m}(\mbox{\mathversion{bold}$v$}_{0},\mbox{\mathversion{bold}$h$}^{\varepsilon})
\big)\}\\[3mm]
& \qquad +
\mbox{\mathversion{bold}$F$}_{m}(\mbox{\mathversion{bold}$v$}_{0},\mbox{\mathversion{bold}$h$}^{\varepsilon})
+
\varepsilon \mbox{\mathversion{bold}$G$}_{m}(\mbox{\mathversion{bold}$v$}_{0},\mbox{\mathversion{bold}$h$}^{\varepsilon})\bigg|_{s=0,1}\\[3mm]
&=
-\left(
(-1)^{l}+ M_{1}\right) 
\big( B(\mbox{\mathversion{bold}$v$}_{0})^{-1}( 
\mbox{\mathversion{bold}$F$}_{m}(\mbox{\mathversion{bold}$v$}_{0},\mbox{\mathversion{bold}$h$}^{\varepsilon})
+
\varepsilon \mbox{\mathversion{bold}$G$}_{m}(\mbox{\mathversion{bold}$v$}_{0},\mbox{\mathversion{bold}$h$}^{\varepsilon})
)\big)\\[3mm]
& \qquad +
\big( (-1)^{l}+\varepsilon M_{1}\big)M_{2}
A(\mbox{\mathversion{bold}$v$}_{0})
\{ B(\mbox{\mathversion{bold}$v$}_{0})^{-1}\big( 
\mbox{\mathversion{bold}$F$}_{m}(\mbox{\mathversion{bold}$v$}_{0},\mbox{\mathversion{bold}$h$}^{\varepsilon})
+
\varepsilon \mbox{\mathversion{bold}$G$}_{m}(\mbox{\mathversion{bold}$v$}_{0},\mbox{\mathversion{bold}$h$}^{\varepsilon})
\big)\}\\[3mm]
& \qquad + \big( (-1)^{l}+M_{1}\big)\bigg\{
\mbox{\mathversion{bold}$v$}_{0}\cdot 
\bigg( B(\mbox{\mathversion{bold}$v$}_{0})^{-1}( 
\mbox{\mathversion{bold}$F$}_{m}(\mbox{\mathversion{bold}$v$}_{0},\mbox{\mathversion{bold}$h$}^{\varepsilon})
+
\varepsilon \mbox{\mathversion{bold}$G$}_{m}(\mbox{\mathversion{bold}$v$}_{0},\mbox{\mathversion{bold}$h$}^{\varepsilon})
) \bigg)\bigg\}\mbox{\mathversion{bold}$v$}_{0}\\[3mm]
& \qquad +
\mbox{\mathversion{bold}$F$}_{m}(\mbox{\mathversion{bold}$v$}_{0},\mbox{\mathversion{bold}$h$}^{\varepsilon})
+
\varepsilon \mbox{\mathversion{bold}$G$}_{m}(\mbox{\mathversion{bold}$v$}_{0},\mbox{\mathversion{bold}$h$}^{\varepsilon})\bigg|_{s=0,1}\\[3mm]
& =
-B(\mbox{\mathversion{bold}$v$}_{0})
B(\mbox{\mathversion{bold}$v$}_{0})^{-1}\big( 
\mbox{\mathversion{bold}$F$}_{m}(\mbox{\mathversion{bold}$v$}_{0},\mbox{\mathversion{bold}$h$}^{\varepsilon})
+
\varepsilon \mbox{\mathversion{bold}$G$}_{m}(\mbox{\mathversion{bold}$v$}_{0},\mbox{\mathversion{bold}$h$}^{\varepsilon})\big)\\[3mm]
& \qquad + \big( (-1)^{l}+M_{1}\big)\bigg\{
\mbox{\mathversion{bold}$v$}_{0}\cdot 
\bigg( B(\mbox{\mathversion{bold}$v$}_{0})^{-1}( 
\mbox{\mathversion{bold}$F$}_{m}(\mbox{\mathversion{bold}$v$}_{0},\mbox{\mathversion{bold}$h$}^{\varepsilon})
+
\varepsilon \mbox{\mathversion{bold}$G$}_{m}(\mbox{\mathversion{bold}$v$}_{0},\mbox{\mathversion{bold}$h$}^{\varepsilon})
) \bigg)\bigg\}\mbox{\mathversion{bold}$v$}_{0}\\[3mm]
& \qquad +
\mbox{\mathversion{bold}$F$}_{m}(\mbox{\mathversion{bold}$v$}_{0},\mbox{\mathversion{bold}$h$}^{\varepsilon})
+
\varepsilon \mbox{\mathversion{bold}$G$}_{m}(\mbox{\mathversion{bold}$v$}_{0},\mbox{\mathversion{bold}$h$}^{\varepsilon})\bigg|_{s=0,1}\\[3mm]
&= 
\big( (-1)^{l}+M_{1}\big)\bigg\{
\mbox{\mathversion{bold}$v$}_{0}\cdot 
\bigg( B(\mbox{\mathversion{bold}$v$}_{0})^{-1}( 
\mbox{\mathversion{bold}$F$}_{m}(\mbox{\mathversion{bold}$v$}_{0},\mbox{\mathversion{bold}$h$}^{\varepsilon})
+
\varepsilon \mbox{\mathversion{bold}$G$}_{m}(\mbox{\mathversion{bold}$v$}_{0},\mbox{\mathversion{bold}$h$}^{\varepsilon})
) \bigg)\bigg\}\mbox{\mathversion{bold}$v$}_{0}\bigg|_{s=0,1}.
\end{align*}
Since \( B(\mbox{\mathversion{bold}$v$}_{0})^{-1}=
I_{3}+D(\mbox{\mathversion{bold}$v$}_{0})\), we have
\begin{align*}
\mbox{\mathversion{bold}$v$}_{0}\cdot 
\bigg( B( \mbox{\mathversion{bold}$v$}_{0})^{-1}\big( 
\mbox{\mathversion{bold}$F$}_{m}(\mbox{\mathversion{bold}$v$}_{0},\mbox{\mathversion{bold}$h$}^{\varepsilon})
&+
\varepsilon \mbox{\mathversion{bold}$G$}_{m}(\mbox{\mathversion{bold}$v$}_{0},\mbox{\mathversion{bold}$h$}^{\varepsilon})
\big) \bigg)\bigg|_{s=0,1}\\[3mm]
&=
\mbox{\mathversion{bold}$v$}_{0}\cdot 
\bigg( 
\mbox{\mathversion{bold}$F$}_{m}(\mbox{\mathversion{bold}$v$}_{0},\mbox{\mathversion{bold}$h$}^{\varepsilon})
+
\varepsilon \mbox{\mathversion{bold}$G$}_{m}(\mbox{\mathversion{bold}$v$}_{0},\mbox{\mathversion{bold}$h$}^{\varepsilon})
) \bigg)\\[3mm]
& \qquad +
\mbox{\mathversion{bold}$v$}_{0}\cdot D(\mbox{\mathversion{bold}$v$}_{0})
\bigg( 
\mbox{\mathversion{bold}$F$}_{m}(\mbox{\mathversion{bold}$v$}_{0},\mbox{\mathversion{bold}$h$}^{\varepsilon})
+
\varepsilon \mbox{\mathversion{bold}$G$}_{m}(\mbox{\mathversion{bold}$v$}_{0},\mbox{\mathversion{bold}$h$}^{\varepsilon})
) \bigg)\bigg|_{s=0,1}.
\end{align*}
From (\ref{vperp2}), the first term on the right-hand side is \( 0\). 
Furthermore, since
\begin{eqnarray*}
D(\mbox{\mathversion{bold}$v$}_{0})=
\sum^{\infty}_{k=1}M_{2}^{k}A(\mbox{\mathversion{bold}$v$}_{0})^{k}
=
A(\mbox{\mathversion{bold}$v$}_{0})
\bigg( \sum^{\infty}_{k=0}M_{2}^{k+1}A(\mbox{\mathversion{bold}$v$}_{0})^{k} \bigg),
\end{eqnarray*}
and \( A(\mbox{\mathversion{bold}$v$}_{0})\mbox{\mathversion{bold}$W$}=
\mbox{\mathversion{bold}$v$}_{0}\times \mbox{\mathversion{bold}$W$}\)
for any vector \( \mbox{\mathversion{bold}$W$}\), 
the second term on the right-hand side is also \( 0 \). 
Hence we have
\begin{eqnarray*}
\mbox{\mathversion{bold}$g$}^{\varepsilon}_{m}
( \mbox{\mathversion{bold}$v$}^{\varepsilon}_{0})|_{s=0,1}
=\mbox{\mathversion{bold}$0$},
\end{eqnarray*}
which means that \( \mbox{\mathversion{bold}$v$}^{\varepsilon}_{0} \)
satisfies the \( m\)-th order compatibility condition. From the assumption of induction,
\( \mbox{\mathversion{bold}$F$}_{m}(\mbox{\mathversion{bold}$v$}_{0},\mbox{\mathversion{bold}$h$}^{\varepsilon})
+
\varepsilon \mbox{\mathversion{bold}$G$}_{m}(\mbox{\mathversion{bold}$v$}_{0},\mbox{\mathversion{bold}$h$}^{\varepsilon})\) is \( O(\varepsilon )\) and 
\( B( \mbox{\mathversion{bold}$v$}_{0})^{-1}\) is \( O(1)\), which implies that
\( \partial^{2m}_{s}\mbox{\mathversion{bold}$h$}^{\varepsilon}|_{s=0,1} \)
is \( O(\varepsilon)\). This finishes the proof for odd \( m\).
The case \( m=2l \ (l\geq 1) \) can be treated in the same way, and one can prove that
there exists \( \varepsilon_{\ast ,2}>0\) such that for any \( \varepsilon \in 
(0,\varepsilon_{\ast ,2}]\),
choosing 
\begin{eqnarray*}
\partial^{2m}_{s}\mbox{\mathversion{bold}$h$}^{\varepsilon}|_{s=0,1}
=\widetilde{B}(\mbox{\mathversion{bold}$v$}_{0})^{-1}
\big( \mbox{\mathversion{bold}$F$}_{m}(\mbox{\mathversion{bold}$v$}_{0},\mbox{\mathversion{bold}$h$}^{\varepsilon})
+
\varepsilon \mbox{\mathversion{bold}$G$}_{m}(\mbox{\mathversion{bold}$v$}_{0},\mbox{\mathversion{bold}$h$}^{\varepsilon})\big)|_{s=0,1}
\end{eqnarray*}
is sufficient. Here, 
\begin{align*}
\widetilde{B}(\mbox{\mathversion{bold}$v$}_{0})
=
\bigg( (-1)^{l}+\varepsilon \sum^{l-1}_{j=0}(-1)^{j}
c_{m,2j}\varepsilon^{2j}\bigg)I_{3}
-\bigg(\varepsilon \sum^{l-1}_{j=0}(-1)^{j}c_{m,2j+1}\varepsilon^{2j+1}\bigg)
A(\mbox{\mathversion{bold}$v$}_{0}).
\end{align*}
Finally, after setting \( \varepsilon _{\ast}:= \min \{\varepsilon _{\ast,1},
\varepsilon _{\ast,2}\} \), for any \( \varepsilon 
\in (0,\varepsilon _{\ast}] \), we define \( \mbox{\mathversion{bold}$h$}^{\varepsilon}\)
on \( I\) as
\begin{eqnarray*}
\mbox{\mathversion{bold}$h$}^{\varepsilon}(s)=
\psi_{0}(s)\bigg(
\sum^{m}_{j=0}\frac{\partial^{2j}_{s}\mbox{\mathversion{bold}$h$}^{\varepsilon}(0)}{
(2j)!}s^{2j}\bigg)
+
\psi_{1}(s)\bigg(
\sum^{m}_{j=0}\frac{\partial^{2j}_{s}\mbox{\mathversion{bold}$h$}^{\varepsilon}(1)}{
(2j)!}(s-1)^{2j}\bigg),
\end{eqnarray*}
where \( \psi_{0} \) and \( \psi_{1}\) are smooth cut-off functions satisfying
\begin{eqnarray*}
\begin{array}{ll}
\psi_{0}(s) = 1, \ \psi_{1}(s)=0, & \mbox{for} \ s\in [0,\frac{1}{3}], \\[3mm]
\psi_{0}(s) = 0, \ \psi_{1}(s)=1, & \mbox{for} \ s\in [\frac{2}{3},1].
\end{array}
\end{eqnarray*}
Since the differential coefficients \( \partial^{2j}_{s}\mbox{\mathversion{bold}$h$}^{\varepsilon}|_{s=0,1} \
(0 \leq j \leq m)\)
were chosen to be \( O(\varepsilon )\), \( \mbox{\mathversion{bold}$h$}^{\varepsilon} \to 
\mbox{\mathversion{bold}$0$}\) as \( \varepsilon \to +0 \) in \( H^{2m+1}(I) \) and as a consequence,
\( \mbox{\mathversion{bold}$v$}^{\varepsilon}_{0} \to 
 \mbox{\mathversion{bold}$v$}_{0} \) in \( H^{2m+1}(I) \).
When \( \mbox{\mathversion{bold}$v$}_{0}\in H^{2m}(I)\), we can regard 
\( \mbox{\mathversion{bold}$v$}_{0}\) as an initial datum belonging to 
\( H^{2m-1}(I)\) and apply our arguements above and from the 
explicit form of \( \mbox{\mathversion{bold}$h$}^{\varepsilon }\), 
we see that 
\( \mbox{\mathversion{bold}$v$}_{0}^{\varepsilon } \to 
\mbox{\mathversion{bold}$v$}_{0}\) in \( H^{2m}(I)\).

Furthermore, by a similar calculation and also utilizing the arguments given in 
Rauch and Massey \cite{1}, we can prove that for any \( N\geq 1 \), we can construct a 
smooth approximating series \( \{ \mbox{\mathversion{bold}$v$}^{n}_{0}\}_{n=1}^{\infty} \)
such that \( \mbox{\mathversion{bold}$v$}^{n}_{0}\in 
H^{2(m+N)+1}(I) \ (n\in \mathbf{N})\), \( \mbox{\mathversion{bold}$v$}^{n}_{0}\) satisfies the compatibility conditions for 
(\ref{vslant}) up to order \( m+N \), and 
\( \mbox{\mathversion{bold}$v$}^{n}_{0}\to \mbox{\mathversion{bold}$v$}_{0}\) in 
\( H^{2m+1}(I) \) as \( n\to \infty \). The same is true for 
initial datum belonging to a Soblev space with even indices.

\medskip

We summarize the conclusions of this subsection in the following propositions.
\begin{pr}
Let \( l\) be an arbitrary non-negative integer. 
For \( \mbox{\mathversion{bold}$v$}_{0}\in H^{l+3}(I) \) 
satisfying the compatibility conditions for {\rm (\ref{vslant})} up to order 
\( [\frac{l+3}{2}]\),
there exists \( \varepsilon_{\ast}>0 \) such that for any \( \varepsilon \in (0,\varepsilon_{\ast}] \),
there exists a corrected initial datum \( \mbox{\mathversion{bold}$v$}^{\varepsilon}_{0}\) such that 
\( |\mbox{\mathversion{bold}$v$}^{\varepsilon}_{0}|\equiv 1\), 
\( \mbox{\mathversion{bold}$v$}^{\varepsilon}_{0}\in H^{l+3}(I)\), 
\( \mbox{\mathversion{bold}$v$}^{\varepsilon}_{0}\) satisfies the compatibility conditions for 
{\rm (\ref{rnl})} up to order \( [\frac{l+3}{2}]\), and
\begin{eqnarray*}
\mbox{\mathversion{bold}$v$}^{\varepsilon}_{0}\to \mbox{\mathversion{bold}$v$}_{0} \ \mbox{in} \ H^{l+3}(I)
\end{eqnarray*}
as \( \varepsilon \to +0\). \( \mbox{\mathversion{bold}$v$}^{\varepsilon}_{0}\) also satisfies
\begin{align*}
\|\mbox{\mathversion{bold}$v$}^{\varepsilon}_{0}\|_{l+3}\leq 
C\| \mbox{\mathversion{bold}$v$}_{0}\|_{l+3},
\end{align*}
where \( C>0\) is independent of \( \varepsilon \), which follows 
from the convergence of
\( \{ \mbox{\mathversion{bold}$v$}^{\varepsilon}_{0}\}_{0< \varepsilon \leq
\varepsilon_{\ast}}\).
\label{id1}
\end{pr}

\begin{pr}
Let \( l\) and \( N \) be non-negative integers. 
For \( \mbox{\mathversion{bold}$v$}_{0}\in H^{l+3}(I) \)
satisfying the compatibility conditions for {\rm (\ref{vslant})} up to order 
\( [\frac{l+3}{2}]\), 
there exists \( \{\mbox{\mathversion{bold}$v$}^{n}_{0}\}^{\infty}_{n=1} \) such that for any \( n\in \mathbf{N}\),
\( |\mbox{\mathversion{bold}$v$}^{n}_{0}|\equiv 1 \),
\( \mbox{\mathversion{bold}$v$}^{n}_{0}\in H^{l+3+N}(I) \), 
\( \mbox{\mathversion{bold}$v$}^{n}_{0}\)
satisfies the compatibility conditions for {\rm (\ref{vslant})} up to order 
\( [\frac{l+3+N}{2}] \), and 
\begin{eqnarray*}
\mbox{\mathversion{bold}$v$}^{n}_{0}\to \mbox{\mathversion{bold}$v$}_{0} \ \mbox{in} \ 
H^{l+3}(I)
\end{eqnarray*}
as \( n\to \infty \).
\label{id2}
\end{pr}
Hence, by combining the above two propositions, we see that given a \( \mbox{\mathversion{bold}$v$}_{0}\)
satisfying the compatibility conditions for (\ref{vslant}), we can construct a smoother
initial datum satisfying the necessary compatibility conditions for (\ref{rnl}).


\section{Construction of the solution to (\ref{rnl})}
\setcounter{equation}{0}

We construct the solution of (\ref{rnl}) based on an iteration scheme for the following 
linearized problem.
\begin{eqnarray}
\left\{
\begin{array}{ll}
\mbox{\mathversion{bold}$u$}_{t}=\varepsilon \mbox{\mathversion{bold}$u$}_{ss} + 
\mbox{\mathversion{bold}$b$}\times \mbox{\mathversion{bold}$u$}_{ss} + \mbox{\mathversion{bold}$f$}, & s\in I, \ t>0, \\[3mm]
\mbox{\mathversion{bold}$u$}(s,0)=\mbox{\mathversion{bold}$u$}_{0}(s), & s\in I, \\[3mm]
\mbox{\mathversion{bold}$u$}(0,t)=\mbox{\mathversion{bold}$a$}, \ 
\mbox{\mathversion{bold}$u$}(1,t) = \mbox{\mathversion{bold}$e$}_{3}, & t>0,
\end{array}\right.
\label{linear}
\end{eqnarray}
where \( \mbox{\mathversion{bold}$u$}_{0}= \mbox{\mathversion{bold}$u$}_{0}(s)\), 
\( \mbox{\mathversion{bold}$b$}=\mbox{\mathversion{bold}$b$}(s,t) \), and 
\( \mbox{\mathversion{bold}$f$}=\mbox{\mathversion{bold}$f$}(s,t)\) are given vector valued functions.
We make some notations to define the compatibility conditions for 
(\ref{linear}).
Set \( \mbox{\mathversion{bold}$L$}_{1}(\mbox{\mathversion{bold}$u$},\mbox{\mathversion{bold}$b$},
\mbox{\mathversion{bold}$f$}):= 
\varepsilon \mbox{\mathversion{bold}$u$}_{ss} + 
\mbox{\mathversion{bold}$b$}\times \mbox{\mathversion{bold}$u$}_{ss} + \mbox{\mathversion{bold}$f$}\)
and 
\begin{eqnarray*}
\mbox{\mathversion{bold}$L$}_{m}
:=
\varepsilon \partial^{2}_{s}\mbox{\mathversion{bold}$L$}_{m-1}+
\sum^{m-1}_{j=0}
\left(
\begin{array}{c}
m-1\\
j
\end{array}
\right)
\partial^{j}_{t}\mbox{\mathversion{bold}$b$}\times \partial^{2}_{s}\mbox{\mathversion{bold}$L$}_{m-1-j}
+
\partial^{m-1}_{t}\mbox{\mathversion{bold}$f$},
\end{eqnarray*}
for \( m \geq 2 \). Here, \( \mbox{\mathversion{bold}$L$}_{k}=\mbox{\mathversion{bold}$L$}_{k}
(\mbox{\mathversion{bold}$u$},\mbox{\mathversion{bold}$b$},
\mbox{\mathversion{bold}$f$})\) for \( 1\leq k\leq m \).
The compatibility conditions are defined as follows.
\begin{df}(Compatibility conditions for (\ref{linear})).
For a non-negative integer \( m\),, 
we say that \( \mbox{\mathversion{bold}$u$}_{0}\in H^{2m+1}(I)\), \( \mbox{\mathversion{bold}$b$}\in Y_{T}^{m}(I)\), and
\( \mbox{\mathversion{bold}$f$}\in Z_{T}^{m}(I) \) satisfy the \( m\)-th order compatibility condition
for {\rm (\ref{linear})} if
\begin{eqnarray*}
\mbox{\mathversion{bold}$u$}_{0}(0)=\mbox{\mathversion{bold}$a$}, \
\mbox{\mathversion{bold}$u$}_{0}(1)=\mbox{\mathversion{bold}$e$}_{3}
\end{eqnarray*}
when \( m=0\), and
\begin{eqnarray*}
L_{m}(\mbox{\mathversion{bold}$u$}_{0},\mbox{\mathversion{bold}$b$},
\mbox{\mathversion{bold}$f$})|_{s=0,1} = \mbox{\mathversion{bold}$0$}
\end{eqnarray*}
when \( m\geq 1 \).
We also say that \( \mbox{\mathversion{bold}$u$}_{0}\),\(\mbox{\mathversion{bold}$b$}\), and 
\(\mbox{\mathversion{bold}$f$}\) satisfy the compatibility conditions for {\rm (\ref{linear})} up to order 
\( m\) if they satisfy the \( k\)-th order compatibility condition for all \( k \) with \( 0\leq k \leq m\).
\end{df}
Here, \( Y^{m}_{T}(I) \) and \( Z^{m}_{T}(I) \) are function spaces defined as follows.
\begin{align*}
Y^{m}_{T}(I) &:= \bigcap ^{m}_{j=0}C^{j}\big( [0,T];H^{2(m-j)+1}(I)\big)\\[3mm]
Z^{m}_{T}(I) &:= 
\bigcap ^{m}_{j=0}C^{j}\big([0,T]; H^{2(m-j)}(I)\big).
\end{align*}

Fix \( m\geq 1 \).
The solution of (\ref{rnl}) is constructed by the following iteration scheme. We define
\( \mbox{\mathversion{bold}$u$}^{n}\) as the solution of
\begin{eqnarray}
\left\{
\begin{array}{ll}
\mbox{\mathversion{bold}$u$}^{n}_{t}=\varepsilon \mbox{\mathversion{bold}$u$}_{ss}
+\mbox{\mathversion{bold}$u$}^{n-1}\times \mbox{\mathversion{bold}$u$}^{n}_{ss}
+ \varepsilon |\mbox{\mathversion{bold}$u$}^{n-1}_{s}|^{2}\mbox{\mathversion{bold}$u$}^{n-1}, & s\in I, \ t>0, \\[3mm]
\mbox{\mathversion{bold}$u$}^{n}(s,0) = \mbox{\mathversion{bold}$v$}^{\varepsilon}_{0}(s), & s\in I, \\[3mm]
\mbox{\mathversion{bold}$u$}^{n}(0,t)=\mbox{\mathversion{bold}$a$}, \ 
\mbox{\mathversion{bold}$u$}^{n}(1,t)=\mbox{\mathversion{bold}$e$}_{3}, & t>0,
\end{array}\right.
\label{iter}
\end{eqnarray}
for \( n \geq 2 \). Here, \( \mbox{\mathversion{bold}$v$}^{\varepsilon}_{0}\) is the corrected initial datum obtained in 
Section 3.
Now, we must define 
\( \mbox{\mathversion{bold}$u$}^{1} \) appropriately so that the necessary compatibility conditions are satisfied at 
each step of the iteration. This is accomplished by choosing
\begin{eqnarray*}
\mbox{\mathversion{bold}$u$}^{1}(s,t) =
\mbox{\mathversion{bold}$v$}^{\varepsilon}_{0}(s) +
\sum^{m}_{j=1}\frac{t^{j}}{j!}\mbox{\mathversion{bold}$Q$}_{j}(\mbox{\mathversion{bold}$v$}^{\varepsilon}_{0}),
\end{eqnarray*}
where \( \mbox{\mathversion{bold}$Q$}_{j}\) was defined in Section 3.
From Proposition \ref{id1} and \ref{id2}, we can assume that \( \mbox{\mathversion{bold}$v$}^{\varepsilon}_{0}\)
is smooth, and thus we assume that \( \mbox{\mathversion{bold}$v$}^{\varepsilon}_{0}\) is smooth enough that
\( \mbox{\mathversion{bold}$u$}^{1}\in Y^{m+1}_{T}(I) \).
The solvability of (\ref{linear}), and as a consequence, the fact that \( \{\mbox{\mathversion{bold}$u$}^{n}\}^{\infty}_{n=1}\)
is well-defined, can be proved by utilizing the
Sobolev--Slobodetskii space \( W^{2m+1+\alpha , m+(1+\alpha)/2}_{2}(I\times [0,T])\) 
with \( \alpha \in (0,\frac{1}{2}) \) following the arguments in Section 2 of 
Nishiyama \cite{10}. 
The convergence of \( \{ \mbox{\mathversion{bold}$u$}^{n}\}^{\infty}_{n=1} \) can also be proved 
in the same Sobolev--Slobodetskii space, based on the estimate of the solution for (\ref{linear}).
Also see 
Solonnikov \cite{2} for the definition and properties of 
Sobolev--Slobodetskii spaces.

The limit of \( \{ \mbox{\mathversion{bold}$u$}^{n}\}^{\infty}_{n=1}\) is the desired solution of (\ref{rnl}) and we arrive at 
the following existence theorem for (\ref{rnl}).
\begin{pr}
Let \( m\geq 1 \) be an integer and \( T>0 \). There exists a unique solution 
\( \mbox{\mathversion{bold}$v$}^{\varepsilon}\in Y^{m}_{T}(I)\) to {\rm (\ref{rnl})} with a smooth initial datum
\( \mbox{\mathversion{bold}$v$}^{\varepsilon}_{0}\).
\end{pr}


\section{Unique solvability of problem (\ref{vslant})}

We construct the solution of (\ref{vslant}) by taking the limit \( \varepsilon \to +0 \) in (\ref{rnl}). 
We first derive estimates for the solution of (\ref{rnl}) uniform in \( \varepsilon \), and then prove the convergence of 
the solution as \( \varepsilon \to +0 \).
In this section and for the rest of the paper, \( C\) denotes generic positive constants which may be different from 
line to line. What \( C\) depends on will be stated when ever it is necessary.

\subsection{Uniform estimate of \( \mbox{\mathversion{bold}$v$}^{\varepsilon}\) with respect to \( \varepsilon \)}

We first prove the following.
\begin{lm}
For any \( m\geq 1 \), a solution \( \mbox{\mathversion{bold}$v$}^{\varepsilon}\in Y^{m}_{T}(I)\) of {\rm (\ref{rnl})}
satisfies \( |\mbox{\mathversion{bold}$v$}^{\varepsilon}|= 1 \) in \( I\times [0,T]\).
\label{nstr}
\end{lm}
{\it Proof.} 
Setting \( h(s,t):= |\mbox{\mathversion{bold}$v$}^{\varepsilon}(s,t)|^{2}-1 \), we see that \( h\) satisfies
\begin{eqnarray*}
\left\{
\begin{array}{ll}
h_{t}= \varepsilon h_{ss} + 2\varepsilon |\mbox{\mathversion{bold}$v$}^{\varepsilon}_{s}|^{2}h, & s\in I, \ t>0, \\[3mm]
h(s,0)=0, & s\in I, \\[3mm]
h(0,t)=h(1,t)=0, & t>0,
\end{array}\right.
\end{eqnarray*}
where \( |\mbox{\mathversion{bold}$a$}|=|\mbox{\mathversion{bold}$e$}_{3}|=1 \) was used.
From the Sobolev's embedding theorem,
\begin{align*}
\frac{1}{2}\frac{{\rm d}}{{\rm d}t}\| h\| ^{2}=
\varepsilon (h,h_{ss}) + 2\varepsilon (h, |\mbox{\mathversion{bold}$v$}^{\varepsilon}_{s}|^{2}h)
\leq 
-\varepsilon \|h_{s}\|^{2}+C\|h\|^{2}
\end{align*}
holds. Here, \( C>0\) depends on \( 
\| \mbox{\mathversion{bold}$v$}_{s}^{\varepsilon}\|_{1}\), and from Gronwall's inequality 
we see that \( h=0 \) in \( I\times [0,T] \), and this finishes the proof. \hfill \( \Box .\)

\medskip
We introduce a property which will be used throughout this section.
When \( \mbox{\mathversion{bold}$v$}^{\varepsilon}_{s}\neq \mbox{\mathversion{bold}$0$}\), 
\( \mbox{\mathversion{bold}$v$}^{\varepsilon}\),
\( \frac{\mbox{\mathversion{bold}$v$}^{\varepsilon}_{s}}{|\mbox{\mathversion{bold}$v$}^{\varepsilon}_{s}|}\),
and 
\( \frac{\mbox{\mathversion{bold}$v$}^{\varepsilon}\times \mbox{\mathversion{bold}$v$}^{\varepsilon}_{s}}
{|\mbox{\mathversion{bold}$v$}^{\varepsilon}_{s}|}\) form a orthonormal basis of 
\( \mathbf{R}^{3}\). 
From this we have for \( n\geq 2\),
\begin{eqnarray}
\mbox{\mathversion{bold}$v$}^{\varepsilon}_{s}\times \partial^{n}_{s}\mbox{\mathversion{bold}$v$}^{\varepsilon}
=
-[\mbox{\mathversion{bold}$v$}^{\varepsilon}\cdot \partial^{n}_{s}\mbox{\mathversion{bold}$v$}^{\varepsilon}]
\mbox{\mathversion{bold}$v$}^{\varepsilon}\times 
\mbox{\mathversion{bold}$v$}^{\varepsilon}_{s}
+
[(\mbox{\mathversion{bold}$v$}^{\varepsilon}\times \mbox{\mathversion{bold}$v$}^{\varepsilon}_{s})\cdot
\partial^{n}_{s}\mbox{\mathversion{bold}$v$}^{\varepsilon}]\mbox{\mathversion{bold}$v$}^{\varepsilon}.
\label{decom}
\end{eqnarray}
The above equation is also true when \( \mbox{\mathversion{bold}$v$}^{\varepsilon}_{s}=
\mbox{\mathversion{bold}$0$}\). We also have from Lemma \ref{nstr}, 
\begin{eqnarray}
\mbox{\mathversion{bold}$v$}^{\varepsilon}\cdot \partial ^{n}_{s}\mbox{\mathversion{bold}$v$}^{\varepsilon}
=
-\frac{1}{2}\sum^{n-1}_{j=1}
\left(
\begin{array}{c}
n\\
j
\end{array}\right)
\partial^{j}_{s}\mbox{\mathversion{bold}$v$}^{\varepsilon}\cdot
\partial^{n-j}_{s}\mbox{\mathversion{bold}$v$}^{\varepsilon}
\label{one}
\end{eqnarray}
for \( n\geq 2\), and \( \mbox{\mathversion{bold}$v$}^{\varepsilon}\cdot \mbox{\mathversion{bold}$v$}^{\varepsilon}_{s}
= 0\).
We further introduce two auxilary lemmas which we will use to 
prove the uniform estimate.
\begin{lm}
For \( m\geq 2 \), we have
\begin{align}
\partial^{m}_{t}\mbox{\mathversion{bold}$v$}^{\varepsilon}
&=
\sum^{m}_{j=0}a_{m,j}\varepsilon^{j}
A_{m-j}\partial^{2m}_{s}\mbox{\mathversion{bold}$v$}^{\varepsilon}
+
\sum^{m-1}_{j=0}
\sum^{m-1-j}_{k=0}
e_{m,j,k}\varepsilon^{j}
A_{m-1-j-k}\mbox{\mathversion{bold}$v$}^{\varepsilon}_{s}\times 
\big( A_{k}\partial^{2m-1}_{s}\mbox{\mathversion{bold}$v$}^{\varepsilon}\big)
\label{2m2m-1}\\[3mm]
& \qquad +
\sum^{m}_{j=1}b_{m,j}\varepsilon^{m+1-j}
[\mbox{\mathversion{bold}$v$}^{\varepsilon}_{s}\cdot 
(A_{j-1}\partial^{2m-1}_{s}\mbox{\mathversion{bold}$v$}^{\varepsilon})]
\mbox{\mathversion{bold}$v$}^{\varepsilon}
+
\mbox{\mathversion{bold}$U$}^{\varepsilon}_{m}, \nonumber 
\end{align}
where \( \mbox{\mathversion{bold}$U$}^{\varepsilon}_{m} \) are terms 
that can be estimated as
\begin{align*}
\big| \mbox{\mathversion{bold}$U$}^{\varepsilon}_{m} |_{s=0,1} \big|
\leq 
C\|\mbox{\mathversion{bold}$v$}^{\varepsilon}_{s}\|_{2(m-1)},
\end{align*}
with \( C>0\) depending on
\( \| \mbox{\mathversion{bold}$v$}^{\varepsilon}_{s}\|_{2(m-2)}\), and 
\( A_{k} = A_{k}(\mbox{\mathversion{bold}$v$}^{\varepsilon})\) for 
\( 0\leq k\leq m\) is defined in Lemma {\rm \ref{fhighlm}}.
\( a_{m,j}\), \(e_{m,j,k}\), and \( b_{m,j}\) are absolute constants 
independent of \( \varepsilon \), and \( a_{m,0}=1 \).
\label{lm2m2m-1}
\end{lm}
\begin{lm}
For \( m\geq 2\) we have
\begin{align}
\mbox{\mathversion{bold}$v$}^{\varepsilon}\times 
\partial^{2m}_{s}\mbox{\mathversion{bold}$v$}^{\varepsilon}|_{s=0,1}
&=
(-1)^{m}\varepsilon A_{m-1}
\big\{ \partial^{m-1}_{t}\mbox{\mathversion{bold}$v$}^{\varepsilon}_{ss}
+
2(\mbox{\mathversion{bold}$v$}^{\varepsilon}_{s}\cdot 
\partial^{m-1}_{t}\mbox{\mathversion{bold}$v$}^{\varepsilon}_{s})
\mbox{\mathversion{bold}$v$}^{\varepsilon}\big\}
\label{vc2m}\\[3mm]
& \qquad +
\varepsilon 
\bigg\{
\sum^{m-2}_{l=0}(-1)^{l}
A_{l+2}\partial^{l}_{t}
\partial^{2(m-l)}_{s}\mbox{\mathversion{bold}$v$}^{\varepsilon}
\bigg\}
+
\varepsilon \mbox{\mathversion{bold}$W$}^{\varepsilon}_{m}\big|_{s=0,1},
\nonumber
\end{align}
\begin{align}
\partial^{j}_{s}\mbox{\mathversion{bold}$v$}^{\varepsilon}\cdot 
\partial^{k}_{s}\mbox{\mathversion{bold}$v$}^{\varepsilon}|_{s=0,1}
&=
\varepsilon 
Y_{j,k}|_{s=0,1} \quad \mbox{for} \ j+k = 2(m-1)+1.
\label{oddinner}
\end{align}
where \( \mbox{\mathversion{bold}$W$}^{\varepsilon}_{m}\) are terms 
satisfying the estimate
\begin{align*}
\big|
\mbox{\mathversion{bold}$W$}^{\varepsilon}_{m}|_{s=0,1}
\big|
\leq 
C\| \mbox{\mathversion{bold}$v$}^{\varepsilon}_{s}\|_{2(m-1)}
\end{align*}
with \( C>0\) depending on
\( \| \mbox{\mathversion{bold}$v$}^{\varepsilon}_{s}\|_{2(m-2)}\), and
\( Y_{j,k}\) are terms satisfying the estimate
\begin{align*}
\big| Y_{j,k}|_{s=0,1}\big|
\leq
C\|\mbox{\mathversion{bold}$v$}^{\varepsilon}_{s}\|^{2}_{2(m-1)},
\end{align*}
when \( j+k = 2(m-1)+1 \).
Note that from the equation in 
{\rm (\ref{rnl})}, the right-hand side of {\rm (\ref{vc2m})} can be estimated from above by
\( \varepsilon C\| \mbox{\mathversion{bold}$v$}^{\varepsilon}_{s}\|_{2m}\), where
\( C>0\) depends on 
\( \| \mbox{\mathversion{bold}$v$}^{\varepsilon}_{s}\|_{2(m-1)}\).
\label{lmvc2m}
\end{lm}
Lemma \ref{lm2m2m-1} can be proved by induction
in a similar way as Lemma \ref{fhighlm}
and Lemma \ref{rmhigh} and hence the details are omitted.\\[3mm]
\noindent
{\it Proof of Lemma \ref{lmvc2m}} \ We prove (\ref{vc2m}) and (\ref{oddinner}) simultaneously by induction.
From \( \mbox{\mathversion{bold}$v$}^{\varepsilon}_{t}|_{s=0,1}=
\mbox{\mathversion{bold}$0$}\), we have from the equation
\begin{align*}
\mbox{\mathversion{bold}$v$}^{\varepsilon}\times 
\mbox{\mathversion{bold}$v$}^{\varepsilon}_{ss} +
\varepsilon \mbox{\mathversion{bold}$v$}^{\varepsilon}_{ss}
+\varepsilon |\mbox{\mathversion{bold}$v$}^{\varepsilon}_{s}|^{2}
\mbox{\mathversion{bold}$v$}^{\varepsilon}|_{s=0,1}
= \mbox{\mathversion{bold}$0$}.
\end{align*}
Acting the exterior product \( \mbox{\mathversion{bold}$v$}^{\varepsilon}\times \)
from the left-side in the above equation along with (\ref{one}) yields
\begin{align}
\mbox{\mathversion{bold}$v$}^{\varepsilon}_{ss}|_{s=0,1}
=
-\varepsilon \mbox{\mathversion{bold}$v$}^{\varepsilon}\times
\mbox{\mathversion{bold}$v$}^{\varepsilon}_{ss}
-
|\mbox{\mathversion{bold}$v$}^{\varepsilon}_{s}|^{2}
\mbox{\mathversion{bold}$v$}^{\varepsilon}
|_{s=0,1}.
\label{bv2}
\end{align}
Hence we have
\begin{align*}
\mbox{\mathversion{bold}$v$}^{\varepsilon}_{s}\cdot 
\mbox{\mathversion{bold}$v$}^{\varepsilon}_{ss}|_{s=0,1}
=
-\varepsilon 
(\mbox{\mathversion{bold}$v$}^{\varepsilon}_{s}\cdot 
(\mbox{\mathversion{bold}$v$}^{\varepsilon}\times 
\mbox{\mathversion{bold}$v$}^{\varepsilon}_{ss}))|_{s=0,1}
\end{align*}
and from (\ref{one}) we have
\begin{align*}
\mbox{\mathversion{bold}$v$}^{\varepsilon}\cdot 
\mbox{\mathversion{bold}$v$}^{\varepsilon}_{sss}|_{s=0,1}
=
-3\mbox{\mathversion{bold}$v$}^{\varepsilon}_{s}\cdot 
\mbox{\mathversion{bold}$v$}^{\varepsilon}_{ss}|_{s=0,1}
=
3\varepsilon 
(\mbox{\mathversion{bold}$v$}^{\varepsilon}_{s}\cdot 
(\mbox{\mathversion{bold}$v$}^{\varepsilon}\times 
\mbox{\mathversion{bold}$v$}^{\varepsilon}_{ss}))|_{s=0,1},
\end{align*}
which proves (\ref{oddinner}) for \( m=2\).
From the equation in (\ref{rnl}) and 
\( \mbox{\mathversion{bold}$v$}^{\varepsilon}_{t}|_{s=0,1} =
\mbox{\mathversion{bold}$0$}\), we have
\begin{align*}
\mbox{\mathversion{bold}$v$}^{\varepsilon}\times 
\mbox{\mathversion{bold}$v$}^{\varepsilon}_{ss}|_{s=0,1}
=
-\varepsilon( \mbox{\mathversion{bold}$v$}^{\varepsilon}_{ss}+
|\mbox{\mathversion{bold}$v$}^{\varepsilon}_{s}|^{2}
\mbox{\mathversion{bold}$v$}^{\varepsilon})|_{s=0,1}.
\end{align*}
Taking the \( t\) derivative of the above yields
\begin{align*}
\mbox{\mathversion{bold}$v$}^{\varepsilon}\times 
\mbox{\mathversion{bold}$v$}^{\varepsilon}_{tss}|_{s=0,1}
&=
-\varepsilon 
\{
\mbox{\mathversion{bold}$v$}^{\varepsilon}_{tss}+
2(\mbox{\mathversion{bold}$v$}^{\varepsilon}_{s}\cdot 
\mbox{\mathversion{bold}$v$}^{\varepsilon}_{ts})
\mbox{\mathversion{bold}$v$}^{\varepsilon} \}|_{s=0,1}.
\end{align*}
On the other hand, by substituting the equation for 
\( \mbox{\mathversion{bold}$v$}^{\varepsilon}\), we have
\begin{align*}
\mbox{\mathversion{bold}$v$}^{\varepsilon}\times 
\mbox{\mathversion{bold}$v$}^{\varepsilon}_{tss}|_{s=0,1}
&=
\mbox{\mathversion{bold}$v$}^{\varepsilon}\times 
(\mbox{\mathversion{bold}$v$}^{\varepsilon}\times 
\partial^{4}_{s}\mbox{\mathversion{bold}$v$}^{\varepsilon})
+
2\mbox{\mathversion{bold}$v$}^{\varepsilon}\times 
(\mbox{\mathversion{bold}$v$}^{\varepsilon}_{s}\times 
\mbox{\mathversion{bold}$v$}^{\varepsilon}_{sss})
+
\varepsilon \mbox{\mathversion{bold}$v$}^{\varepsilon}\times 
\partial^{4}_{s}\mbox{\mathversion{bold}$v$}^{\varepsilon}
\\[3mm]
& \qquad +
2\varepsilon
(\mbox{\mathversion{bold}$v$}^{\varepsilon}_{s}\cdot 
\mbox{\mathversion{bold}$v$}^{\varepsilon}_{sss})
\mbox{\mathversion{bold}$v$}^{\varepsilon}
+
\varepsilon 
\mbox{\mathversion{bold}$W$}^{\varepsilon}_{2}|_{s=0,1}\\[3mm]
&=
(\mbox{\mathversion{bold}$v$}^{\varepsilon}\cdot 
\partial^{4}_{s}\mbox{\mathversion{bold}$v$}^{\varepsilon})
\mbox{\mathversion{bold}$v$}^{\varepsilon}
-
\partial^{4}_{s}\mbox{\mathversion{bold}$v$}^{\varepsilon}
+
\varepsilon
\mbox{\mathversion{bold}$v$}^{\varepsilon}\times 
\partial^{4}_{s}\mbox{\mathversion{bold}$v$}^{\varepsilon}
+
2\varepsilon (\mbox{\mathversion{bold}$v$}^{\varepsilon}_{s}\cdot 
\mbox{\mathversion{bold}$v$}^{\varepsilon}_{sss})
\mbox{\mathversion{bold}$v$}^{\varepsilon}\\[3mm]
& \qquad +
2 \varepsilon Y_{0,3}\mbox{\mathversion{bold}$v$}^{\varepsilon}_{s}
+
\varepsilon \mbox{\mathversion{bold}$W$}^{\varepsilon}_{2}|_{s=0,1}.
\end{align*}
Since the exact form of 
\( \mbox{\mathversion{bold}$W$}^{\varepsilon}_{m}\) is not
needed for the proof and the application of the lemma, 
\( \mbox{\mathversion{bold}$W$}^{\varepsilon}_{2}\) will denote 
the collection of terms satisfying the property stated in the lemma, and 
may change from line to line. For example, 
the term \( 2 \varepsilon Y_{0,3}\mbox{\mathversion{bold}$v$}^{\varepsilon}_{s}\)
can be included in 
\( \varepsilon \mbox{\mathversion{bold}$W$}^{\varepsilon}_{2}\).
Combining the above two equations yields
\begin{align*}
(\mbox{\mathversion{bold}$v$}^{\varepsilon}\cdot 
\partial^{4}_{s}\mbox{\mathversion{bold}$v$}^{\varepsilon})
\mbox{\mathversion{bold}$v$}^{\varepsilon}
&-
\partial^{4}_{s}\mbox{\mathversion{bold}$v$}^{\varepsilon}
+
\varepsilon
\mbox{\mathversion{bold}$v$}^{\varepsilon}\times 
\partial^{4}_{s}\mbox{\mathversion{bold}$v$}^{\varepsilon}
+
2\varepsilon (\mbox{\mathversion{bold}$v$}^{\varepsilon}_{s}\cdot 
\mbox{\mathversion{bold}$v$}^{\varepsilon}_{sss})
\mbox{\mathversion{bold}$v$}^{\varepsilon}
+
\varepsilon \mbox{\mathversion{bold}$W$}^{\varepsilon}_{2}|_{s=0,1}\\[3mm]
&=
-\varepsilon 
\{
\mbox{\mathversion{bold}$v$}^{\varepsilon}_{tss}+
2(\mbox{\mathversion{bold}$v$}^{\varepsilon}_{s}\cdot 
\mbox{\mathversion{bold}$v$}^{\varepsilon}_{ts})
\mbox{\mathversion{bold}$v$}^{\varepsilon} \}|_{s=0,1}.
\end{align*}
Taking the exterior product with \( \mbox{\mathversion{bold}$v$}^{\varepsilon}\)
from the left-side and rearranging the terms yield
\begin{align*}
\mbox{\mathversion{bold}$v$}^{\varepsilon}\times 
\partial^{4}_{s}\mbox{\mathversion{bold}$v$}^{\varepsilon}|_{s=0,1}
&=
\varepsilon 
\mbox{\mathversion{bold}$v$}^{\varepsilon}\times 
\{
\mbox{\mathversion{bold}$v$}^{\varepsilon}_{tss}+
2(\mbox{\mathversion{bold}$v$}^{\varepsilon}_{s}\cdot 
\mbox{\mathversion{bold}$v$}^{\varepsilon}_{ts})
\mbox{\mathversion{bold}$v$}^{\varepsilon} \}
+
\varepsilon \mbox{\mathversion{bold}$v$}^{\varepsilon}\times 
(\mbox{\mathversion{bold}$v$}^{\varepsilon}\times 
\partial^{4}_{s}\mbox{\mathversion{bold}$v$}^{\varepsilon})
+
\varepsilon \mbox{\mathversion{bold}$W$}^{\varepsilon}_{2}|_{s=0,1}\\[3mm]
&=
\varepsilon 
A_{1} 
\{
\mbox{\mathversion{bold}$v$}^{\varepsilon}_{tss}+
2(\mbox{\mathversion{bold}$v$}^{\varepsilon}_{s}\cdot 
\mbox{\mathversion{bold}$v$}^{\varepsilon}_{ts})
\mbox{\mathversion{bold}$v$}^{\varepsilon} \}
+
\varepsilon \mbox{\mathversion{bold}$v$}^{\varepsilon}\times 
(\mbox{\mathversion{bold}$v$}^{\varepsilon}\times 
\partial^{4}_{s}\mbox{\mathversion{bold}$v$}^{\varepsilon})
+
\varepsilon \mbox{\mathversion{bold}$W$}^{\varepsilon}_{2}|_{s=0,1},
\end{align*}
which proves (\ref{vc2m}) with \( m=2 \).
Suppose that (\ref{vc2m}) and (\ref{oddinner}) hold up to 
\( m-1\) for some \( m \geq 3 \).
From the assumption of induction, we have
\begin{align*}
\mbox{\mathversion{bold}$v$}^{\varepsilon}\times 
\partial^{2(m-1)}_{s}\mbox{\mathversion{bold}$v$}^{\varepsilon}|_{s=0,1}
&=
(-1)^{m-1}\varepsilon A_{m-2}\big\{
\partial^{m-2}_{t}\mbox{\mathversion{bold}$v$}^{\varepsilon}_{ss}
+
2(\mbox{\mathversion{bold}$v$}^{\varepsilon}_{s}\cdot 
\partial^{m-2}_{t}\mbox{\mathversion{bold}$v$}^{\varepsilon}_{s})
\mbox{\mathversion{bold}$v$}^{\varepsilon}
\big\}\\[3mm]
&+
\varepsilon 
\bigg\{
\sum^{m-3}_{l=0}(-1)^{l}
A_{l+2}
\partial^{l}_{t}\partial^{2(m-1-l)}_{s}
\mbox{\mathversion{bold}$v$}^{\varepsilon}\bigg\}
+
\varepsilon \mbox{\mathversion{bold}$W$}^{\varepsilon}_{m-1}\big| _{s=0,1}. 
\end{align*}
Taking the \( t\) derivative of the above yields
\begin{align}
\mbox{\mathversion{bold}$v$}^{\varepsilon}\times
\partial^{2(m-1)}_{s}\mbox{\mathversion{bold}$v$}^{\varepsilon}_{t}|_{s=0,1}
&=
(-1)^{m-1}\varepsilon 
A_{m-2}\big\{
\partial^{m-1}_{t}\mbox{\mathversion{bold}$v$}^{\varepsilon}_{ss}
+
2(\mbox{\mathversion{bold}$v$}^{\varepsilon}_{s}\cdot 
\partial^{m-1}_{t}\mbox{\mathversion{bold}$v$}^{\varepsilon}_{s})
\mbox{\mathversion{bold}$v$}^{\varepsilon}
\big\}\label{mt}\\[3mm]
& \qquad +
\varepsilon \bigg\{
\sum^{m-3}_{l=0}(-1)^{l}
A_{l+2}
\partial^{l+1}_{t}\partial^{2(m-1-l)}_{s}
\mbox{\mathversion{bold}$v$}^{\varepsilon}
\bigg\}
+
\varepsilon 
\mbox{\mathversion{bold}$W$}^{\varepsilon}_{m}\big|_{s=0,1},
\nonumber
\end{align}
where \( 
\mbox{\mathversion{bold}$W$}^{\varepsilon}_{m}\) 
denotes the collection of terms satisfying the condition of the lemma and
as before, the contents of \( 
\mbox{\mathversion{bold}$W$}^{\varepsilon}_{m}\) will change from 
line to line. First we observe
\begin{align*}
\partial^{2(m-1)}_{s}\mbox{\mathversion{bold}$v$}^{\varepsilon}_{t}
&=
\mbox{\mathversion{bold}$v$}^{\varepsilon}\times 
\partial^{2m}_{s}\mbox{\mathversion{bold}$v$}^{\varepsilon}
+
2(m-1)\mbox{\mathversion{bold}$v$}^{\varepsilon}_{s}\times 
\partial^{2m-1}_{s}\mbox{\mathversion{bold}$v$}^{\varepsilon}
+
\varepsilon \partial^{2m}_{s}\mbox{\mathversion{bold}$v$}^{\varepsilon}
+
2\varepsilon (\mbox{\mathversion{bold}$v$}^{\varepsilon}_{s}\cdot 
\partial^{2m-1}_{s}\mbox{\mathversion{bold}$v$}^{\varepsilon})
\mbox{\mathversion{bold}$v$}^{\varepsilon}\\[3mm]
& \qquad +
\sum^{2(m-1)-2}_{j=2}
\left(
\begin{array}{c}
2(m-1)\\
j
\end{array}\right)
\partial^{j}_{s}\mbox{\mathversion{bold}$v$}^{\varepsilon}
\times \partial^{2(m-1)-j+2}_{s}\mbox{\mathversion{bold}$v$}^{\varepsilon}
+
\varepsilon \mbox{\mathversion{bold}$W$}^{\varepsilon}_{m}. 
\end{align*}
Taking the exterior product with \( \mbox{\mathversion{bold}$v$}^{\varepsilon}\)
from the left yields
\begin{align}
\label{m-1}
\mbox{\mathversion{bold}$v$}^{\varepsilon}\times 
\partial^{2(m-1)}_{s}\mbox{\mathversion{bold}$v$}^{\varepsilon}_{t}
&=(
\mbox{\mathversion{bold}$v$}^{\varepsilon}\cdot 
\partial^{2m}_{s}\mbox{\mathversion{bold}$v$}^{\varepsilon})
\mbox{\mathversion{bold}$v$}^{\varepsilon}
-
\partial^{2m}_{s}\mbox{\mathversion{bold}$v$}^{\varepsilon}
+
2(m-1)(\mbox{\mathversion{bold}$v$}^{\varepsilon}\cdot 
\partial^{2m-1}_{s}\mbox{\mathversion{bold}$v$}^{\varepsilon})
\mbox{\mathversion{bold}$v$}^{\varepsilon}
+
\varepsilon \mbox{\mathversion{bold}$v$}^{\varepsilon}\times
\partial^{2m}_{s}\mbox{\mathversion{bold}$v$}^{\varepsilon}\\[3mm]
& \qquad 
+
\sum^{2(m-1)-2}_{j=2}
\left(
\begin{array}{c}
2(m-1)\\
j
\end{array}\right)
\mbox{\mathversion{bold}$v$}^{\varepsilon}\times 
(\partial^{j}_{s}\mbox{\mathversion{bold}$v$}^{\varepsilon}\times 
\partial^{2(m-1)-j+2}_{s}\mbox{\mathversion{bold}$v$}^{\varepsilon})
+
\varepsilon \mbox{\mathversion{bold}$W$}^{\varepsilon}_{m}. \nonumber
\end{align}
Since the terms in the summation can be calculated as
\begin{align*}
\mbox{\mathversion{bold}$v$}^{\varepsilon}\times 
(\partial^{j}_{s}\mbox{\mathversion{bold}$v$}^{\varepsilon}\times 
\partial^{2(m-1)-j+2}_{s}\mbox{\mathversion{bold}$v$}^{\varepsilon})
=
(\mbox{\mathversion{bold}$v$}^{\varepsilon}\cdot 
\partial^{2(m-1)-j+2}_{s}\mbox{\mathversion{bold}$v$}^{\varepsilon})
\mbox{\mathversion{bold}$v$}^{\varepsilon}
-
(\mbox{\mathversion{bold}$v$}^{\varepsilon}\cdot 
\partial^{j}_{s}\mbox{\mathversion{bold}$v$}^{\varepsilon})
\mbox{\mathversion{bold}$v$}^{\varepsilon},
\end{align*}
if \( j\) is odd, then \( 2(m-1)-j+2 \) is also odd and we have from the 
assumption of induction for (\ref{oddinner}) that
\begin{align*}
\mbox{\mathversion{bold}$v$}^{\varepsilon}\cdot 
\partial^{2(m-1)-j+2}_{s}\mbox{\mathversion{bold}$v$}^{\varepsilon}|_{s=0,1}
&=
\varepsilon Y_{0,2(m-1)-j+2}|_{s=0,1}\\[3mm]
\mbox{\mathversion{bold}$v$}^{\varepsilon}\cdot 
\partial^{j}_{s}\mbox{\mathversion{bold}$v$}^{\varepsilon}|_{s=0,1}
&=
\varepsilon Y_{0,j}|_{s=0,1}
\end{align*}
holds, which implies that the terms in the summation
of (\ref{m-1}) with odd \( j\) can be 
included into \( \varepsilon \mbox{\mathversion{bold}$W$}^{\varepsilon}_{m}\).
When \( j\) is even, \( 2(m-1)-j+2 \) is also even and we have
\begin{align*}
\partial^{j}_{s}\mbox{\mathversion{bold}$v$}^{\varepsilon}&\times 
\partial^{2(m-1)-j+2}_{s}\mbox{\mathversion{bold}$v$}^{\varepsilon}\\[3mm]
&=
\big\{
(\mbox{\mathversion{bold}$v$}^{\varepsilon}\cdot 
\partial^{j}_{s}\mbox{\mathversion{bold}$v$}^{\varepsilon})
\mbox{\mathversion{bold}$v$}^{\varepsilon}
-
\mbox{\mathversion{bold}$v$}^{\varepsilon}\times 
(\mbox{\mathversion{bold}$v$}^{\varepsilon}\times 
\partial^{j}_{s}\mbox{\mathversion{bold}$v$}^{\varepsilon})\big\}\\[3mm]
& \qquad \qquad \times
\big\{
(\mbox{\mathversion{bold}$v$}^{\varepsilon}\cdot 
\partial^{2(m-1)-j+2}_{s}\mbox{\mathversion{bold}$v$}^{\varepsilon})
\mbox{\mathversion{bold}$v$}^{\varepsilon}
-
\mbox{\mathversion{bold}$v$}^{\varepsilon}\times 
(\mbox{\mathversion{bold}$v$}^{\varepsilon}\times 
\partial^{2(m-1)-j+2}_{s}\mbox{\mathversion{bold}$v$}^{\varepsilon})\big\}\\[3mm]
&=
-(\mbox{\mathversion{bold}$v$}^{\varepsilon}\cdot 
\partial^{j}_{s}\mbox{\mathversion{bold}$v$}^{\varepsilon})
\mbox{\mathversion{bold}$v$}^{\varepsilon}\times
\{
\mbox{\mathversion{bold}$v$}^{\varepsilon}\times 
(\mbox{\mathversion{bold}$v$}^{\varepsilon}\times 
\partial^{2(m-1)-j+2}_{s}\mbox{\mathversion{bold}$v$}^{\varepsilon})\}\\[3mm]
& \qquad -
(\mbox{\mathversion{bold}$v$}^{\varepsilon}\cdot 
\partial^{2(m-1)-j+2}_{s}\mbox{\mathversion{bold}$v$}^{\varepsilon})
[\mbox{\mathversion{bold}$v$}^{\varepsilon}\times 
(\mbox{\mathversion{bold}$v$}^{\varepsilon}\times 
\partial^{j}_{s}\mbox{\mathversion{bold}$v$}^{\varepsilon})]
\times \mbox{\mathversion{bold}$v$}^{\varepsilon}\\[3mm]
& \qquad +
[\mbox{\mathversion{bold}$v$}^{\varepsilon}\times 
(\mbox{\mathversion{bold}$v$}^{\varepsilon}\times 
\partial^{j}_{s}\mbox{\mathversion{bold}$v$}^{\varepsilon})]\times 
[\mbox{\mathversion{bold}$v$}^{\varepsilon}\times 
(\mbox{\mathversion{bold}$v$}^{\varepsilon}\times 
\partial^{2(m-1)-j+2}_{s}\mbox{\mathversion{bold}$v$}^{\varepsilon})],
\end{align*}
where we have used the decomposition
\begin{eqnarray*}
\mbox{\mathversion{bold}$W$}=
(\mbox{\mathversion{bold}$v$}^{\varepsilon}\cdot 
\mbox{\mathversion{bold}$W$})\mbox{\mathversion{bold}$v$}^{\varepsilon}
-
\mbox{\mathversion{bold}$v$}^{\varepsilon}\times 
(\mbox{\mathversion{bold}$v$}^{\varepsilon}\times 
\mbox{\mathversion{bold}$W$})
\end{eqnarray*}
for any vector \( \mbox{\mathversion{bold}$W$}\), 
which follows from \( |\mbox{\mathversion{bold}$v$}^{\varepsilon}|\equiv 1 \).
From the assumption of induction for (\ref{vc2m}), we see that
\( \mbox{\mathversion{bold}$v$}^{\varepsilon}\times 
\partial^{k}_{s}\mbox{\mathversion{bold}$v$}^{\varepsilon} \)
can be 
included into \( \varepsilon \mbox{\mathversion{bold}$W$}^{\varepsilon}_{m}\)
for \( k=j \) and \( k=2(m-1)-j+2 \), and thus the terms in the summation
of (\ref{m-1}) with 
even \( j\) can also be included into \( \varepsilon 
\mbox{\mathversion{bold}$W$}^{\varepsilon}_{m}\).
Hence, we have
\begin{align*}
\mbox{\mathversion{bold}$v$}^{\varepsilon}\times 
\partial^{2(m-1)}_{s}\mbox{\mathversion{bold}$v$}^{\varepsilon}_{t}
&=(
\mbox{\mathversion{bold}$v$}^{\varepsilon}\cdot 
\partial^{2m}_{s}\mbox{\mathversion{bold}$v$}^{\varepsilon})
\mbox{\mathversion{bold}$v$}^{\varepsilon}
-
\partial^{2m}_{s}\mbox{\mathversion{bold}$v$}^{\varepsilon}
+
2(m-1)(\mbox{\mathversion{bold}$v$}^{\varepsilon}\cdot 
\partial^{2m-1}_{s}\mbox{\mathversion{bold}$v$}^{\varepsilon})
\mbox{\mathversion{bold}$v$}^{\varepsilon}
+
\varepsilon \mbox{\mathversion{bold}$v$}^{\varepsilon}\times
\partial^{2m}_{s}\mbox{\mathversion{bold}$v$}^{\varepsilon}\\[3mm]
& \qquad +
\varepsilon \mbox{\mathversion{bold}$W$}^{\varepsilon}_{m}.
\end{align*}
Substituting the above to (\ref{mt}), taking the exterior product 
with \( \mbox{\mathversion{bold}$v$}^{\varepsilon}\) from the left, 
and mutiplying by \( -1\) yields
\begin{align*}
\mbox{\mathversion{bold}$v$}^{\varepsilon}\times 
\partial^{2m}_{s}\mbox{\mathversion{bold}$v$}^{\varepsilon}|_{s=0,1}
&=
(-1)^{m}\varepsilon 
A_{m-1}\big\{
\partial^{m-1}_{t}\mbox{\mathversion{bold}$v$}^{\varepsilon}_{ss}
+
2(\mbox{\mathversion{bold}$v$}^{\varepsilon}_{s}\cdot 
\partial^{m-1}_{t}\mbox{\mathversion{bold}$v$}^{\varepsilon}_{s})
\mbox{\mathversion{bold}$v$}^{\varepsilon}
\big\}\\[3mm]
& \qquad +
\varepsilon \mbox{\mathversion{bold}$v$}^{\varepsilon}\times 
(\mbox{\mathversion{bold}$v$}^{\varepsilon}\times 
\partial^{2m}_{s}\mbox{\mathversion{bold}$v$}^{\varepsilon})
+
\varepsilon \bigg\{
\sum^{m-3}_{l=0}(-1)^{l+1}
A_{l+3}
\partial^{l+1}_{t}\partial^{2(m-1-l)}_{s}
\mbox{\mathversion{bold}$v$}^{\varepsilon}
\bigg\}\\[3mm]
& \qquad +
\varepsilon 
\mbox{\mathversion{bold}$W$}^{\varepsilon}_{m}\big|_{s=0,1}\\[3mm]
&=
(-1)^{m}\varepsilon 
A_{m-1}\big\{
\partial^{m-1}_{t}\mbox{\mathversion{bold}$v$}^{\varepsilon}_{ss}
+
2(\mbox{\mathversion{bold}$v$}^{\varepsilon}_{s}\cdot 
\partial^{m-1}_{t}\mbox{\mathversion{bold}$v$}^{\varepsilon}_{s})
\mbox{\mathversion{bold}$v$}^{\varepsilon}
\big\}\\[3mm]
& \qquad +
\varepsilon
\bigg\{
\sum^{m-2}_{l=0}
A_{l+2}\partial^{l}_{t}\partial^{2(m-l)}\mbox{\mathversion{bold}$v$}^{\varepsilon}\bigg\}
+
\varepsilon \mbox{\mathversion{bold}$W$}^{\varepsilon}_{m}\big|_{s=0,1},
\end{align*}
and this proves (\ref{vc2m}). Set the right-hand side of (\ref{vc2m})
as \( \varepsilon 
\mbox{\mathversion{bold}$Z$}^{\varepsilon}_{m}|_{s=0,1}\). 

To prove (\ref{oddinner}), take \( j\) and \( k\) such that
\( j+k = 2(m-1)+1 \). 
We assume without loss of generality that 
\( j\) is even and set \( j=2l \), which gives 
\( k= 2(m-1-l)+1 \).
We first consider the case \( j\neq 0\).
When \( l=1\), we have
\begin{align*}
\mbox{\mathversion{bold}$v$}^{\varepsilon}\times 
\mbox{\mathversion{bold}$v$}^{\varepsilon}_{ss}|_{s=0,1}
=
-\varepsilon (\mbox{\mathversion{bold}$v$}^{\varepsilon}_{ss}
+
|\mbox{\mathversion{bold}$v$}^{\varepsilon}_{s}|^{2}
\mbox{\mathversion{bold}$v$}^{\varepsilon})
\end{align*}
Taking the exterior product with 
\( \partial^{2(m-2)+1}_{s}\mbox{\mathversion{bold}$v$}^{\varepsilon}\) from the 
left yields
\begin{align*}
\partial^{2(m-2)+1}_{s}\mbox{\mathversion{bold}$v$}^{\varepsilon}\times 
(\mbox{\mathversion{bold}$v$}^{\varepsilon}\times 
\mbox{\mathversion{bold}$v$}^{\varepsilon}_{ss})|_{s=0,1}
=
\varepsilon 
\partial^{2(m-2)+1}_{s}
\mbox{\mathversion{bold}$v$}^{\varepsilon}\times 
(\mbox{\mathversion{bold}$v$}^{\varepsilon}_{ss}+
|\mbox{\mathversion{bold}$v$}^{\varepsilon}_{s}|^{2}
\mbox{\mathversion{bold}$v$}^{\varepsilon})|_{s=0,1}.
\end{align*}
Expanding the exterior products and rearranging the terms yield
\begin{align*}
(\mbox{\mathversion{bold}$v$}^{\varepsilon}_{ss}\cdot 
\partial^{2(m-2)+1}_{s}\mbox{\mathversion{bold}$v$}^{\varepsilon})
\mbox{\mathversion{bold}$v$}^{\varepsilon}|_{s=0,1}
&=
(\mbox{\mathversion{bold}$v$}^{\varepsilon}\cdot 
\partial^{2(m-2)+1}_{s}\mbox{\mathversion{bold}$v$}^{\varepsilon})
\mbox{\mathversion{bold}$v$}^{\varepsilon}_{ss}
+
\varepsilon 
\partial^{2(m-2)+1}_{s}
\mbox{\mathversion{bold}$v$}^{\varepsilon}\times 
(\mbox{\mathversion{bold}$v$}^{\varepsilon}_{ss}+
|\mbox{\mathversion{bold}$v$}^{\varepsilon}_{s}|^{2}
\mbox{\mathversion{bold}$v$}^{\varepsilon})|_{s=0,1}\\[3mm]
&=
\varepsilon Y_{0,2(m-2)+1}\mbox{\mathversion{bold}$v$}^{\varepsilon}_{ss}
+
\varepsilon 
\partial^{2(m-2)+1}_{s}
\mbox{\mathversion{bold}$v$}^{\varepsilon}\times 
(\mbox{\mathversion{bold}$v$}^{\varepsilon}_{ss}+
|\mbox{\mathversion{bold}$v$}^{\varepsilon}_{s}|^{2}
\mbox{\mathversion{bold}$v$}^{\varepsilon})|_{s=0,1}.
\end{align*}
Taking the inner product with \( \mbox{\mathversion{bold}$v$}^{\varepsilon}\) tells us that
\begin{align*}
\mbox{\mathversion{bold}$v$}^{\varepsilon}_{ss}
& \cdot 
\partial^{2(m-2)+1}_{s}\mbox{\mathversion{bold}$v$}^{\varepsilon}|_{s=0,1}\\[3mm]
&=
\varepsilon \big\{
Y_{0,2(m-2)+1}(\mbox{\mathversion{bold}$v$}^{\varepsilon}\cdot
\mbox{\mathversion{bold}$v$}^{\varepsilon}_{ss})
+
\mbox{\mathversion{bold}$v$}^{\varepsilon}\cdot
\big( 
\partial^{2(m-2)+1}_{s}
\mbox{\mathversion{bold}$v$}^{\varepsilon}\times 
(\mbox{\mathversion{bold}$v$}^{\varepsilon}_{ss}+
|\mbox{\mathversion{bold}$v$}^{\varepsilon}_{s}|^{2}
\mbox{\mathversion{bold}$v$}^{\varepsilon})
\big) \big\}|_{s=0,1},
\end{align*}
and setting \( Y_{2,2(m-2)+1} \) as the expression in brackets in the 
right-hand side proves the case \( l=1\).
From (\ref{vc2m}) we have for \( 2\leq l \leq m-1 \),
\begin{align*}
\mbox{\mathversion{bold}$v$}^{\varepsilon}\times 
\partial^{2l}_{s}\mbox{\mathversion{bold}$v$}^{\varepsilon}|_{s=0,1}
=
\varepsilon \mbox{\mathversion{bold}$Z$}^{\varepsilon}_{l}|_{s=0,1},
\end{align*}
and taking the exterior product with 
\( \partial^{2(m-1-l)+1}_{s}\mbox{\mathversion{bold}$v$}^{\varepsilon}\) from the left 
yields
\begin{align*}
\partial^{2(m-1-l)+1}_{s}\mbox{\mathversion{bold}$v$}^{\varepsilon}
\times 
(\mbox{\mathversion{bold}$v$}^{\varepsilon}\times 
\partial^{2l}_{s}\mbox{\mathversion{bold}$v$}^{\varepsilon})|_{s=0,1}
=
\varepsilon 
\partial^{2(m-1-l)+1}_{s}\mbox{\mathversion{bold}$v$}^{\varepsilon}\times
\mbox{\mathversion{bold}$Z$}^{\varepsilon}_{l}|_{s=0,1}.
\end{align*}
Expanding the exterior product as before, we have
\begin{align*}
(
\partial^{2(m-1-l)+1}_{s}\mbox{\mathversion{bold}$v$}^{\varepsilon}\cdot 
\partial^{2l}_{s}\mbox{\mathversion{bold}$v$}^{\varepsilon})
\mbox{\mathversion{bold}$v$}^{\varepsilon}|_{s=0,1}
=
\varepsilon \big\{
Y_{0,2(m-1-l)+1}
\partial^{2l}_{s}\mbox{\mathversion{bold}$v$}^{\varepsilon}
+
\partial^{2(m-1-l)+1}_{s}\mbox{\mathversion{bold}$v$}^{\varepsilon}\times
\mbox{\mathversion{bold}$Z$}^{\varepsilon}_{l}\big\}|_{s=0,1}.
\end{align*}
After taking the inner product of the above with 
\( \mbox{\mathversion{bold}$v$}^{\varepsilon}\), we see that setting
\begin{align*}
Y_{2l,2(m-1-l)+1}:= 
Y_{0,2(m-1-l)+1}
(\mbox{\mathversion{bold}$v$}^{\varepsilon}\cdot 
\partial^{2l}_{s}\mbox{\mathversion{bold}$v$}^{\varepsilon})
+
\mbox{\mathversion{bold}$v$}^{\varepsilon}\cdot
(
\partial^{2(m-1-l)+1}_{s}\mbox{\mathversion{bold}$v$}^{\varepsilon}\times
\mbox{\mathversion{bold}$Z$}^{\varepsilon}_{l}
)
\end{align*}
proves the case \( l\geq 2\). Finally, when \( l=0 \) we have
from (\ref{one}) that
\begin{align*}
\mbox{\mathversion{bold}$v$}^{\varepsilon}\cdot 
\partial^{2(m-1)+1}_{s}\mbox{\mathversion{bold}$v$}^{\varepsilon}
|_{s=0,1}
&=
-\frac{1}{2}
\sum^{2(m-1)}_{i=1}
\left(
\begin{array}{c}
2(m-1)+1\\
i
\end{array}\right)
\partial^{i}_{s}\mbox{\mathversion{bold}$v$}^{\varepsilon}\cdot 
\partial^{2(m-1)+1-i}_{s}\mbox{\mathversion{bold}$v$}^{\varepsilon}\big| _{s=0,1}\\[3mm]
&=
-\frac{1}{2}\sum^{m-1}_{l=1}
\left(
\begin{array}{c}
2(m-1)+1\\
2l
\end{array}\right)
\partial^{2l}_{s}\mbox{\mathversion{bold}$v$}^{\varepsilon}\cdot 
\partial^{2(m-1-l)+1}_{s}\mbox{\mathversion{bold}$v$}^{\varepsilon}\\[3mm]
& \qquad -
\frac{1}{2}\sum^{m-1}_{l=1}
\left(
\begin{array}{c}
2(m-1)+1\\
2l-1
\end{array}\right)
\partial^{2l-1}_{s}\mbox{\mathversion{bold}$v$}^{\varepsilon}\cdot 
\partial^{2(m-1-l)}_{s}\mbox{\mathversion{bold}$v$}^{\varepsilon}\big|_{s=0,1}\\[3mm]
&=
\varepsilon 
\bigg\{
-\frac{1}{2}\sum^{m-1}_{l=1}
\left(
\begin{array}{c}
2(m-1)+1\\
2l
\end{array}\right)
Y_{2l,2(m-1-l)+1}\\[3mm]
& \qquad -
\frac{1}{2}\sum^{m-1}_{l=1}
\left(
\begin{array}{c}
2(m-1)+1\\
2l-1
\end{array}\right)
Y_{2l-1,2(m-1-l)}
\bigg\}
\big|_{s=0,1},
\end{align*}
which shows that 
\begin{align*}
Y_{0,2(m-1)+1}&:=
-\frac{1}{2}\sum^{m-1}_{l=1}
\left(
\begin{array}{c}
2(m-1)+1\\
2l
\end{array}\right)
Y_{2l,2(m-1-l)+1}\\[3mm]
& \qquad -
\frac{1}{2}\sum^{m-1}_{l=1}
\left(
\begin{array}{c}
2(m-1)+1\\
2l-1
\end{array}\right)
Y_{2l-1,2(m-1-l)}
\end{align*}
proves the case \( l=0\) and this finishes the proof of (\ref{oddinner}) and the 
proof of the lemma. \hfill \( \Box \).

\bigskip
Now we prove the following uniform estimate.
\begin{pr}
\label{ue}
Let \( m\geq 1 \) be an integer, \( T>0\), and \( \mbox{\mathversion{bold}$v$}^{\varepsilon}\in 
Y^{m+2}_{T}(I)\) be the solution of {\rm (\ref{rnl})}. There exist \( T_{0}\in (0,T]\) and \( \varepsilon_{0}\in 
(0,\varepsilon_{\ast}] \) such that
for any \( \varepsilon \in (0,\varepsilon _{0}] \),
\begin{eqnarray*}
\|\mbox{\mathversion{bold}$v$}^{\varepsilon}\|_{Y^{m}_{T_{0}}(I)} \leq 
c_{\ast},
\end{eqnarray*}
where \( c_{\ast}>0\) depends on \( \| \mbox{\mathversion{bold}$v$}_{0s}\|_{2m}\) and \( T_{0}\), and
\( T_{0}\) depends on \( \| \mbox{\mathversion{bold}$v$}_{0s}\|_{2} \).
Here, \( \varepsilon_{\ast}>0 \) is given in Proposition {\rm \ref{id1}}.
\end{pr}
Note that the assumption 
\( \mbox{\mathversion{bold}$v$}^{\varepsilon}\in Y^{m+2}_{T}(I)\) is
not essential because we can approximate the initial datum as smooth as we need
and obtain solutions as smooth as necessary to justify the calculations below.

\noindent{\it Proof.}
By the interpolation inequality and integration by parts, we have
\begin{align*}
\frac{1}{2}\frac{{\rm d}}{{\rm d}t}\|\mbox{\mathversion{bold}$v$}^{\varepsilon}_{s}\|^{2}
=
(\mbox{\mathversion{bold}$v$}^{\varepsilon}_{s},\mbox{\mathversion{bold}$v$}^{\varepsilon}_{st})
&=
-(\mbox{\mathversion{bold}$v$}^{\varepsilon}_{ss},\mbox{\mathversion{bold}$v$}^{\varepsilon}_{t})+
[\mbox{\mathversion{bold}$v$}^{\varepsilon}_{s}\cdot \mbox{\mathversion{bold}$v$}^{\varepsilon}_{t}]^{1}_{s=0}\\[3mm]
&=-\varepsilon\|\mbox{\mathversion{bold}$v$}^{\varepsilon}_{ss}\|^{2}
-\varepsilon(\mbox{\mathversion{bold}$v$}^{\varepsilon}_{ss},|\mbox{\mathversion{bold}$v$}^{\varepsilon}_{s}|^{2}
\mbox{\mathversion{bold}$v$}^{\varepsilon})\\[3mm]
&\leq -\varepsilon\|\mbox{\mathversion{bold}$v$}^{\varepsilon}_{ss}\|^{2} 
+ \frac{\varepsilon}{2}\|\mbox{\mathversion{bold}$v$}^{\varepsilon}_{s}\|^{2}_{1} 
+ C\|\mbox{\mathversion{bold}$v$}^{\varepsilon}_{s}\|^{6},
\end{align*}
where \( [ \ \cdot \ ]^{1}_{s=0}\) is defined by \( [f]^{1}_{s=0} = f|_{s=1}
-f|_{s=0}\),
and we have used \( \mbox{\mathversion{bold}$v$}^{\varepsilon}_{t}|_{s=0,1}
= \mbox{\mathversion{bold}$0$}\).
We further calculate
\begin{align*}
\frac{1}{2}
\frac{{\rm d}}{{\rm d}t}\|\mbox{\mathversion{bold}$v$}^{\varepsilon}_{ss}\|^{2}
&= -( \mbox{\mathversion{bold}$v$}^{\varepsilon}_{sss},
\mbox{\mathversion{bold}$v$}^{\varepsilon}_{ts})
+
[\mbox{\mathversion{bold}$v$}^{\varepsilon}_{ss}\cdot \mbox{\mathversion{bold}$v$}^{\varepsilon}_{ts}]^{1}_{s=0}\\[3mm]
&=-(\mbox{\mathversion{bold}$v$}^{\varepsilon}_{sss},\mbox{\mathversion{bold}$v$}^{\varepsilon}_{s}\times \mbox{\mathversion{bold}$v$}^{\varepsilon}_{ss})
-\varepsilon 
\| \mbox{\mathversion{bold}$v$}^{\varepsilon}_{sss}\|^{2}
-2\varepsilon
(\mbox{\mathversion{bold}$v$}^{\varepsilon}_{sss},
(\mbox{\mathversion{bold}$v$}^{\varepsilon}_{s}\cdot 
\mbox{\mathversion{bold}$v$}^{\varepsilon}_{ss})
\mbox{\mathversion{bold}$v$}^{\varepsilon})\\[3mm]
& \qquad + 
-\varepsilon(\mbox{\mathversion{bold}$v$}^{\varepsilon}_{sss},
|\mbox{\mathversion{bold}$v$}^{\varepsilon}_{s}|^{2}
\mbox{\mathversion{bold}$v$}^{\varepsilon}_{s})
+
[\mbox{\mathversion{bold}$v$}^{\varepsilon}_{ss}\cdot 
\mbox{\mathversion{bold}$v$}^{\varepsilon}_{ts}]^{1}_{s=0}.
\end{align*}
From equation (\ref{decom}) and (\ref{one}), we have
\begin{align*}
\mbox{\mathversion{bold}$v$}^{\varepsilon}_{sss}\cdot 
(\mbox{\mathversion{bold}$v$}^{\varepsilon}_{s}\times
\mbox{\mathversion{bold}$v$}^{\varepsilon}_{ss})
&=
-\mbox{\mathversion{bold}$v$}^{\varepsilon}_{ss}\cdot 
(\mbox{\mathversion{bold}$v$}^{\varepsilon}_{s}\times 
\mbox{\mathversion{bold}$v$}^{\varepsilon}_{sss})\\[3mm]
&= -\mbox{\mathversion{bold}$v$}^{\varepsilon}_{ss}\cdot 
\big\{
-(\mbox{\mathversion{bold}$v$}^{\varepsilon}\cdot 
\mbox{\mathversion{bold}$v$}^{\varepsilon}_{sss})
\mbox{\mathversion{bold}$v$}^{\varepsilon}\times
\mbox{\mathversion{bold}$v$}^{\varepsilon}_{s}
+
[(\mbox{\mathversion{bold}$v$}^{\varepsilon}\times
\mbox{\mathversion{bold}$v$}^{\varepsilon}_{s})\cdot
\mbox{\mathversion{bold}$v$}^{\varepsilon}_{sss}]
\mbox{\mathversion{bold}$v$}^{\varepsilon}\big\}\\[3mm]
&=
-\mbox{\mathversion{bold}$v$}^{\varepsilon}_{ss}\cdot
\big\{
3(\mbox{\mathversion{bold}$v$}^{\varepsilon}_{s}\cdot 
\mbox{\mathversion{bold}$v$}^{\varepsilon}_{ss})
\mbox{\mathversion{bold}$v$}^{\varepsilon}\times 
\mbox{\mathversion{bold}$v$}^{\varepsilon}_{s}
+
[(\mbox{\mathversion{bold}$v$}^{\varepsilon}\times
\mbox{\mathversion{bold}$v$}^{\varepsilon}_{s}\cdot
\mbox{\mathversion{bold}$v$}^{\varepsilon}_{sss}]
\mbox{\mathversion{bold}$v$}^{\varepsilon} \big\}\\[3mm]
&=-3(\mbox{\mathversion{bold}$v$}^{\varepsilon}_{s}\cdot
\mbox{\mathversion{bold}$v$}^{\varepsilon}_{ss})
\mbox{\mathversion{bold}$v$}^{\varepsilon}_{ss}\cdot
(\mbox{\mathversion{bold}$v$}^{\varepsilon}\times \mbox{\mathversion{bold}$v$}^{\varepsilon}_{s})
+
|\mbox{\mathversion{bold}$v$}^{\varepsilon}_{s}|^{2}
(\mbox{\mathversion{bold}$v$}^{\varepsilon}\times 
\mbox{\mathversion{bold}$v$}^{\varepsilon}_{s})\cdot
\mbox{\mathversion{bold}$v$}^{\varepsilon}_{sss},
\end{align*}
which allows us to further calculate
\begin{align*}
\frac{1}{2}
\frac{{\rm d}}{{\rm d}t}\|\mbox{\mathversion{bold}$v$}^{\varepsilon}_{ss}\|^{2}
&\leq 
C\|\mbox{\mathversion{bold}$v$}^{\varepsilon}_{s}\|_{1}^{2}
(1+\|\mbox{\mathversion{bold}$v$}^{\varepsilon}_{s}\|_{1}^{2})
-\frac{\varepsilon}{2}\| \mbox{\mathversion{bold}$v$}^{\varepsilon}_{sss}\|^{2}
-(\mbox{\mathversion{bold}$v$}^{\varepsilon}_{sss},
|\mbox{\mathversion{bold}$v$}^{\varepsilon}_{s}|^{2}
\mbox{\mathversion{bold}$v$}^{\varepsilon}\times 
\mbox{\mathversion{bold}$v$}^{\varepsilon}_{s})
+[\mbox{\mathversion{bold}$v$}^{\varepsilon}_{ss}\cdot 
\mbox{\mathversion{bold}$v$}^{\varepsilon}_{ts}]^{1}_{s=0}\\[3mm]
&=C\|\mbox{\mathversion{bold}$v$}^{\varepsilon}_{s}\|^{2}
(1+\|\mbox{\mathversion{bold}$v$}^{\varepsilon}_{s}\|^{2})
-\frac{\varepsilon}{2}\| \mbox{\mathversion{bold}$v$}^{\varepsilon}_{sss}\|^{2}
+2(\mbox{\mathversion{bold}$v$}^{\varepsilon}_{ss},
(\mbox{\mathversion{bold}$v$}^{\varepsilon}_{s}\cdot 
\mbox{\mathversion{bold}$v$}^{\varepsilon}_{ss})
\mbox{\mathversion{bold}$v$}^{\varepsilon}\times 
\mbox{\mathversion{bold}$v$}^{\varepsilon}_{s})\\[3mm]
& \qquad +[|\mbox{\mathversion{bold}$v$}^{\varepsilon}_{s}|^{2}
\mbox{\mathversion{bold}$v$}^{\varepsilon}_{ss}\cdot 
(\mbox{\mathversion{bold}$v$}^{\varepsilon}\times 
\mbox{\mathversion{bold}$v$}^{\varepsilon}_{s})]^{1}_{s=0}
+
[\mbox{\mathversion{bold}$v$}^{\varepsilon}_{ss}\cdot 
\mbox{\mathversion{bold}$v$}^{\varepsilon}_{ts}]^{1}_{s=0}.
\end{align*}
We have seen that
\begin{align*}
\mbox{\mathversion{bold}$v$}^{\varepsilon}_{ss}|_{s=0,1}
=
-\varepsilon \mbox{\mathversion{bold}$v$}^{\varepsilon}\times
\mbox{\mathversion{bold}$v$}^{\varepsilon}_{ss}
-
|\mbox{\mathversion{bold}$v$}^{\varepsilon}_{s}|^{2}
\mbox{\mathversion{bold}$v$}^{\varepsilon}
|_{s=0,1}.
\end{align*}
Hence, we have
\begin{align*}
\big| 
[|\mbox{\mathversion{bold}$v$}^{\varepsilon}_{s}|^{2}
\mbox{\mathversion{bold}$v$}^{\varepsilon}_{ss}\cdot 
(\mbox{\mathversion{bold}$v$}^{\varepsilon}\times 
\mbox{\mathversion{bold}$v$}^{\varepsilon}_{s})]^{1}_{s=0}
\big|
&=
\big|
[
-\varepsilon |\mbox{\mathversion{bold}$v$}^{\varepsilon}_{s}|^{2}
(\mbox{\mathversion{bold}$v$}^{\varepsilon}\times 
\mbox{\mathversion{bold}$v$}^{\varepsilon}_{s})\cdot 
(\mbox{\mathversion{bold}$v$}^{\varepsilon}\times 
\mbox{\mathversion{bold}$v$}^{\varepsilon}_{ss})
]^{1}_{s=0}
\big| \\[3mm]
&\leq
\frac{\varepsilon}{8}\|\mbox{\mathversion{bold}$v$}^{\varepsilon}_{ss}\|^{2}_{1}
+
C\|\mbox{\mathversion{bold}$v$}^{\varepsilon}_{s}\|^{6}_{1}
\end{align*}
Furthermore, we have
\begin{align*}
\mbox{\mathversion{bold}$v$}^{\varepsilon}_{ts}=
\mbox{\mathversion{bold}$v$}^{\varepsilon}\times 
\mbox{\mathversion{bold}$v$}^{\varepsilon}_{sss}
+ \varepsilon \mbox{\mathversion{bold}$v$}^{\varepsilon}_{sss}
+2\varepsilon (\mbox{\mathversion{bold}$v$}^{\varepsilon}_{s}\cdot
\mbox{\mathversion{bold}$v$}^{\varepsilon}_{ss})
\mbox{\mathversion{bold}$v$}^{\varepsilon}
+
\varepsilon |\mbox{\mathversion{bold}$v$}^{\varepsilon}_{s}|^{2}
\mbox{\mathversion{bold}$v$}^{\varepsilon}_{s},
\end{align*}
and combining with (\ref{decom}), (\ref{one}), and (\ref{bv2}), we have
\begin{align*}
\mbox{\mathversion{bold}$v$}^{\varepsilon}_{ss}\cdot
\mbox{\mathversion{bold}$v$}^{\varepsilon}_{ts}|_{s=0,1}
&=
-\varepsilon (\mbox{\mathversion{bold}$v$}^{\varepsilon}\times 
\mbox{\mathversion{bold}$v$}^{\varepsilon}_{ss})\cdot 
(\mbox{\mathversion{bold}$v$}^{\varepsilon}\times 
\mbox{\mathversion{bold}$v$}^{\varepsilon}_{sss})
-
\varepsilon^{2}(\mbox{\mathversion{bold}$v$}^{\varepsilon}\times 
\mbox{\mathversion{bold}$v$}^{\varepsilon}_{ss})\cdot 
\mbox{\mathversion{bold}$v$}^{\varepsilon}_{sss}
+
\mbox{\mathversion{bold}$\eta $}_{2}^{\varepsilon }|_{s=0,1}
\end{align*}
where \( \mbox{\mathversion{bold}$\eta $}_{2}^{\varepsilon}\) are terms 
which can be estimated as
\begin{align*}
\big| \mbox{\mathversion{bold}$\eta $}_{2}^{\varepsilon}|_{s=0,1}\big|
\leq
\frac{\varepsilon}{8}
\| \mbox{\mathversion{bold}$v$}^{\varepsilon}_{s}\|^{2}_{2}
+
\| \mbox{\mathversion{bold}$v$}^{\varepsilon}_{s}\|^{6}_{1},
\end{align*}
and we arrive at
\begin{align*}
\frac{1}{2}
\frac{{\rm d}}{{\rm d}t}\|\mbox{\mathversion{bold}$v$}^{\varepsilon}_{ss}\|^{2}
&\leq 
C\|
\mbox{\mathversion{bold}$v$}^{\varepsilon}_{s}\|_{1}^{2}
(1+\|\mbox{\mathversion{bold}$v$}^{\varepsilon}_{s}\|_{1}^{4})
-
\frac{\varepsilon}{4}\|
\mbox{\mathversion{bold}$v$}^{\varepsilon}_{sss}\|^{2}\\[3mm]
& \qquad -
[\varepsilon (\mbox{\mathversion{bold}$v$}^{\varepsilon}\times 
\mbox{\mathversion{bold}$v$}^{\varepsilon}_{ss})\cdot 
(\mbox{\mathversion{bold}$v$}^{\varepsilon}\times 
\mbox{\mathversion{bold}$v$}^{\varepsilon}_{sss})]^{1}_{s=0}
-
[\varepsilon^{2}(\mbox{\mathversion{bold}$v$}^{\varepsilon}\times 
\mbox{\mathversion{bold}$v$}^{\varepsilon}_{ss})\cdot 
\mbox{\mathversion{bold}$v$}^{\varepsilon}_{sss}]^{1}_{s=0}.
\end{align*}
We further estimate \( \mbox{\mathversion{bold}$v$}^{\varepsilon}_{sss}\) to 
close the estimate.
\begin{align*}
\frac{1}{2}
\frac{{\rm d}}{{\rm d}t}
\| \mbox{\mathversion{bold}$v$}^{\varepsilon}_{sss}\|^{2}
&=
-( \partial^{4}_{s}\mbox{\mathversion{bold}$v$}^{\varepsilon},
\mbox{\mathversion{bold}$v$}^{\varepsilon}_{tss})
+
[
\mbox{\mathversion{bold}$v$}^{\varepsilon}_{sss}\cdot 
\mbox{\mathversion{bold}$v$}^{\varepsilon}_{tss}]^{1}_{s=0}\\[3mm]
&=-2(\partial^{4}_{s}\mbox{\mathversion{bold}$v$}^{\varepsilon},
\mbox{\mathversion{bold}$v$}^{\varepsilon}_{s}\times 
\mbox{\mathversion{bold}$v$}^{\varepsilon}_{sss})
-
\varepsilon \| \partial^{4}_{s}\mbox{\mathversion{bold}$v$}^{\varepsilon} \|^{2}
-
\varepsilon (\partial^{4}_{s}\mbox{\mathversion{bold}$v$}^{\varepsilon},
|\mbox{\mathversion{bold}$v$}^{\varepsilon}_{s}|^{2}
\mbox{\mathversion{bold}$v$}^{\varepsilon}_{ss})\\[3mm]
& \qquad -
2\varepsilon (\partial^{4}_{s}\mbox{\mathversion{bold}$v$}^{\varepsilon},
(\mbox{\mathversion{bold}$v$}^{\varepsilon}_{s}\cdot 
\mbox{\mathversion{bold}$v$}^{\varepsilon}_{sss})
\mbox{\mathversion{bold}$v$}^{\varepsilon})
-
2\varepsilon (\partial ^{4}_{s}\mbox{\mathversion{bold}$v$}^{\varepsilon},
|\mbox{\mathversion{bold}$v$}^{\varepsilon}_{ss}|^{2}
\mbox{\mathversion{bold}$v$}^{\varepsilon})
- 4\varepsilon(\partial^{4}_{s}\mbox{\mathversion{bold}$v$}^{\varepsilon},
(\mbox{\mathversion{bold}$v$}^{\varepsilon}_{s}\cdot 
\mbox{\mathversion{bold}$v$}^{\varepsilon}_{ss})
\mbox{\mathversion{bold}$v$}^{\varepsilon}_{s})\\[3mm]
& \qquad +
[\mbox{\mathversion{bold}$v$}^{\varepsilon}_{sss}\cdot 
\mbox{\mathversion{bold}$v$}^{\varepsilon}_{tss}]^{1}_{s=0}.
\end{align*}
From (\ref{decom}) and (\ref{one}), we see that 
\begin{align*}
(\partial^{4}_{s}\mbox{\mathversion{bold}$v$}^{\varepsilon},
\mbox{\mathversion{bold}$v$}^{\varepsilon}_{s}\times 
\mbox{\mathversion{bold}$v$}^{\varepsilon}_{sss})
&=
-4(\mbox{\mathversion{bold}$v$}^{\varepsilon}_{sss},
(\mbox{\mathversion{bold}$v$}^{\varepsilon}_{s}\cdot 
\mbox{\mathversion{bold}$v$}^{\varepsilon}_{sss})
\mbox{\mathversion{bold}$v$}^{\varepsilon}\times 
\mbox{\mathversion{bold}$v$}^{\varepsilon}_{s})
-
3(\mbox{\mathversion{bold}$v$}^{\varepsilon}_{sss},
|\mbox{\mathversion{bold}$v$}^{\varepsilon}_{ss}|^{2}
\mbox{\mathversion{bold}$v$}^{\varepsilon}\times 
\mbox{\mathversion{bold}$v$}^{\varepsilon}_{s})\\[3mm]
& \qquad +
3(\mbox{\mathversion{bold}$v$}^{\varepsilon}_{ss},
[(\mbox{\mathversion{bold}$v$}^{\varepsilon}\times 
\mbox{\mathversion{bold}$v$}^{\varepsilon}_{s})\cdot 
\partial^{4}_{s}\mbox{\mathversion{bold}$v$}^{\varepsilon}]
\mbox{\mathversion{bold}$v$}^{\varepsilon}_{s}),
\end{align*}
which yields
\begin{align*}
\frac{1}{2}
\frac{{\rm d}}{{\rm d}t}
\| \mbox{\mathversion{bold}$v$}^{\varepsilon}_{sss}\|^{2}
&\leq 
C\|\mbox{\mathversion{bold}$v$}^{\varepsilon}_{s}\|_{2}^{2}
(1+\|\mbox{\mathversion{bold}$v$}^{\varepsilon}_{s}\|_{2}^{2})
-
\frac{\varepsilon}{2}
\|\partial^{4}_{s}\mbox{\mathversion{bold}$v$}^{\varepsilon}\|^{2}
-
6(\mbox{\mathversion{bold}$v$}^{\varepsilon}_{ss},
[(\mbox{\mathversion{bold}$v$}^{\varepsilon}\times 
\mbox{\mathversion{bold}$v$}^{\varepsilon}_{s})\cdot 
\partial^{4}_{s}\mbox{\mathversion{bold}$v$}^{\varepsilon}]
\mbox{\mathversion{bold}$v$}^{\varepsilon}_{s})\\[3mm]
& \qquad +
[\mbox{\mathversion{bold}$v$}^{\varepsilon}_{sss}\cdot 
\mbox{\mathversion{bold}$v$}^{\varepsilon}_{tss}]^{1}_{s=0}\\[3mm]
&\leq 
C\|\mbox{\mathversion{bold}$v$}^{\varepsilon}_{s}\|_{2}^{2}
(1+\|\mbox{\mathversion{bold}$v$}^{\varepsilon}_{s}\|_{2}^{2})
-
\frac{\varepsilon}{2}
\|\partial^{4}_{s}\mbox{\mathversion{bold}$v$}^{\varepsilon}\|^{2}
-
6\big[
(\mbox{\mathversion{bold}$v$}^{\varepsilon}_{s}\cdot 
\mbox{\mathversion{bold}$v$}^{\varepsilon}_{ss})
[(\mbox{\mathversion{bold}$v$}^{\varepsilon}\times 
\mbox{\mathversion{bold}$v$}^{\varepsilon}_{s})\cdot 
\mbox{\mathversion{bold}$v$}^{\varepsilon}_{sss}]
\big]^{1}_{s=0}\\[3mm]
& \qquad + 
[ \mbox{\mathversion{bold}$v$}^{\varepsilon}_{sss}\cdot 
\mbox{\mathversion{bold}$v$}^{\varepsilon}_{tss} ]^{1}_{s=0},
\end{align*}
where integration by parts was used. We see from (\ref{bv2}) 
and \( \mbox{\mathversion{bold}$v$}^{\varepsilon}\cdot 
\mbox{\mathversion{bold}$v$}^{\varepsilon}_{s}\equiv 0 \) that
\begin{align*}
\mbox{\mathversion{bold}$v$}^{\varepsilon}_{s}\cdot 
\mbox{\mathversion{bold}$v$}^{\varepsilon}_{ss}|_{s=0,1}
=
-\varepsilon \mbox{\mathversion{bold}$v$}^{\varepsilon}_{s}\cdot 
(\mbox{\mathversion{bold}$v$}^{\varepsilon}_{s}\times 
\mbox{\mathversion{bold}$v$}^{\varepsilon}_{ss})|_{s=0,1}
\end{align*}
holds, and hence
\begin{align*}
\frac{1}{2}
\frac{{\rm d}}{{\rm d}t}
\| \mbox{\mathversion{bold}$v$}^{\varepsilon}_{sss}\|^{2}
&\leq 
C\|\mbox{\mathversion{bold}$v$}^{\varepsilon}_{s}\|_{2}^{2}
(1+\|\mbox{\mathversion{bold}$v$}^{\varepsilon}_{s}\|_{2}^{6})
-
\frac{\varepsilon}{4}
\|\partial^{4}_{s}\mbox{\mathversion{bold}$v$}^{\varepsilon}\|^{2}
+
[
\mbox{\mathversion{bold}$v$}^{\varepsilon}_{sss}\cdot 
\mbox{\mathversion{bold}$v$}^{\varepsilon}_{tss}]^{1}_{s=0}.
\end{align*}
We further investigate the boundary term.
From \( \mbox{\mathversion{bold}$v$}^{\varepsilon}_{tt}|_{s=0,1}
=
\mbox{\mathversion{bold}$0$}\), we have
\begin{align*}
&\mbox{\mathversion{bold}$v$}^{\varepsilon}\times 
( \mbox{\mathversion{bold}$v$}^{\varepsilon}\times 
\partial^{4}_{s}\mbox{\mathversion{bold}$v$}^{\varepsilon})
+
2\varepsilon \mbox{\mathversion{bold}$v$}^{\varepsilon}\times 
\partial^{4}_{s}\mbox{\mathversion{bold}$v$}^{\varepsilon}
+
\varepsilon^{2}\partial^{4}_{s}\mbox{\mathversion{bold}$v$}^{\varepsilon}
+2\varepsilon \mbox{\mathversion{bold}$v$}^{\varepsilon}_{s}\times 
\mbox{\mathversion{bold}$v$}^{\varepsilon}_{sss}\\[3mm]
& \qquad +
2\varepsilon^{2}(\mbox{\mathversion{bold}$v$}^{\varepsilon}_{s}\cdot 
\mbox{\mathversion{bold}$v$}^{\varepsilon}_{sss})
\mbox{\mathversion{bold}$v$}^{\varepsilon}+
2\varepsilon[\mbox{\mathversion{bold}$v$}^{\varepsilon}_{s}\cdot 
(\mbox{\mathversion{bold}$v$}^{\varepsilon}\times 
\mbox{\mathversion{bold}$v$}^{\varepsilon}_{sss})]
\mbox{\mathversion{bold}$v$}^{\varepsilon}
+
\varepsilon \mbox{\mathversion{bold}$w$}^{\varepsilon}_{2}|_{s=0,1}
=\mbox{\mathversion{bold}$0$},
\end{align*}
from which we further deduce that
\begin{align*}
(1-\varepsilon^{2})\partial^{4}_{s}\mbox{\mathversion{bold}$v$}^{\varepsilon}
&-
2\varepsilon\mbox{\mathversion{bold}$v$}^{\varepsilon}\times 
\partial^{4}_{s}\mbox{\mathversion{bold}$v$}^{\varepsilon} |_{s=0,1}
=
-4(\mbox{\mathversion{bold}$v$}^{\varepsilon}_{s}\cdot 
\mbox{\mathversion{bold}$v$}^{\varepsilon}_{sss})
\mbox{\mathversion{bold}$v$}^{\varepsilon}
-3|\mbox{\mathversion{bold}$v$}^{\varepsilon}_{ss}|^{2}
\mbox{\mathversion{bold}$v$}^{\varepsilon}
-6(\mbox{\mathversion{bold}$v$}^{\varepsilon}_{s}\cdot 
\mbox{\mathversion{bold}$v$}^{\varepsilon}_{ss})
\mbox{\mathversion{bold}$v$}^{\varepsilon}_{s}\\[3mm]
& \qquad +
2\varepsilon \mbox{\mathversion{bold}$v$}^{\varepsilon}_{s}\times 
\mbox{\mathversion{bold}$v$}^{\varepsilon}_{sss}
+
4\varepsilon ^{2}(\mbox{\mathversion{bold}$v$}^{\varepsilon}_{s}\cdot 
\mbox{\mathversion{bold}$v$}^{\varepsilon}_{sss})
\mbox{\mathversion{bold}$v$}^{\varepsilon}
+
2\varepsilon[\mbox{\mathversion{bold}$v$}^{\varepsilon}_{s}\cdot 
(\mbox{\mathversion{bold}$v$}^{\varepsilon}\times 
\mbox{\mathversion{bold}$v$}^{\varepsilon}_{sss})]
\mbox{\mathversion{bold}$v$}^{\varepsilon}
+
\varepsilon \mbox{\mathversion{bold}$w$}^{\varepsilon}_{2}|_{s=0,1}
\end{align*}
where \( \mbox{\mathversion{bold}$w$}^{\varepsilon}_{2}\) are terms 
that can be estimated as
\begin{align*}
\big| 
\mbox{\mathversion{bold}$w$}^{\varepsilon}_{2}|_{s=0,1}
\big|
\leq 
C\| \mbox{\mathversion{bold}$v$}^{\varepsilon}_{s}\|_{2}
( 1 + \| \mbox{\mathversion{bold}$v$}^{\varepsilon}_{s}\|^{2}_{2}),
\end{align*}
and the exact form may change from line to line.
Since we have \( | \mbox{\mathversion{bold}$v$}^{\varepsilon}|\equiv 1\)
, there exists \( \varepsilon_{1}\in (0,1) \) such that 
for any \( \varepsilon \in (0,\varepsilon_{1}]\), the matrix
\( I_{3}-\frac{2\varepsilon}{1-\varepsilon^{2}}
A(\mbox{\mathversion{bold}$v$}^{\varepsilon})\) is reversible 
and the inverse matrix can be expressed as
\begin{align*}
\bigg( I_{3}-\frac{2\varepsilon}{1-\varepsilon^{2}}
A(\mbox{\mathversion{bold}$v$}^{\varepsilon})\bigg)^{-1}
=
I_{3}+D_{4}(\mbox{\mathversion{bold}$v$}^{\varepsilon}),
\end{align*}
where \( D_{4}(\mbox{\mathversion{bold}$v$}^{\varepsilon})\) is given by
\begin{align*}
D_{4}(\mbox{\mathversion{bold}$v$}^{\varepsilon})=
\sum^{\infty}_{j=1}\bigg(\frac{2\varepsilon}{1-\varepsilon^{2}}\bigg)^{j}
A(\mbox{\mathversion{bold}$v$}^{\varepsilon})^{j}.
\end{align*}
Here, \( A(\mbox{\mathversion{bold}$v$}^{\varepsilon})
\mbox{\mathversion{bold}$W$}= 
\mbox{\mathversion{bold}$v$}^{\varepsilon}\times 
\mbox{\mathversion{bold}$W$}\) for a vector \( 
\mbox{\mathversion{bold}$W$}\).
We arrive at
\begin{align*}
\partial^{4}_{s}\mbox{\mathversion{bold}$v$}^{\varepsilon}|_{s=0,1}
&=
\frac{1}{1-\varepsilon^{2}}\big( I_{3}
+D_{4}(\mbox{\mathversion{bold}$v$}^{\varepsilon})\big)
\big\{-4(\mbox{\mathversion{bold}$v$}^{\varepsilon}_{s}\cdot 
\mbox{\mathversion{bold}$v$}^{\varepsilon}_{sss})
\mbox{\mathversion{bold}$v$}^{\varepsilon}
-3|\mbox{\mathversion{bold}$v$}^{\varepsilon}_{ss}|^{2}
\mbox{\mathversion{bold}$v$}^{\varepsilon}
-6(\mbox{\mathversion{bold}$v$}^{\varepsilon}_{s}\cdot 
\mbox{\mathversion{bold}$v$}^{\varepsilon}_{ss})
\mbox{\mathversion{bold}$v$}^{\varepsilon}_{s}\\[3mm]
& \qquad +
2\varepsilon \mbox{\mathversion{bold}$v$}^{\varepsilon}_{s}\times 
\mbox{\mathversion{bold}$v$}^{\varepsilon}_{sss}
+
4\varepsilon ^{2}(\mbox{\mathversion{bold}$v$}^{\varepsilon}_{s}\cdot 
\mbox{\mathversion{bold}$v$}^{\varepsilon}_{sss})
\mbox{\mathversion{bold}$v$}^{\varepsilon}
+
2\varepsilon[\mbox{\mathversion{bold}$v$}^{\varepsilon}_{s}\cdot 
(\mbox{\mathversion{bold}$v$}^{\varepsilon}\times 
\mbox{\mathversion{bold}$v$}^{\varepsilon}_{sss})]
\mbox{\mathversion{bold}$v$}^{\varepsilon}
+
\varepsilon \mbox{\mathversion{bold}$w$}^{\varepsilon}_{2}\big\}
|_{s=0,1}.
\end{align*}
Since,
\begin{align*}
D_{4}(\mbox{\mathversion{bold}$v$}^{\varepsilon})
&=
\frac{2 \varepsilon }{1-\varepsilon ^{2}}
\bigg(
\sum^{\infty}_{j=0}
\bigg(\frac{2\varepsilon}{1-\varepsilon^{2}}\bigg)^{j}
A(\mbox{\mathversion{bold}$v$}^{\varepsilon})^{j}\bigg)
A(\mbox{\mathversion{bold}$v$}^{\varepsilon})\\[3mm]
&=:
\frac{2 \varepsilon }{1-\varepsilon ^{2}}
B_{4}(\mbox{\mathversion{bold}$v$}^{\varepsilon})
A(\mbox{\mathversion{bold}$v$}^{\varepsilon}),
\end{align*}
we have
\begin{align*}
\partial^{4}_{s}\mbox{\mathversion{bold}$v$}^{\varepsilon}|_{s=0,1}
&=
\frac{1}{1-\varepsilon^{2}}
\big\{
-4(\mbox{\mathversion{bold}$v$}^{\varepsilon}_{s}\cdot 
\mbox{\mathversion{bold}$v$}^{\varepsilon}_{sss})
\mbox{\mathversion{bold}$v$}^{\varepsilon}
-3|\mbox{\mathversion{bold}$v$}^{\varepsilon}_{ss}|^{2}
\mbox{\mathversion{bold}$v$}^{\varepsilon}
-6(\mbox{\mathversion{bold}$v$}^{\varepsilon}_{s}\cdot 
\mbox{\mathversion{bold}$v$}^{\varepsilon}_{ss})
\mbox{\mathversion{bold}$v$}^{\varepsilon}_{s}\\[3mm]
& \qquad +
2\varepsilon \mbox{\mathversion{bold}$v$}^{\varepsilon}_{s}\times 
\mbox{\mathversion{bold}$v$}^{\varepsilon}_{sss}
+
4\varepsilon ^{2}(\mbox{\mathversion{bold}$v$}^{\varepsilon}_{s}\cdot 
\mbox{\mathversion{bold}$v$}^{\varepsilon}_{sss})
\mbox{\mathversion{bold}$v$}^{\varepsilon}
+
2\varepsilon[\mbox{\mathversion{bold}$v$}^{\varepsilon}_{s}\cdot 
(\mbox{\mathversion{bold}$v$}^{\varepsilon}\times 
\mbox{\mathversion{bold}$v$}^{\varepsilon}_{sss})]
\mbox{\mathversion{bold}$v$}^{\varepsilon}
\big\}\\[3mm]
& \qquad +
\frac{2\varepsilon }{(1-\varepsilon^{2})^{2}}B_{4}(\mbox{\mathversion{bold}$v$}^{\varepsilon})
\big\{
2\varepsilon \mbox{\mathversion{bold}$v$}^{\varepsilon}\times 
(\mbox{\mathversion{bold}$v$}^{\varepsilon}_{s}\times 
\mbox{\mathversion{bold}$v$}^{\varepsilon}_{sss})
-
6(\mbox{\mathversion{bold}$v$}^{\varepsilon}_{s}\cdot 
\mbox{\mathversion{bold}$v$}^{\varepsilon}_{ss})
\mbox{\mathversion{bold}$v$}^{\varepsilon}_{s}
\big\}\\[3mm]
& \qquad +
\frac{\varepsilon }{1-\varepsilon ^{2}}
\big( I_{3}+ D_{4}(\mbox{\mathversion{bold}$v$}^{\varepsilon})\big)
\mbox{\mathversion{bold}$w$}^{\varepsilon}_{2}
\big|_{s=0,1}.
\end{align*}
From \( \mbox{\mathversion{bold}$v$}^{\varepsilon}_{s}\cdot 
\mbox{\mathversion{bold}$v$}^{\varepsilon}_{ss}|_{s=0,1}
=
-\varepsilon \mbox{\mathversion{bold}$v$}^{\varepsilon}_{s}\cdot 
(\mbox{\mathversion{bold}$v$}^{\varepsilon}\times 
\mbox{\mathversion{bold}$v$}^{\varepsilon}_{ss} )|_{s=0,1}\)
and (\ref{one}), we finally have
\begin{align*}
\partial^{4}_{s}\mbox{\mathversion{bold}$v$}^{\varepsilon}|_{s=0,1}
&=
\frac{1}{1-\varepsilon^{2}}
\big\{
-4(\mbox{\mathversion{bold}$v$}^{\varepsilon}_{s}\cdot 
\mbox{\mathversion{bold}$v$}^{\varepsilon}_{sss})
\mbox{\mathversion{bold}$v$}^{\varepsilon}
-3|\mbox{\mathversion{bold}$v$}^{\varepsilon}_{ss}|^{2}
\mbox{\mathversion{bold}$v$}^{\varepsilon}
+
2\varepsilon \mbox{\mathversion{bold}$v$}^{\varepsilon}_{s}\times 
\mbox{\mathversion{bold}$v$}^{\varepsilon}_{sss}\\[3mm]
& \qquad +
4\varepsilon ^{2}(\mbox{\mathversion{bold}$v$}^{\varepsilon}_{s}\cdot 
\mbox{\mathversion{bold}$v$}^{\varepsilon}_{sss})
\mbox{\mathversion{bold}$v$}^{\varepsilon}
+
2\varepsilon[\mbox{\mathversion{bold}$v$}^{\varepsilon}_{s}\cdot 
(\mbox{\mathversion{bold}$v$}^{\varepsilon}\times 
\mbox{\mathversion{bold}$v$}^{\varepsilon}_{sss})]
\mbox{\mathversion{bold}$v$}^{\varepsilon}
\big\}\\[3mm]
& \qquad +
\frac{2\varepsilon }{(1-\varepsilon^{2})^{2}}B_{4}(\mbox{\mathversion{bold}$v$}^{\varepsilon})
\big\{
2\varepsilon \mbox{\mathversion{bold}$v$}^{\varepsilon}\times 
(\mbox{\mathversion{bold}$v$}^{\varepsilon}_{s}\times 
\mbox{\mathversion{bold}$v$}^{\varepsilon}_{sss})
\big\}
+
\varepsilon 
\mbox{\mathversion{bold}$w$}^{\varepsilon}_{2}
\big|_{s=0,1} \\[3mm]
&=
\frac{1}{1-\varepsilon^{2}}
\big\{
-4(\mbox{\mathversion{bold}$v$}^{\varepsilon}_{s}\cdot 
\mbox{\mathversion{bold}$v$}^{\varepsilon}_{sss})
\mbox{\mathversion{bold}$v$}^{\varepsilon}
-3|\mbox{\mathversion{bold}$v$}^{\varepsilon}_{ss}|^{2}
\mbox{\mathversion{bold}$v$}^{\varepsilon}
+
2\varepsilon \mbox{\mathversion{bold}$v$}^{\varepsilon}_{s}\times 
\mbox{\mathversion{bold}$v$}^{\varepsilon}_{sss}\\[3mm]
& \qquad +
4\varepsilon ^{2}(\mbox{\mathversion{bold}$v$}^{\varepsilon}_{s}\cdot 
\mbox{\mathversion{bold}$v$}^{\varepsilon}_{sss})
\mbox{\mathversion{bold}$v$}^{\varepsilon}
+
2\varepsilon[\mbox{\mathversion{bold}$v$}^{\varepsilon}_{s}\cdot 
(\mbox{\mathversion{bold}$v$}^{\varepsilon}\times 
\mbox{\mathversion{bold}$v$}^{\varepsilon}_{sss})]
\mbox{\mathversion{bold}$v$}^{\varepsilon}
\big\}\\[3mm]
& \qquad +
\frac{4\varepsilon ^{2} }{(1-\varepsilon^{2})^{2}}B_{4}(\mbox{\mathversion{bold}$v$}^{\varepsilon})
\{
[\mbox{\mathversion{bold}$v$}^{\varepsilon}_{s}\cdot
(
\mbox{\mathversion{bold}$v$}^{\varepsilon}\times 
\mbox{\mathversion{bold}$v$}^{\varepsilon}_{sss})]
\mbox{\mathversion{bold}$v$}^{\varepsilon}\}
+
\varepsilon 
\mbox{\mathversion{bold}$w$}^{\varepsilon}_{2}
\big|_{s=0,1},
\end{align*}
where we have used the fact that
\begin{align*}
\|D_{4}(\mbox{\mathversion{bold}$v$}^{\varepsilon})\|_{L^{\infty}(I\times [0,T])} + 
\|B_{4}(\mbox{\mathversion{bold}$v$}^{\varepsilon})\|_{L^{\infty}(I\times [0,T])} \leq C
\end{align*}
holds with \( C>0\) independent of 
\( \varepsilon \) because \( |\mbox{\mathversion{bold}$v$}^{\varepsilon}|\equiv 1 \).
The above expression along with (\ref{one}) and
\begin{align*}
\mbox{\mathversion{bold}$v$}^{\varepsilon}_{tss}|_{s=0,1}
=
\mbox{\mathversion{bold}$v$}^{\varepsilon}\times
\partial^{4}_{s}\mbox{\mathversion{bold}$v$}^{\varepsilon}+
2\mbox{\mathversion{bold}$v$}^{\varepsilon}_{s}\times 
\mbox{\mathversion{bold}$v$}^{\varepsilon}_{sss}+
\varepsilon \partial^{4}_{s}\mbox{\mathversion{bold}$v$}^{\varepsilon}
+2\varepsilon(\mbox{\mathversion{bold}$v$}^{\varepsilon}_{s}\cdot 
\mbox{\mathversion{bold}$v$}^{\varepsilon}_{sss})
\mbox{\mathversion{bold}$v$}^{\varepsilon}
+
\varepsilon \tilde{\mbox{\mathversion{bold}$w$}^{\varepsilon}}_{2},
\end{align*}
where \( \tilde{\mbox{\mathversion{bold}$w$}}^{\varepsilon}_{2}\) are 
terms that satisfy the same form of estimate as 
\( \mbox{\mathversion{bold}$w$}^{\varepsilon}_{2}\),
yields the following estimate. 
\begin{align*}
\big|
[\mbox{\mathversion{bold}$v$}^{\varepsilon}_{sss}\cdot 
\mbox{\mathversion{bold}$v$}^{\varepsilon}_{tss}|_{s=0,1}
\big|
\leq 
\frac{\varepsilon}{8}\|\mbox{\mathversion{bold}$v$}^{\varepsilon}_{sss}\|^{2}_{1}
+
C\|\mbox{\mathversion{bold}$v$}^{\varepsilon}_{s}\|^{6}_{2},
\end{align*}
where we also utilized (\ref{decom}) and (\ref{one}).
The above estimate allows us to close the energy estimate as
\begin{align*}
\frac{1}{2}\frac{{\rm d}}{{\rm d}t}
\| \mbox{\mathversion{bold}$v$}^{\varepsilon}_{sss}\|^{2}
\leq 
C\|\mbox{\mathversion{bold}$v$}^{\varepsilon}_{s}\|^{2}_{s}
(1+ \| \mbox{\mathversion{bold}$v$}^{\varepsilon}_{s}\|^{6}_{2})
-
\frac{\varepsilon}{8}
\|\partial^{4}_{s}\mbox{\mathversion{bold}$v$}^{\varepsilon}\|^{2}.
\end{align*}
Combining the estimates for \( \mbox{\mathversion{bold}$v$}^{\varepsilon}_{s}\),
\( \mbox{\mathversion{bold}$v$}^{\varepsilon}_{ss}\), and 
\( \mbox{\mathversion{bold}$v$}^{\varepsilon}_{sss}\) yields
\begin{align*}
\frac{{\rm d}}{{\rm d}t}
\| \mbox{\mathversion{bold}$v$}^{\varepsilon}_{s}\|^{2}_{2}
&\leq 
C\|\mbox{\mathversion{bold}$v$}^{\varepsilon}_{s}\|^{2}_{2}
(1+\|\mbox{\mathversion{bold}$v$}^{\varepsilon}\|^{6}_{2})\\[3mm]
&\leq 
C(1+\|\mbox{\mathversion{bold}$v$}^{\varepsilon}\|^{2}_{2})^{4}.
\end{align*}
Comparing \( \|\mbox{\mathversion{bold}$v$}^{\varepsilon}_{s}\|^{2}_{2}\)
with the solution of the ordinary differential equation given by
\begin{align*}
\left\{
\begin{array}{ll}
R_{t}=C(1+R)^{4}, & t>0, \\[3mm]
R(0)=\|\mbox{\mathversion{bold}$v$}^{\varepsilon}_{0}\|_{2}^{2}, & \ 
\end{array}\right.
\end{align*}
we see that
there exists a \( T_{0} \in (0,T) \) depending on 
\( \| \mbox{\mathversion{bold}$v$}_{0s}\|_{2}\) such that 
\begin{align*}
\| \mbox{\mathversion{bold}$v$}^{\varepsilon}_{s}\|_{2}
\leq 
c_{0},
\end{align*}
where \( c_{0}>0 \) is a constant depending on 
\( \| \mbox{\mathversion{bold}$v$}^{\varepsilon}_{0s}\|_{2} \) and 
\( T_{0}\).
Finally, we see that
\begin{align*}
\| \mbox{\mathversion{bold}$v$}^{\varepsilon}_{t}\|
+
\| \mbox{\mathversion{bold}$v$}^{\varepsilon}_{ts}\|
\leq 
c_{1},
\end{align*}
holds from the equation of (\ref{rnl}), where the 
dependence of \( c_{1}\) is the same as \( c_{0}\).
This finishes the proof of the proposition for \( m=1\). 

We proceed by induction.
Suppose the statement of the proposition is true up to 
\( m-1 \) for some \( m\geq 2\). From this assumption of induction, 
a solution \( \mbox{\mathversion{bold}$v$}^{\varepsilon}\in 
Y^{m+2}_{T}(I) \)
satisfies
\begin{align*}
\| \mbox{\mathversion{bold}$v$}^{\varepsilon}\|_{Y^{m-1}_{T_{0}}(I)}
\leq c_{\ast},
\end{align*}
where \( c_{\ast}>0\) depends on 
\( \| \mbox{\mathversion{bold}$v$}_{0}\|_{2(m-1)}\) and 
\( T_{0}\), and \( T_{0}\) depends on 
\( \| \mbox{\mathversion{bold}$v$}_{0s}\|_{2}\).
Since the estimates for the \( t\) derivatives of 
\( \mbox{\mathversion{bold}$v$}^{\varepsilon}\) can be obtained 
from the equation, it is sufficient to prove the uniform estimate 
for \( \partial^{2m}_{s}\mbox{\mathversion{bold}$v$}^{\varepsilon}\) and 
\( \partial^{2m+1}_{s}\mbox{\mathversion{bold}$v$}^{\varepsilon}\).
We have for \( k_{0}=2m \ \mbox{or} \ 2m+1 \),
\begin{align*}
\frac{1}{2}\frac{{\rm d}}{{\rm d}t}
\|
\partial^{k_{0}}_{s}\mbox{\mathversion{bold}$v$}^{\varepsilon}
\|^{2}
&=
(\partial^{k_{0}}_{s}\mbox{\mathversion{bold}$v$}^{\varepsilon},
\partial^{k_{0}}_{s}\mbox{\mathversion{bold}$v$}^{\varepsilon}_{t})\\[3mm]
&=
(\partial^{k_{0}}_{s}\mbox{\mathversion{bold}$v$}^{\varepsilon},
\mbox{\mathversion{bold}$v$}^{\varepsilon}\times 
\partial^{k_{0}+2}_{s}\mbox{\mathversion{bold}$v$}^{\varepsilon})
+
k_{0}( \partial^{k_{0}}_{s}\mbox{\mathversion{bold}$v$}^{\varepsilon},
\mbox{\mathversion{bold}$v$}^{\varepsilon}_{s}\times
\partial^{k_{0}+1}_{s}\mbox{\mathversion{bold}$v$}^{\varepsilon})\\[3mm]
& \qquad +
\varepsilon(\partial^{k_{0}}_{s}\mbox{\mathversion{bold}$v$}^{\varepsilon},
\partial^{k_{0}+2}_{s}\mbox{\mathversion{bold}$v$}^{\varepsilon}))
+
k_{0}\varepsilon
(\partial^{k_{0}}_{s}\mbox{\mathversion{bold}$v$}^{\varepsilon},
(\mbox{\mathversion{bold}$v$}^{\varepsilon}_{s}\cdot 
\partial^{k_{0}+1}_{s}\mbox{\mathversion{bold}$v$}^{\varepsilon})
\mbox{\mathversion{bold}$v$}^{\varepsilon})
+
(\partial^{k_{0}}_{s}\mbox{\mathversion{bold}$v$}^{\varepsilon},
\mbox{\mathversion{bold}$V$}^{\varepsilon}_{k_{0}})\\[3mm]
&=
(k_{0}-1)(\partial^{k_{0}}_{s}\mbox{\mathversion{bold}$v$}^{\varepsilon},
\mbox{\mathversion{bold}$v$}^{\varepsilon}_{s}\times 
\partial^{k_{0}+1}_{s}\mbox{\mathversion{bold}$v$}^{\varepsilon})
-
\varepsilon 
\|\partial^{k_{0}+1}_{s}\mbox{\mathversion{bold}$v$}^{\varepsilon}\|^{2}\\[3mm]
& \qquad +
k_{0}\varepsilon 
(\partial^{k_{0}}_{s}\mbox{\mathversion{bold}$v$}^{\varepsilon},
(\mbox{\mathversion{bold}$v$}^{\varepsilon}_{s}\cdot 
\partial^{k_{0}+1}_{s}\mbox{\mathversion{bold}$v$}^{\varepsilon})
\mbox{\mathversion{bold}$v$}^{\varepsilon})
+
[
\partial^{k_{0}}_{s}\mbox{\mathversion{bold}$v$}^{\varepsilon}\cdot
(\mbox{\mathversion{bold}$v$}^{\varepsilon}\times 
\partial^{k_{0}+1}_{s}\mbox{\mathversion{bold}$v$}^{\varepsilon})
]^{1}_{s=0}\\[3mm]
& \qquad +
\varepsilon 
[
\partial^{k_{0}}_{s}\mbox{\mathversion{bold}$v$}^{\varepsilon}\cdot 
\partial^{k_{0}+1}_{s}\mbox{\mathversion{bold}$v$}^{\varepsilon}
]^{1}_{s=0}
+
(\partial^{k_{0}}_{s}\mbox{\mathversion{bold}$v$}^{\varepsilon},
\mbox{\mathversion{bold}$V$}^{\varepsilon}_{k_{0}}),
\end{align*}
where \( \mbox{\mathversion{bold}$V$}^{\varepsilon}_{k_{0}}\) are terms that can 
be estimated as
\begin{eqnarray*}
\| \mbox{\mathversion{bold}$V$}^{\varepsilon}_{k_{0}}\|
\leq 
C\| \mbox{\mathversion{bold}$v$}^{\varepsilon}_{s}\|_{k_{0}-1},
\end{eqnarray*}
with \( C>0\) depending on
\( \| \mbox{\mathversion{bold}$v$}^{\varepsilon}_{s}\|_{2(m-1)} \), and 
are harmless for the current estimate.
From (\ref{decom}) and (\ref{one}), we have
\begin{align*}
\mbox{\mathversion{bold}$v$}^{\varepsilon}_{s}\times 
\partial^{k_{0}+1}_{s}\mbox{\mathversion{bold}$v$}^{\varepsilon}
&=
-(\mbox{\mathversion{bold}$v$}^{\varepsilon}\cdot
\partial^{k_{0}+1}_{s}\mbox{\mathversion{bold}$v$}^{\varepsilon})
\mbox{\mathversion{bold}$v$}^{\varepsilon}\times 
\mbox{\mathversion{bold}$v$}^{\varepsilon}_{s}
+
[(\mbox{\mathversion{bold}$v$}^{\varepsilon}\times 
\mbox{\mathversion{bold}$v$}^{\varepsilon}_{s})\cdot 
\partial^{k_{0}+1}_{s}\mbox{\mathversion{bold}$v$}^{\varepsilon}
]
\mbox{\mathversion{bold}$v$}^{\varepsilon}\\[3mm]
&=
\frac{1}{2}\sum^{k_{0}}_{j=1}
\left(
\begin{array}{c}
k_{0}+1\\
j
\end{array}\right)
(
\partial^{j}_{s}\mbox{\mathversion{bold}$v$}^{\varepsilon}\cdot 
\partial^{k_{0}+1-j}_{s}\mbox{\mathversion{bold}$v$}^{\varepsilon})
\mbox{\mathversion{bold}$v$}^{\varepsilon}\times 
\mbox{\mathversion{bold}$v$}^{\varepsilon}_{s}
+
[(\mbox{\mathversion{bold}$v$}^{\varepsilon}\times 
\mbox{\mathversion{bold}$v$}^{\varepsilon}_{s})\cdot 
\partial^{k_{0}+1}_{s}\mbox{\mathversion{bold}$v$}^{\varepsilon}
]
\mbox{\mathversion{bold}$v$}^{\varepsilon},
\end{align*}
which gives 
\begin{align*}
\frac{1}{2}\frac{{\rm d}}{{\rm d}t}
\|
\partial^{k_{0}}_{s}\mbox{\mathversion{bold}$v$}^{\varepsilon}
\|^{2}
&=
\frac{(k_{0}-1)}{2}\sum^{k_{0}}_{j=1}
\left(
\begin{array}{c}
k_{0}+1\\
j
\end{array}\right)
(\partial^{k_{0}}_{s}\mbox{\mathversion{bold}$v$}^{\varepsilon},
(\partial^{j}_{s}\mbox{\mathversion{bold}$v$}^{\varepsilon}\cdot 
\partial^{k_{0}+1-j}_{s}\mbox{\mathversion{bold}$v$}^{\varepsilon})
\mbox{\mathversion{bold}$v$}^{\varepsilon}\times 
\mbox{\mathversion{bold}$v$}^{\varepsilon}_{s})\\[3mm]
& \qquad +
(k_{0}-1)(\partial^{k_{0}}_{s}\mbox{\mathversion{bold}$v$}^{\varepsilon},
[(\mbox{\mathversion{bold}$v$}^{\varepsilon}\times 
\mbox{\mathversion{bold}$v$}^{\varepsilon}_{s})\cdot 
\partial^{k_{0}+1}_{s}\mbox{\mathversion{bold}$v$}^{\varepsilon}]
\mbox{\mathversion{bold}$v$}^{\varepsilon})
-
\varepsilon
\| \partial^{k_{0}+1}_{s}\mbox{\mathversion{bold}$v$}^{\varepsilon}\|^{2}\\[3mm]
& \qquad +
k_{0}\varepsilon ( \partial^{k_{0}}_{s}\mbox{\mathversion{bold}$v$}^{\varepsilon},
(\mbox{\mathversion{bold}$v$}^{\varepsilon}_{s}\cdot 
\partial^{k_{0}+1}_{s}\mbox{\mathversion{bold}$v$}^{\varepsilon})
\mbox{\mathversion{bold}$v$}^{\varepsilon})
+
[\partial^{k_{0}}_{s}\mbox{\mathversion{bold}$v$}^{\varepsilon}\cdot 
(\mbox{\mathversion{bold}$v$}^{\varepsilon}\times 
\partial^{k_{0}+1}_{s}\mbox{\mathversion{bold}$v$}^{\varepsilon})]^{1}_{s=0}\\[3mm]
& \qquad +
\varepsilon
[\partial^{k_{0}}_{s}\mbox{\mathversion{bold}$v$}^{\varepsilon}\cdot 
\partial^{k_{0}+1}_{s}\mbox{\mathversion{bold}$v$}^{\varepsilon}]^{1}_{s=0}
+
(\partial^{k_{0}}_{s}\mbox{\mathversion{bold}$v$}^{\varepsilon},
\mbox{\mathversion{bold}$V$}^{\varepsilon}_{k_{0}})\\[3mm]
& \leq
C\|\partial^{k_{0}}_{s}\mbox{\mathversion{bold}$v$}^{\varepsilon}\|^{2}
-
\frac{\varepsilon}{2}
\| \partial^{k_{0}+1}_{s}\mbox{\mathversion{bold}$v$}^{\varepsilon}\|^{2}
+
(k_{0}-1)(\partial^{k_{0}}_{s}\mbox{\mathversion{bold}$v$}^{\varepsilon},
[(\mbox{\mathversion{bold}$v$}^{\varepsilon}\times 
\mbox{\mathversion{bold}$v$}^{\varepsilon}_{s})\cdot 
\partial^{k_{0}+1}_{s}\mbox{\mathversion{bold}$v$}^{\varepsilon})]
\mbox{\mathversion{bold}$v$}^{\varepsilon})\\[3mm]
& \qquad +
[\partial^{k_{0}}_{s}\mbox{\mathversion{bold}$v$}^{\varepsilon}\cdot 
(\mbox{\mathversion{bold}$v$}^{\varepsilon}\times 
\partial^{k_{0}+1}_{s}\mbox{\mathversion{bold}$v$}^{\varepsilon})]^{1}_{s=0}
+
\varepsilon
[
\partial^{k_{0}}_{s}\mbox{\mathversion{bold}$v$}^{\varepsilon}\cdot 
\partial^{k_{0}+1}_{s}\mbox{\mathversion{bold}$v$}^{\varepsilon}]^{1}_{s=0}
+
(\partial^{k_{0}}_{s}\mbox{\mathversion{bold}$v$}^{\varepsilon},
\mbox{\mathversion{bold}$V$}^{\varepsilon}_{k_{0}}),
\end{align*}
where \( C>0\) depends on 
\( \|\mbox{\mathversion{bold}$v$}^{\varepsilon}_{s}\|_{2(m-1)}\).
From (\ref{one}), we further calculate \newpage
\begin{align*}
\frac{1}{2}\frac{{\rm d}}{{\rm d}t}
\|
\partial^{k_{0}}_{s}\mbox{\mathversion{bold}$v$}^{\varepsilon}
\|^{2}
&\leq 
C\|\partial^{k_{0}}_{s}\mbox{\mathversion{bold}$v$}^{\varepsilon}\|^{2}
-
\frac{\varepsilon}{2}
\| \partial^{k_{0}+1}_{s}\mbox{\mathversion{bold}$v$}^{\varepsilon}\|^{2}\\[3mm]
& \qquad -
\frac{(k_{0}-1)}{2}\sum^{k_{0}-1}_{j=1}
\left(
\begin{array}{c}
k_{0}\\
j
\end{array}\right)
(
\partial^{k_{0}-j}_{s}\mbox{\mathversion{bold}$v$}^{\varepsilon},
[(\mbox{\mathversion{bold}$v$}^{\varepsilon}\times
\mbox{\mathversion{bold}$v$}^{\varepsilon}_{s})\cdot
\partial^{k_{0}+1}_{s}\mbox{\mathversion{bold}$v$}^{\varepsilon}]
\partial^{j}_{s}\mbox{\mathversion{bold}$v$}^{\varepsilon})\\[3mm]
& \qquad +
[\partial^{k_{0}}_{s}\mbox{\mathversion{bold}$v$}^{\varepsilon}\cdot 
(\mbox{\mathversion{bold}$v$}^{\varepsilon}\times 
\partial^{k_{0}+1}_{s}\mbox{\mathversion{bold}$v$}^{\varepsilon})]^{1}_{s=0}
+
\varepsilon [\partial^{k_{0}}_{s}\mbox{\mathversion{bold}$v$}^{\varepsilon}\cdot 
\partial^{k_{0}+1}_{s}\mbox{\mathversion{bold}$v$}^{\varepsilon}]^{1}_{s=0}
+
(\partial^{k_{0}}_{s}\mbox{\mathversion{bold}$v$}^{\varepsilon},
\mbox{\mathversion{bold}$V$}^{\varepsilon}_{k_{0}})\\[3mm]
&=
C\|\partial^{k_{0}}_{s}\mbox{\mathversion{bold}$v$}^{\varepsilon}\|^{2}
-
\frac{\varepsilon}{2}
\| \partial^{k_{0}+1}_{s}\mbox{\mathversion{bold}$v$}^{\varepsilon}\|^{2}\\[3mm]
& \qquad +
\frac{(k_{0}-1)}{2}\sum^{k_{0}-1}_{j=1}
\left(
\begin{array}{c}
k_{0}\\
j
\end{array}\right)
(\partial^{k_{0}-j+1}_{s}\mbox{\mathversion{bold}$v$}^{\varepsilon},
[(\mbox{\mathversion{bold}$v$}^{\varepsilon}\times 
\mbox{\mathversion{bold}$v$}^{\varepsilon}_{s})\cdot 
\partial^{k_{0}}_{s}\mbox{\mathversion{bold}$v$}^{\varepsilon}]
\partial^{j}_{s}\mbox{\mathversion{bold}$v$}^{\varepsilon})\\[3mm]
& \qquad +
\frac{(k_{0}-1)}{2}\sum^{k_{0}s-1}_{j=1}
\left(
\begin{array}{c}
k_{0}\\
j
\end{array}\right)
(\partial^{k_{0}-j}_{s}\mbox{\mathversion{bold}$v$}^{\varepsilon},
[(\mbox{\mathversion{bold}$v$}^{\varepsilon}\times 
\mbox{\mathversion{bold}$v$}^{\varepsilon}_{ss})\cdot 
\partial^{k_{0}}_{s}\mbox{\mathversion{bold}$v$}^{\varepsilon}]
\partial^{j}_{s}\mbox{\mathversion{bold}$v$}^{\varepsilon})\\[3mm]
& \qquad +
\frac{(k_{0}-1)}{2}\sum^{k_{0}-1}_{j=1}
\left(
\begin{array}{c}
k_{0}\\
j
\end{array}\right)
(\partial^{k_{0}-j}_{s}\mbox{\mathversion{bold}$v$}^{\varepsilon},
[(\mbox{\mathversion{bold}$v$}^{\varepsilon}\times 
\mbox{\mathversion{bold}$v$}^{\varepsilon}_{s})\cdot 
\partial^{k_{0}}_{s}\mbox{\mathversion{bold}$v$}^{\varepsilon}]
\partial^{j+1}_{s}\mbox{\mathversion{bold}$v$}^{\varepsilon})\\[3mm]
& \qquad -
\frac{(k_{0}-1)}{2}\sum^{k_{0}-1}_{j=1}
\left(
\begin{array}{c}
k_{0}\\
j
\end{array}\right)
\big[
(\partial^{k_{0}-j}_{s}\mbox{\mathversion{bold}$v$}^{\varepsilon}\cdot 
\partial^{j}_{s}\mbox{\mathversion{bold}$v$}^{\varepsilon})
[(\mbox{\mathversion{bold}$v$}^{\varepsilon}\times 
\mbox{\mathversion{bold}$v$}^{\varepsilon}_{s})\cdot 
\partial^{k_{0}}_{s}\mbox{\mathversion{bold}$v$}^{\varepsilon}]\big]^{1}_{s=0}\\[3mm]
& \qquad +
[
\partial^{k_{0}}_{s}\mbox{\mathversion{bold}$v$}^{\varepsilon}\cdot 
(\mbox{\mathversion{bold}$v$}^{\varepsilon}\times 
\partial^{k_{0}+1}_{s}\mbox{\mathversion{bold}$v$}^{\varepsilon})]^{1}_{s=0}
+
\varepsilon [\partial^{k_{0}}_{s}\mbox{\mathversion{bold}$v$}^{\varepsilon}\cdot 
\partial^{k_{0}+1}_{s}\mbox{\mathversion{bold}$v$}^{\varepsilon}]^{1}_{s=0}
+
(\partial^{k_{0}}_{s}\mbox{\mathversion{bold}$v$}^{\varepsilon},
\mbox{\mathversion{bold}$V$}^{\varepsilon}_{k_{0}})\\[3mm]
& \leq 
C\|\mbox{\mathversion{bold}$v$}^{\varepsilon}_{s}\|_{k_{0}-1}^{2}
-
\frac{\varepsilon}{2}
\|\partial^{k_{0}+1}_{s}\mbox{\mathversion{bold}$v$}^{\varepsilon}\|^{2}\\[3mm]
& \qquad -
\frac{(k_{0}-1)}{2}\sum^{k_{0}-1}_{j=1}
\left(
\begin{array}{c}
k_{0}\\
j
\end{array}\right)
\big[
(\partial^{k_{0}-j}_{s}\mbox{\mathversion{bold}$v$}^{\varepsilon}\cdot 
\partial^{j}_{s}\mbox{\mathversion{bold}$v$}^{\varepsilon})
[(\mbox{\mathversion{bold}$v$}^{\varepsilon}\times 
\mbox{\mathversion{bold}$v$}^{\varepsilon}_{s})\cdot 
\partial^{k_{0}}_{s}\mbox{\mathversion{bold}$v$}^{\varepsilon}]\big]^{1}_{s=0}\\[3mm]
& \qquad +
[
\partial^{k_{0}}_{s}\mbox{\mathversion{bold}$v$}^{\varepsilon}\cdot 
(\mbox{\mathversion{bold}$v$}^{\varepsilon}\times 
\partial^{k_{0}+1}_{s}\mbox{\mathversion{bold}$v$}^{\varepsilon})]^{1}_{s=0}
+
\varepsilon [\partial^{k_{0}}_{s}\mbox{\mathversion{bold}$v$}^{\varepsilon}\cdot 
\partial^{k_{0}+1}_{s}\mbox{\mathversion{bold}$v$}^{\varepsilon}]^{1}_{s=0}
+
(\partial^{k_{0}}_{s}\mbox{\mathversion{bold}$v$}^{\varepsilon},
\mbox{\mathversion{bold}$V$}^{\varepsilon}_{k_{0}}),
\end{align*}
where \( C>0\) depends on 
\( \| \mbox{\mathversion{bold}$v$}^{\varepsilon}_{s}\|_{2(m-1)} \).
We make use of Lemma \ref{lm2m2m-1} and \ref{lmvc2m} to estimate the 
boundary terms. 
First we focus on \(
\varepsilon [\partial^{k_{0}}_{s}\mbox{\mathversion{bold}$v$}^{\varepsilon}\cdot 
\partial^{k_{0}+1}_{s}\mbox{\mathversion{bold}$v$}^{\varepsilon}]^{1}_{s=0}
\).
Setting \( s=0,1\) in (\ref{2m2m-1}) with \( m\) as \( m+1\), we have
\begin{align*}
\mbox{\mathversion{bold}$0$}
&=
\sum^{m+1}_{j=0}a_{m+1,j}\varepsilon^{j}
A_{m+1-j}\partial^{2(m+1)}_{s}\mbox{\mathversion{bold}$v$}^{\varepsilon}
+
\sum^{m}_{j=0}
\sum^{m-j}_{k=0}
e_{m+1,j,k}\varepsilon^{j}
A_{m-j-k}\mbox{\mathversion{bold}$v$}^{\varepsilon}_{s}\times 
\big( A_{k}\partial^{2(m+1)-1}_{s}\mbox{\mathversion{bold}$v$}^{\varepsilon}\big)\\[3mm]
& \qquad +
\sum^{m+1}_{j=1}b_{m+1,j}\varepsilon^{m+2-j}
[\mbox{\mathversion{bold}$v$}^{\varepsilon}_{s}\cdot 
(A_{j-1}\partial^{2(m+1)-1}_{s}\mbox{\mathversion{bold}$v$}^{\varepsilon})]
\mbox{\mathversion{bold}$v$}^{\varepsilon}
+
\mbox{\mathversion{bold}$U$}^{\varepsilon}_{m+1}
\big|_{s=0,1}.
\end{align*}
As we showed in Section 3, for \( i\geq 1 \) we have
\begin{align*}
A_{i}\partial^{2(m+1)-1}_{s}\mbox{\mathversion{bold}$v$}^{\varepsilon}
=
\left\{
\begin{array}{ll}
(-1)^{n}\mbox{\mathversion{bold}$v$}^{\varepsilon}\times 
\partial^{2(m+1)-1}_{s}\mbox{\mathversion{bold}$v$}^{\varepsilon}, & i=2n+1,\\[3mm]
(-1)^{n+1}\big(
(\mbox{\mathversion{bold}$v$}^{\varepsilon}\cdot \partial^{2(m+1)-1}_{s}
\mbox{\mathversion{bold}$v$}^{\varepsilon})\mbox{\mathversion{bold}$v$}^{\varepsilon}
-
\partial^{2(m+1)-1}_{s}\mbox{\mathversion{bold}$v$}^{\varepsilon}\big), & i=2n,
\end{array}\right.
\end{align*}
which gives
\begin{align*}
\mbox{\mathversion{bold}$v$}^{\varepsilon}_{s}\times 
(A_{i}\partial^{2(m+1)-1}_{s}\mbox{\mathversion{bold}$v$}^{\varepsilon})
=
\left\{
\begin{array}{ll}
(-1)^{n}(\mbox{\mathversion{bold}$v$}^{\varepsilon}_{s}\cdot 
\partial^{2(m+1)-1}_{s}\mbox{\mathversion{bold}$v$}^{\varepsilon})
\mbox{\mathversion{bold}$v$}^{\varepsilon}, & 
i=2n+1,\\[3mm]
(-1)^{n+1}\big\{
(\mbox{\mathversion{bold}$v$}^{\varepsilon}\cdot \partial^{2(m+1)-1}_{s}
\mbox{\mathversion{bold}$v$}^{\varepsilon})\mbox{\mathversion{bold}$v$}^{\varepsilon}_{s}\times 
\mbox{\mathversion{bold}$v$}^{\varepsilon} & \ \\[3mm]
\hspace*{4.5cm} -
\mbox{\mathversion{bold}$v$}^{\varepsilon}_{s}\times 
\partial^{2(m+1)-1}_{s}\mbox{\mathversion{bold}$v$}^{\varepsilon}\big\}, &
i=2n.
\end{array}\right.
\end{align*}
Finally, if we take the exterior product of the above with 
\( \mbox{\mathversion{bold}$v$}^{\varepsilon}\) from the left, we see that 
\(\mbox{\mathversion{bold}$v$}^{\varepsilon}\times \big(
\mbox{\mathversion{bold}$v$}^{\varepsilon}_{s}\times 
(A_{i}\partial^{2(m+1)-1}_{s}\mbox{\mathversion{bold}$v$}^{\varepsilon})\big)\)
is identically zero. Hence we have shown that
\begin{align*}
A_{k}\mbox{\mathversion{bold}$v$}^{\varepsilon}_{s}\times 
(A_{i}\partial^{2(m+1)-1}_{s}\mbox{\mathversion{bold}$v$}^{\varepsilon})
=
\mbox{\mathversion{bold}$0$}
\end{align*}
if \( k,i\geq 1\). Hence we have
\begin{align*}
\mbox{\mathversion{bold}$0$}
&=
\sum^{m+1}_{j=0}a_{m+1,j}\varepsilon^{j}
A_{m+1-j}\partial^{2m+2}_{s}\mbox{\mathversion{bold}$v$}^{\varepsilon}\\[3mm]
& \qquad +
\sum^{m}_{j=0}\varepsilon^{j}
\big\{
e_{m+1,j,0}A_{m-j}(\mbox{\mathversion{bold}$v$}^{\varepsilon}_{s}\times
\partial^{2m+1}_{s}\mbox{\mathversion{bold}$v$}^{\varepsilon})
+
e_{m+1,j,m-j}\mbox{\mathversion{bold}$v$}^{\varepsilon}_{s}\times 
(A_{m-j}\partial^{2m+1}_{s}\mbox{\mathversion{bold}$v$}^{\varepsilon})
\big\}
\\[3mm]
& \qquad +
\sum^{m+1}_{j=1}b_{m+1,j}\varepsilon^{m+2-j}
[\mbox{\mathversion{bold}$v$}^{\varepsilon}_{s}\cdot 
(A_{j-1}\partial^{2m+1}_{s}\mbox{\mathversion{bold}$v$}^{\varepsilon})]
\mbox{\mathversion{bold}$v$}^{\varepsilon}
+
\mbox{\mathversion{bold}$U$}^{\varepsilon}_{m+1}
\big|_{s=0,1}.
\end{align*}
Suppose \( m+1 \) is even and set \( m+1=2n \).
The first term on the right-hand side can be calculated as
\begin{align*}
\sum^{m+1}_{j=0}a_{m+1,j}\varepsilon^{j}
A_{m+1-j}\partial^{2m+2}_{s}\mbox{\mathversion{bold}$v$}^{\varepsilon}
&=
\sum^{n}_{k=0}a_{m+1,2k}\varepsilon^{2k}(-1)^{k+1}
\big(
(\mbox{\mathversion{bold}$v$}^{\varepsilon}\cdot 
\partial^{2m+2}_{s}\mbox{\mathversion{bold}$v$}^{\varepsilon})
\mbox{\mathversion{bold}$v$}^{\varepsilon}
-
\partial^{2m+2}_{s}\mbox{\mathversion{bold}$v$}^{\varepsilon}\big)\\[3mm]
& \qquad +
\sum^{n-1}_{k=0}a_{m+1,2k+1}\varepsilon^{2k+1}
(-1)^{k}\mbox{\mathversion{bold}$v$}^{\varepsilon}\times 
\partial^{2m+2}_{s}\mbox{\mathversion{bold}$v$}^{\varepsilon}\\[3mm]
&=
\bigg(
1+\sum^{n}_{k=1}a_{m+1,2k}\varepsilon^{2k}(-1)^{k}
\bigg)
\partial^{2m+2}_{s}\mbox{\mathversion{bold}$v$}^{\varepsilon}\\[3mm]
& \qquad +
\bigg(
\sum^{n-1}_{k=0}a_{m+1,2k+1}\varepsilon^{2k+1}
(-1)^{k}
\bigg)
\mbox{\mathversion{bold}$v$}^{\varepsilon}\times 
\partial^{2m+2}_{s}\mbox{\mathversion{bold}$v$}^{\varepsilon}\\[3mm]
& \qquad +
\sum^{n}_{k=0}a_{m+1,2k}\varepsilon^{2k}(-1)^{k+1}
(\mbox{\mathversion{bold}$v$}^{\varepsilon}\cdot 
\partial^{2m+2}_{s}\mbox{\mathversion{bold}$v$}^{\varepsilon})
\mbox{\mathversion{bold}$v$}^{\varepsilon},
\end{align*}
where we used \( a_{m+1,0}=1 \). Hence we have
\begin{align*}
&\bigg(
1+\sum^{n}_{k=1}a_{m+1,2k}\varepsilon^{2k}(-1)^{k}
\bigg)
\partial^{2m+2}_{s}\mbox{\mathversion{bold}$v$}^{\varepsilon}
+
\bigg(
\sum^{n-1}_{k=0}a_{m+1,2k+1}\varepsilon^{2k+1}
(-1)^{k}
\bigg)
\mbox{\mathversion{bold}$v$}^{\varepsilon}\times 
\partial^{2m+2}_{s}\mbox{\mathversion{bold}$v$}^{\varepsilon}\big|_{s=0,1} \\[3mm]
& \qquad =
-\sum^{n}_{k=0}a_{m+1,2k}\varepsilon^{2k}(-1)^{k+1}
(\mbox{\mathversion{bold}$v$}^{\varepsilon}\cdot 
\partial^{2m+2}_{s}\mbox{\mathversion{bold}$v$}^{\varepsilon})
\mbox{\mathversion{bold}$v$}^{\varepsilon}\\[3mm]
& \qquad \qquad -\sum^{m}_{j=0}\varepsilon^{j}
\big\{
e_{m+1,j,0}A_{m-j}(\mbox{\mathversion{bold}$v$}^{\varepsilon}_{s}\times
\partial^{2m+1}_{s}\mbox{\mathversion{bold}$v$}^{\varepsilon})
+
e_{m+1,j,m-j}\mbox{\mathversion{bold}$v$}^{\varepsilon}_{s}\times 
(A_{m-j}\partial^{2m+1}_{s}\mbox{\mathversion{bold}$v$}^{\varepsilon})
\big\}
\\[3mm]
& \qquad \qquad -
\sum^{m+1}_{j=1}b_{m+1,j}\varepsilon^{m+2-j}
[\mbox{\mathversion{bold}$v$}^{\varepsilon}_{s}\cdot 
(A_{j-1}\partial^{2m+1}_{s}\mbox{\mathversion{bold}$v$}^{\varepsilon})]
\mbox{\mathversion{bold}$v$}^{\varepsilon}
+
\mbox{\mathversion{bold}$U$}^{\varepsilon}_{m+1}
\big|_{s=0,1}.
\end{align*}
Then there exists a \( \varepsilon_{m+1}\in (0,1)\) such that for 
any \( \varepsilon \in (0,\varepsilon_{m+1}) \), the matrix
\begin{align*}
\bigg(
1+\sum^{n}_{k=1}a_{m+1,2k}\varepsilon^{2k}(-1)^{k}
\bigg)I_{3}
-
\varepsilon \bigg(
\sum^{n-1}_{k=0}a_{m+1,2k+1}\varepsilon^{2k}
(-1)^{k+1}
\bigg)A(\mbox{\mathversion{bold}$v$}^{\varepsilon})
\end{align*}
is reversible and the inverse can be expressed as
\begin{align*}
\bigg\{
\bigg(
1+\sum^{n}_{k=1}a_{m+1,2k}\varepsilon^{2k}(-1)^{k}
\bigg)I_{3}
&-
\varepsilon \bigg(
\sum^{n-1}_{k=0}a_{m+1,2k+1}\varepsilon^{2k}
(-1)^{k+1}
\bigg)A(\mbox{\mathversion{bold}$v$}^{\varepsilon})
\bigg\}^{-1}\\[3mm]
& =
\bigg(
1+\sum^{n}_{k=1}a_{m+1,2k}\varepsilon^{2k}(-1)^{k}
\bigg)
(
I_{3}
+
D_{m+1}(\mbox{\mathversion{bold}$v$}^{\varepsilon})
),
\end{align*}
where
\begin{align*}
D_{m+1}(\mbox{\mathversion{bold}$v$}^{\varepsilon})
=
\sum^{\infty}_{j=1}\varepsilon^{j}M^{j}
A(\mbox{\mathversion{bold}$v$}^{\varepsilon})^{j},
\end{align*}
with 
\begin{align*}
M=
\bigg(
\sum^{n-1}_{k=0}a_{m+1,2k+1}\varepsilon^{2k}
(-1)^{k+1}
\bigg)
\bigg(
1+\sum^{n}_{k=1}a_{m+1,2k}\varepsilon^{2k}(-1)^{k}
\bigg)^{-1}.
\end{align*}
Hence we have
\begin{align*}
&\partial^{2m+2}_{s}\mbox{\mathversion{bold}$v$}^{\varepsilon} \big|_{s=0,1}
=
\bigg(
1+\sum^{n}_{k=1}a_{m+1,2k}\varepsilon^{2k}(-1)^{k}
\bigg)
(
I_{3}
+
D_{m+1}(\mbox{\mathversion{bold}$v$}^{\varepsilon})
)\\[3mm]
& \qquad 
\bigg\{
-\sum^{n}_{k=0}a_{m+1,2k}\varepsilon^{2k}(-1)^{k+1}
(\mbox{\mathversion{bold}$v$}^{\varepsilon}\cdot 
\partial^{2m+2}_{s}\mbox{\mathversion{bold}$v$}^{\varepsilon})
\mbox{\mathversion{bold}$v$}^{\varepsilon}\\[3mm]
& \qquad \quad  -\sum^{m}_{j=0}\varepsilon^{j}
\big\{
e_{m+1,j,0}A_{m-j}(\mbox{\mathversion{bold}$v$}^{\varepsilon}_{s}\times
\partial^{2m+1}_{s}\mbox{\mathversion{bold}$v$}^{\varepsilon})
+
e_{m+1,j,m-j}\mbox{\mathversion{bold}$v$}^{\varepsilon}_{s}\times 
(A_{m-j}\partial^{2m+1}_{s}\mbox{\mathversion{bold}$v$}^{\varepsilon})
\big\}
\\[3mm]
& \qquad \quad  -
\sum^{m+1}_{j=1}b_{m+1,j}\varepsilon^{m+2-j}
[\mbox{\mathversion{bold}$v$}^{\varepsilon}_{s}\cdot 
(A_{j-1}\partial^{2m+1}_{s}\mbox{\mathversion{bold}$v$}^{\varepsilon})]
\mbox{\mathversion{bold}$v$}^{\varepsilon}
+
\mbox{\mathversion{bold}$U$}^{\varepsilon}_{m+1}
\bigg\}
\bigg|_{s=0,1}.
\end{align*}
Since
\begin{align*}
D_{m+1}(\mbox{\mathversion{bold}$v$}^{\varepsilon})
&=
\varepsilon 
\bigg(
\sum^{\infty}_{j=0}\varepsilon^{j}M^{j+1}
A(\mbox{\mathversion{bold}$v$}^{\varepsilon})^{j}
\bigg)
A(\mbox{\mathversion{bold}$v$}^{\varepsilon})\\[3mm]
&=:
\varepsilon 
B_{m+1}(\mbox{\mathversion{bold}$v$}^{\varepsilon})
A(\mbox{\mathversion{bold}$v$}^{\varepsilon})
\end{align*}
we have
\begin{align}
\label{mm}
\partial^{2m+2}_{s}\mbox{\mathversion{bold}$v$}^{\varepsilon}|_{s=0,1}
&=
C_{m+1}
\bigg\{
-\sum^{n}_{k=0}a_{m+1,2k}\varepsilon^{2k}(-1)^{k+1}
(\mbox{\mathversion{bold}$v$}^{\varepsilon}\cdot 
\partial^{2m+2}_{s}\mbox{\mathversion{bold}$v$}^{\varepsilon})
\mbox{\mathversion{bold}$v$}^{\varepsilon}\\[3mm]
& \qquad -\sum^{m}_{j=0}\varepsilon^{j}
\big\{
e_{m+1,j,0}A_{m-j}(\mbox{\mathversion{bold}$v$}^{\varepsilon}_{s}\times
\partial^{2m+1}_{s}\mbox{\mathversion{bold}$v$}^{\varepsilon}) \nonumber \\[3mm]
& \qquad \qquad +
e_{m+1,j,m-j}\mbox{\mathversion{bold}$v$}^{\varepsilon}_{s}\times 
(A_{m-j}\partial^{2m+1}_{s}\mbox{\mathversion{bold}$v$}^{\varepsilon})
\big\}\nonumber
\\[3mm]
& \qquad -
\sum^{m+1}_{j=1}b_{m+1,j}\varepsilon^{m+2-j}
[\mbox{\mathversion{bold}$v$}^{\varepsilon}_{s}\cdot 
(A_{j-1}\partial^{2m+1}_{s}\mbox{\mathversion{bold}$v$}^{\varepsilon})]
\mbox{\mathversion{bold}$v$}^{\varepsilon}
\bigg\}\nonumber\\[3mm]
& \qquad +
C_{m+1}
B_{m+1}(\mbox{\mathversion{bold}$v$}^{\varepsilon})
\bigg\{
-\sum^{m}_{j=0}\varepsilon^{j}
\big\{
e_{m+1,j,0}A_{m+1-j}(\mbox{\mathversion{bold}$v$}^{\varepsilon}_{s}\times
\partial^{2m+1}_{s}\mbox{\mathversion{bold}$v$}^{\varepsilon})\bigg\} \nonumber \\[3mm]
& \qquad +
\varepsilon \big( C_{m+1}I_{3}+\varepsilon B_{m+1}
(\mbox{\mathversion{bold}$v$}^{\varepsilon})
\big) 
\mbox{\mathversion{bold}$U$}^{\varepsilon}_{m+1}
\bigg|_{s=0,1},\nonumber
\end{align}
where we have set
\begin{align*}
C_{m+1}=\left( 
1+\sum^{n}_{k=1}a_{m+1,2k}\varepsilon^{2k}(-1)^{k}
\right)^{-1},
\end{align*}
Lastly, since
\begin{align*}
A_{k}(\mbox{\mathversion{bold}$v$}^{\varepsilon}_{s}\times 
\partial^{2m+1}_{s}\mbox{\mathversion{bold}$v$}^{\varepsilon})
=
\left\{
\begin{array}{ll}
(-1)^{l+1}(\mbox{\mathversion{bold}$v$}^{\varepsilon}\cdot \partial^{2m+1}_{s}\mbox{\mathversion{bold}$v$}^{\varepsilon})
\mbox{\mathversion{bold}$v$}^{\varepsilon}_{s}, & k=2l-1, \\[3mm]
(-1)^{l+1}(\mbox{\mathversion{bold}$v$}^{\varepsilon}\cdot \partial^{2m+1}_{s}\mbox{\mathversion{bold}$v$}^{\varepsilon})
\mbox{\mathversion{bold}$v$}^{\varepsilon}\times \mbox{\mathversion{bold}$v$}^{\varepsilon}_{s}, & k=2l,
\end{array}\right.
\end{align*}
This shows that combining (\ref{mm}) and (\ref{one}), we have the estimate
\begin{align*}
\varepsilon 
\big| [\partial^{2m+1}_{s}\mbox{\mathversion{bold}$v$}^{\varepsilon}
\cdot
\partial^{2m+2}_{s}\mbox{\mathversion{bold}$v$}^{\varepsilon}]^{1}_{s=0}
\big|
&\leq
\frac{\varepsilon}{16}\|\partial^{2m+2}_{s}
\mbox{\mathversion{bold}$v$}^{\varepsilon}\|^{2}
+
C\|\mbox{\mathversion{bold}$v$}^{\varepsilon}_{s}\|_{2m}^{2},
\end{align*}
holds, where \( C>0\) depends on \( \| \mbox{\mathversion{bold}$v$}^{\varepsilon}_{s}\|_{2(m-1)}\)
and we also utilized the fact that 
\begin{align*}
\| B_{m+1}(\mbox{\mathversion{bold}$v$}^{\varepsilon})\|_{L^{\infty}(I\times [0,T])}
\leq C
\end{align*}
with \( C>0\) independent of \( \varepsilon \).

Next, we consider 
\( [
\partial^{k}_{s}\mbox{\mathversion{bold}$v$}^{\varepsilon}\cdot 
(\mbox{\mathversion{bold}$v$}^{\varepsilon}\times 
\partial^{k+1}_{s}\mbox{\mathversion{bold}$v$}^{\varepsilon})]^{1}_{s=0}
\).
We further calculate the right-hand side of (\ref{vc2m}).
By Lemma \ref{lm2m2m-1} with \( m=l\) for 
\( 2\leq l \leq m\), we have 
\begin{align*}
\partial^{2(m+1-l)}_{s}
\partial^{l}_{t}\mbox{\mathversion{bold}$v$}^{\varepsilon}
&=
\sum^{l}_{j=0}a_{l,j}\varepsilon^{j}A_{l-j}
\partial^{2m+2}_{s}\mbox{\mathversion{bold}$v$}^{\varepsilon}
+
2m\sum^{l-1}_{j=0}a_{l,j}\varepsilon^{j}
\sum^{l-1-j}_{i=0}A_{l-1-j-i}
\mbox{\mathversion{bold}$v$}^{\varepsilon}_{s}\times 
(A_{i}\partial^{2m+1}_{s}\mbox{\mathversion{bold}$v$}^{\varepsilon})\\[3mm]
& \qquad +
\sum^{l-1}_{j=0}\sum^{l-1-j}_{k=0}
e_{l,j,k}\varepsilon^{j}A_{l-1-j-k}
\mbox{\mathversion{bold}$v$}^{\varepsilon}_{s}\times 
(A_{k}\partial^{2m+1}_{s}\mbox{\mathversion{bold}$v$}^{\varepsilon})\\[3mm]
& \qquad +
\sum^{l}_{j=1}b_{l,j}\varepsilon^{l+1-j}
[\mbox{\mathversion{bold}$v$}^{\varepsilon}_{s}\cdot 
(A_{j-1}\partial^{2m+1}_{s}\mbox{\mathversion{bold}$v$}^{\varepsilon})]
\mbox{\mathversion{bold}$v$}^{\varepsilon}
+
\mbox{\mathversion{bold}$U$}^{\varepsilon}_{m+1}\\[3mm]
\partial^{m}_{t}\mbox{\mathversion{bold}$v$}^{\varepsilon}_{s}
&=
\sum^{m}_{j=0}a_{m,j}\varepsilon ^{j}A_{m-j}\partial^{2m+1}_{s}
\mbox{\mathversion{bold}$v$}^{\varepsilon}
+
\mbox{\mathversion{bold}$U$}^{\varepsilon}_{m+1},
\end{align*}
where as in the proof of Lemma \ref{vc2m}, \( \mbox{\mathversion{bold}$U$}^{\varepsilon}_{m}\) 
denotes the collection of terms that satisfy the estimate stated in Lemma \ref{2m2m-1}, and may 
change from line to line. 
This yields
\begin{align*}
\partial^{2(m+1-l)}_{s}\partial^{l}_{t}\mbox{\mathversion{bold}$v$}^{\varepsilon}
&=
\sum^{l}_{j=0}a_{l,j}\varepsilon^{j}A_{l-j}\partial^{2m+2}_{s}
\mbox{\mathversion{bold}$v$}^{\varepsilon}\\[3mm]
& \qquad +
2m\sum^{l-1}_{j=0}a_{l,j}\varepsilon^{j}
\big(
A_{l-1-j}(\mbox{\mathversion{bold}$v$}^{\varepsilon}_{s}\times 
\partial^{2m+1}_{s}\mbox{\mathversion{bold}$v$}^{\varepsilon})
+
\mbox{\mathversion{bold}$v$}^{\varepsilon}_{s}\times 
(A_{l-1-j}\partial^{2m+1}_{s}\mbox{\mathversion{bold}$v$}^{\varepsilon})
\big)\\[3mm]
& \qquad 
+
\sum^{l-1}_{j=0}\varepsilon^{j}
\big(
e_{l,j,0}
A_{l-1-j}(\mbox{\mathversion{bold}$v$}^{\varepsilon}_{s}\times 
\partial^{2m+1}_{s}\mbox{\mathversion{bold}$v$}^{\varepsilon})
+
e_{l,j,l-1-j}\mbox{\mathversion{bold}$v$}^{\varepsilon}_{s}\times 
(A_{l-1-j}\partial^{2m+1}_{s}\mbox{\mathversion{bold}$v$}^{\varepsilon})
\big)\\[3mm]
& \qquad +
\sum^{l}_{j=1}b_{l,j}\varepsilon^{l+1-j}
[\mbox{\mathversion{bold}$v$}^{\varepsilon}_{s}\cdot 
(A_{j-1}\partial^{2m+1}_{s}\mbox{\mathversion{bold}$v$}^{\varepsilon})]
\mbox{\mathversion{bold}$v$}^{\varepsilon}
+
\mbox{\mathversion{bold}$U$}^{\varepsilon}_{m+1}.
\end{align*}
Substituting these expressions for 
\( \partial^{2(m+1-l)}_{s}\partial^{l}_{t}\mbox{\mathversion{bold}$v$}^{\varepsilon} \)
and 
\( \partial^{m}_{t}\mbox{\mathversion{bold}$v$}^{\varepsilon}_{s}\)
into the 
right-hand side of (\ref{vc2m}) with \( m \) as \( m+1 \) yields
\begin{align}
\label{mmm}
\mbox{\mathversion{bold}$v$}^{\varepsilon}\times 
\partial^{2m+2}_{s}
\mbox{\mathversion{bold}$v$}^{\varepsilon}|_{s=0,1}
& =
(-1)^{m+1}\varepsilon A_{m}
\bigg\{
\sum^{m}_{j=0}a_{m,j}\varepsilon ^{j}\partial^{2m+2}_{s}\mbox{\mathversion{bold}$v$}^{\varepsilon} \\[3mm]
& \qquad +
2m\sum^{m-1}_{j=0}a_{m,j}\varepsilon^{j}
\big(
A_{m-1-j}(\mbox{\mathversion{bold}$v$}^{\varepsilon}_{s}\times 
\partial^{2m+1}_{s}\mbox{\mathversion{bold}$v$}^{\varepsilon})
+
\mbox{\mathversion{bold}$v$}^{\varepsilon}_{s}\times
(A_{m-1-j}\partial^{2m+1}_{s}\mbox{\mathversion{bold}$v$}^{\varepsilon})
\big)\nonumber\\[3mm]
& \qquad +
\sum^{m-1}_{j=0}\varepsilon^{j}
\big\{
e_{m,j,0}A_{m-1-j}(\mbox{\mathversion{bold}$v$}^{\varepsilon}_{s}\times 
\partial^{2m+1}_{s}\mbox{\mathversion{bold}$v$}^{\varepsilon})\nonumber\\[3mm]
& \qquad \qquad \qquad +
e_{m,j,m-1-j}\mbox{\mathversion{bold}$v$}^{\varepsilon}_{s}\times 
(A_{m-1-j}\partial^{2m+1}_{s}\mbox{\mathversion{bold}$v$}^{\varepsilon})
\big\}\nonumber \\[3mm]
& \qquad +
\sum^{m}_{j=1}b_{m,j}\varepsilon^{m+1-j}
[
\mbox{\mathversion{bold}$v$}^{\varepsilon}_{s}\cdot (
A_{j-1}\partial^{2m+1}_{s}\mbox{\mathversion{bold}$v$}^{\varepsilon})]
\mbox{\mathversion{bold}$v$}^{\varepsilon}\nonumber\\[3mm]
& \qquad +
2(
\mbox{\mathversion{bold}$v$}^{\varepsilon}_{s}\cdot 
\big[
\sum^{m}_{j=0}\varepsilon^{j}A_{m-j}
\partial^{2m+1}_{s}\mbox{\mathversion{bold}$v$}^{\varepsilon}\big])
\mbox{\mathversion{bold}$v$}^{\varepsilon}
\bigg\}\nonumber \\[3mm]
& \qquad + \varepsilon 
\bigg\{
\sum^{m-1}_{l=0}(-1)^{l}A_{l+2}
\big[
\sum^{l}_{j=0}a_{l,j}\varepsilon^{j}A_{l-j}\partial^{2m+2}_{s}
\mbox{\mathversion{bold}$v$}^{\varepsilon}\nonumber\\[3mm]
& \qquad +
2m\sum^{l-1}_{j=0}a_{l,j}\varepsilon^{j}
\big(
A_{l-1-j}(\mbox{\mathversion{bold}$v$}^{\varepsilon}_{s}\times 
\partial^{2m+1}_{s}\mbox{\mathversion{bold}$v$}^{\varepsilon})
+
\mbox{\mathversion{bold}$v$}^{\varepsilon}_{s}\times 
(A_{l-1-j}\partial^{2m+1}_{s}\mbox{\mathversion{bold}$v$}^{\varepsilon})
\big)\nonumber\\[3mm]
& \qquad 
+
\sum^{l-1}_{j=0}\varepsilon^{j}
\big[
e_{l,j,0}
A_{l-1-j}(\mbox{\mathversion{bold}$v$}^{\varepsilon}_{s}\times 
\partial^{2m+1}_{s}\mbox{\mathversion{bold}$v$}^{\varepsilon})\nonumber\\[3mm]
& \qquad \qquad \qquad +
e_{l,j,l-1-j}\mbox{\mathversion{bold}$v$}^{\varepsilon}_{s}\times 
(A_{l-1-j}\partial^{2m+1}_{s}\mbox{\mathversion{bold}$v$}^{\varepsilon})
\big]\nonumber\\[3mm]
& \qquad +
\sum^{l}_{j=1}b_{l,j}\varepsilon^{l+1-j}
[\mbox{\mathversion{bold}$v$}^{\varepsilon}_{s}\cdot 
(A_{j-1}\partial^{2m+1}_{s}\mbox{\mathversion{bold}$v$}^{\varepsilon})]
\mbox{\mathversion{bold}$v$}^{\varepsilon}
\big]\bigg\}
+
\varepsilon \mbox{\mathversion{bold}$U$}^{\varepsilon}_{m+1}\big|_{s=0,1}.\nonumber
\end{align}
Substituting (\ref{mm}) in the right-hand side, re-writing terms of the form 
\( A_{k}\mbox{\mathversion{bold}$W$}\) according to \( k\), and
utilizing (\ref{one}) and (\ref{oddinner}), 
we can see that the following estimates hold.
\begin{align*}
\big|
[
\partial^{2m+1}_{s}\mbox{\mathversion{bold}$v$}^{\varepsilon}\cdot 
(\mbox{\mathversion{bold}$v$}^{\varepsilon}\times 
\partial^{2m+2}_{s}\mbox{\mathversion{bold}$v$}^{\varepsilon})]^{1}_{s=0}
\big|
&\leq 
\frac{\varepsilon}{16}\|\partial^{2m+2}_{s}\mbox{\mathversion{bold}$v$}^{\varepsilon}\|^{2}
+
C\| \mbox{\mathversion{bold}$v$}^{\varepsilon}_{s}\|^{2}_{2m}\\[3mm]
\big|
[
\partial^{2m}_{s}\mbox{\mathversion{bold}$v$}^{\varepsilon}\cdot 
(\mbox{\mathversion{bold}$v$}^{\varepsilon}\times 
\partial^{2m+1}_{s}\mbox{\mathversion{bold}$v$}^{\varepsilon})]^{1}_{s=0}
\big|
&=
\big|
[
\partial^{2m+1}_{s}\mbox{\mathversion{bold}$v$}^{\varepsilon}\cdot 
(\mbox{\mathversion{bold}$v$}^{\varepsilon}\times 
\partial^{2m}_{s}\mbox{\mathversion{bold}$v$}^{\varepsilon})]^{1}_{s=0}
\big| \\[3mm]
& \leq 
\frac{\varepsilon }{16}
\| \partial^{2m+2}_{s}\mbox{\mathversion{bold}$v$}^{\varepsilon}\|^{2}
+C\|\mbox{\mathversion{bold}$v$}^{\varepsilon}_{s}\|^{2}_{2m}
\end{align*}
where \( C>0\) depends on \( \| \mbox{\mathversion{bold}$v$}^{\varepsilon}_{s}\|_{2(m-1)}\), but 
not on \( \varepsilon \). We also mention that the above estimate holds because 
every term on the right-hand side of (\ref{mmm}) contains \( \varepsilon \).
Finally from (\ref{oddinner}), we have
\begin{align*}
\big|
\big[
(\partial^{2m+1-j}_{s}\mbox{\mathversion{bold}$v$}^{\varepsilon}\cdot 
\partial^{j}_{s}\mbox{\mathversion{bold}$v$}^{\varepsilon})
[
(\mbox{\mathversion{bold}$v$}^{\varepsilon}\times 
\mbox{\mathversion{bold}$v$}^{\varepsilon}_{s})\cdot 
\partial^{2m+1}_{s}\mbox{\mathversion{bold}$v$}^{\varepsilon}]
\big]^{1}_{s=0}
\big|
&\leq 
2\varepsilon |Y_{2m+1-j,j}| |\mbox{\mathversion{bold}$v$}^{\varepsilon}_{s}|
|\partial^{2m+1}_{s}\mbox{\mathversion{bold}$v$}^{\varepsilon}| \big|_{s=0,1}\\[3mm]
& \leq 
\frac{\varepsilon }{16}\| \partial^{2m+2}_{s}\mbox{\mathversion{bold}$v$}^{\varepsilon}\|^{2}
+
C\| \mbox{\mathversion{bold}$v$}^{\varepsilon}_{s}\|^{2}_{2m}.
\end{align*}

Combining all of the estimates obtained for the boundary terms 
with the estimate we have derived so far for \( \partial^{k_{0}}_{s}\mbox{\mathversion{bold}$v$}^{\varepsilon}\)
yields
\begin{align*}
\frac{1}{2}\frac{{\rm d}}{{\rm d}t}\|\partial^{2m}_{s}\mbox{\mathversion{bold}$v$}^{\varepsilon}\|^{2}_{1}
\leq 
C\|\mbox{\mathversion{bold}$v$}^{\varepsilon}_{s}\|^{2}_{2m} -\frac{\varepsilon }{8}
\| \partial^{2m+1}_{s}\mbox{\mathversion{bold}$v$}^{\varepsilon}\|^{2}_{1},
\end{align*}
where \( C>0\) depends on \( \| \mbox{\mathversion{bold}$v$}^{\varepsilon}_{s}\|_{2(m-1)}\).
If we set \( \varepsilon _{0}:= \min \{ \varepsilon _{l} \ | \ 1\leq l \leq m+1 \} \),
then from Gronwall's inequality and the assumption of induction, we see that there exists \( c_{\ast}>0 \) depending on 
\( \| \mbox{\mathversion{bold}$v$}_{0s}\|_{2m}\) and \( T_{0}\), and  \( T_{0}\) depending on 
\( \| \mbox{\mathversion{bold}$v$}_{0s}\|_{2}\) such that 
for any \( \varepsilon \in (0,\varepsilon _{0}] \),
\begin{align*}
\sup_{0\leq t\leq T_{0}}\| \mbox{\mathversion{bold}$v$}^{\varepsilon}_{s}(t)\|_{2m} \leq c_{\ast}.
\end{align*}
Estimating the \( t\) derivatives of \( \mbox{\mathversion{bold}$v$}^{\varepsilon}\) 
via the equation in (\ref{rnl}), we obtain the necessary estimates of \( \mbox{\mathversion{bold}$v$}^{\varepsilon}\) in
\( Y^{m}_{T_{0}}(I)\) and this finishes the proof of Proposition \ref{ue}. \hfill \( \Box .\)

\subsection{Taking the limit \( \varepsilon \to +0 \)}

Now we take the limit \( \varepsilon \to +0 \). 
Fix \( m\geq 1\) and for \( \varepsilon \in (0,\varepsilon _{0}]\), let
\( \mbox{\mathversion{bold}$v$}^{\varepsilon }\in Y^{m+2}_{T_{0}}(I)\) be the solution of 
(\ref{rnl}) satisfying the uniform estimate given in Proposition \ref{ue}.
Let
\( \varepsilon,\varepsilon' \in (0,\varepsilon _{0}] \),
\( \varepsilon < \varepsilon '\), and set 
\( \mbox{\mathversion{bold}$Z$}:= \mbox{\mathversion{bold}$v$}^{\varepsilon '}
- \mbox{\mathversion{bold}$v$}^{\varepsilon}\). We see that 
\( \mbox{\mathversion{bold}$Z$}\) satisfies
\begin{align*}
\left\{
\begin{array}{ll}
\mbox{\mathversion{bold}$Z$}_{t}=
\mbox{\mathversion{bold}$v$}^{\varepsilon}\times \mbox{\mathversion{bold}$Z$}_{ss}
+
\mbox{\mathversion{bold}$Z$}\times \mbox{\mathversion{bold}$v$}^{\varepsilon '}_{ss}
+
\varepsilon ' \mbox{\mathversion{bold}$Z$}_{ss}
+
(\varepsilon ' - \varepsilon )\mbox{\mathversion{bold}$v$}^{\varepsilon}_{ss}
+
\varepsilon ' |\mbox{\mathversion{bold}$v$}^{\varepsilon '}_{s}|^{2}
\mbox{\mathversion{bold}$v$}^{\varepsilon '}
-
\varepsilon |\mbox{\mathversion{bold}$v$}^{\varepsilon}_{s}|^{2}
\mbox{\mathversion{bold}$v$}^{\varepsilon}, 
& s>0, \ t>0, \\[3mm]
\mbox{\mathversion{bold}$Z$}(s,0)= \mbox{\mathversion{bold}$0$}, & s>0, \\[3mm]
\mbox{\mathversion{bold}$Z$}(0,t)=\mbox{\mathversion{bold}$Z$}(1,t)
=\mbox{\mathversion{bold}$0$}, & t>0.
\end{array}\right.
\end{align*}
We have
\begin{align*}
\frac{1}{2}\frac{{\rm d}}{{\rm d}t}
\| \mbox{\mathversion{bold}$Z$}\|^{2}
=
(\mbox{\mathversion{bold}$Z$},\mbox{\mathversion{bold}$Z$}_{t})
&\leq 
(\mbox{\mathversion{bold}$Z$},\mbox{\mathversion{bold}$v$}^{\varepsilon }\times 
\mbox{\mathversion{bold}$Z$}_{ss})
+\varepsilon '( \mbox{\mathversion{bold}$Z$},\mbox{\mathversion{bold}$Z$}_{ss})
+
(\varepsilon ' +\varepsilon )c_{\ast}\\[3mm]
& \leq
c_{\ast}
\| \mbox{\mathversion{bold}$Z$}\|^{2}_{1} 
-\varepsilon ' \| \mbox{\mathversion{bold}$Z$}_{s}\|^{2}
+
(\varepsilon ' + \varepsilon )c_{\ast}
\end{align*}
\begin{align*}
\frac{1}{2}\frac{{\rm d}}{{\rm d}t}
\|\mbox{\mathversion{bold}$Z$}_{s}\|^{2}
=
(\mbox{\mathversion{bold}$Z$}_{s},\mbox{\mathversion{bold}$Z$}_{ts})
& =
-(\mbox{\mathversion{bold}$Z$}_{ss},\mbox{\mathversion{bold}$Z$}_{t})\\[3mm]
& =
-(\mbox{\mathversion{bold}$Z$}_{ss},\mbox{\mathversion{bold}$Z$}\times 
\mbox{\mathversion{bold}$v$}^{\varepsilon '}_{ss})
-\varepsilon ' \|\mbox{\mathversion{bold}$Z$}_{ss}\|^{2}
-(\varepsilon '-\varepsilon )(\mbox{\mathversion{bold}$Z$}_{ss},
\mbox{\mathversion{bold}$v$}^{\varepsilon }_{ss})\\[3mm]
& \qquad -
\varepsilon ' (\mbox{\mathversion{bold}$Z$}_{ss},
|\mbox{\mathversion{bold}$v$}^{\varepsilon '}_{s}|\mbox{\mathversion{bold}$v$}^{\varepsilon '})
+
\varepsilon 
(\mbox{\mathversion{bold}$Z$}_{ss},|\mbox{\mathversion{bold}$v$}^{\varepsilon }_{s}|^{2}
\mbox{\mathversion{bold}$v$}^{\varepsilon })\\[3mm]
& \leq
c_{\ast}\| \mbox{\mathversion{bold}$Z$}\|^{2}_{1}
-\frac{\varepsilon '}{2}\| \mbox{\mathversion{bold}$Z$}_{ss}\|^{2}
+c_{\ast}(\varepsilon ' + \varepsilon ),
\end{align*}
where \( c_{\ast}>0\) depends on the uniform estimate of 
\( \mbox{\mathversion{bold}$v$}^{\varepsilon }\) in \( Y^{1}_{T_{0}}(I)\). 
Hence by Gronwall's inequality, we see that
\begin{align*}
\sup_{0\leq t \leq T_{0}}\|\mbox{\mathversion{bold}$Z$}(t)\|^{2}_{1}
\leq
c_{\ast}e^{c_{\ast}T_{0}}(\varepsilon ' + \varepsilon ),
\end{align*}
which proves that 
\( \mbox{\mathversion{bold}$Z$}\to \mbox{\mathversion{bold}$0$}\) in 
\( C\big( [0,T_{0}];H^{1}(I) \big)\). Hence 
there exists \( \mbox{\mathversion{bold}$v$}\in 
C\big( [0,T_{0}];H^{1}(I) \big) \) such that
\( \mbox{\mathversion{bold}$v$}^{\varepsilon }\to \mbox{\mathversion{bold}$v$}\).
The uniform estimate and the interpolation inequality implies that
\( \mbox{\mathversion{bold}$v$}^{\varepsilon }\to \mbox{\mathversion{bold}$v$} \) in
\( \bigcap ^{m}_{j=0}C^{j}\big( [0,T_{0}];H^{2(m-j)}(I)\big) \), and 
this 
\( \mbox{\mathversion{bold}$v$} \)
is the desired solution of (\ref{vslant}). 
Note that from the uniform estimate, \( 
\partial^{j}_{t}\mbox{\mathversion{bold}$v$}\in 
\bigcap ^{m}_{j=0}L^{\infty}\big(0,T_{0};H^{2(m-j)+1}(I)\big)\) 
can be shown from an argument utilizing 
weak\( \ast \) convergence, but this does not benefit us because the 
energy estimate for the solution of the original problem (\ref{vslant}) 
is derived in Sobolev spaces with different indices, which will be shown in the 
next subsection.
\subsection{Energy estimate of the solution to (\ref{vslant})}
Up until now, we have assumed that the solution is smooth as we need it to be,
which is possible because the intial datum can be 
approximated by smooth ones from Propostion \ref{id1} and \ref{id2}.
We derive an apriori estimate of the solution to (\ref{vslant}) to 
justify the approximation argument. We emphasize this point because 
the estimate obtained here and the uniform estimate obtained in the last subsection 
are derived in Sobolev spaces with different indices.
Fix an arbitrary non-negative integer \( l\). We derive the estimate in the 
class \( H^{l+3}(I)\).
For any \( T>0\), let \( \mbox{\mathversion{bold}$v$}\in \bigcap ^{m+2}_{j=0}C^{j}
\big( [0,T];H^{2(m+2-j)}(I)\big)\)
be the solution of (\ref{vslant}) obtained in the last 
subsection such that \( l+5 \leq 2m+2\). 
We see that
\begin{align*}
\frac{1}{2}\frac{{\rm d}}{{\rm d}t}|\mbox{\mathversion{bold}$v$}|^{2}
=
\mbox{\mathversion{bold}$v$}\cdot \mbox{\mathversion{bold}$v$}_{t}
\equiv 0,
\end{align*}
which shows that \( | \mbox{\mathversion{bold}$v$}|\equiv 1 \).
Hence we can prove that (\ref{decom}) and (\ref{one}) also holds even if we replace 
\( \mbox{\mathversion{bold}$v$}^{\varepsilon }\) with 
\( \mbox{\mathversion{bold}$v$}\).
Next, we estimate conserved quantities which were utilized in the analysis of 
the Cauchy problem in Nishiyama and Tani \cite{5}.
\begin{align*}
\frac{1}{2}\frac{{\rm d}}{{\rm d}t}\|\mbox{\mathversion{bold}$v$}_{s}\|^{2}
&=
-(\mbox{\mathversion{bold}$v$}_{ss},\mbox{\mathversion{bold}$v$}_{t})
+[\mbox{\mathversion{bold}$v$}_{s}\cdot \mbox{\mathversion{bold}$v$}_{t}]^{1}_{s=0}
=0\\[3mm]
\frac{{\rm d}}{{\rm d}t}\bigg\{
\| \mbox{\mathversion{bold}$v$}_{ss}\|^{2}-\frac{5}{4}
\big\| |\mbox{\mathversion{bold}$v$}_{s}|^{2}\big\|^{2}
\bigg\}
&=
2(\mbox{\mathversion{bold}$v$}_{ss},\mbox{\mathversion{bold}$v$}_{tss})
-5\int^{1}_{0}|\mbox{\mathversion{bold}$v$}_{s}|^{2}
\mbox{\mathversion{bold}$v$}_{s}\cdot \mbox{\mathversion{bold}$v$}_{ts} {\rm d}s \\[3mm]
&= -2(\mbox{\mathversion{bold}$v$}_{sss},\mbox{\mathversion{bold}$v$}_{ts})
-5\int^{1}_{0}|\mbox{\mathversion{bold}$v$}_{s}|^{2}
\mbox{\mathversion{bold}$v$}_{s}\cdot \mbox{\mathversion{bold}$v$}_{ts} {\rm d}s \\[3mm]
& \qquad +
2[\mbox{\mathversion{bold}$v$}_{ss}\cdot \mbox{\mathversion{bold}$v$}_{ts}]^{1}_{s=0}\\[3mm]
&=
\big[
-3|\mbox{\mathversion{bold}$v$}_{s}|^{2}\mbox{\mathversion{bold}$v$}_{ss}\cdot 
(\mbox{\mathversion{bold}$v$}\times \mbox{\mathversion{bold}$v$}_{s})\big]^{1}_{s=0}
+
2[\mbox{\mathversion{bold}$v$}_{ss}\cdot \mbox{\mathversion{bold}$v$}_{ts}]^{1}_{s=0}\\[3mm]
& =
\big[
3|\mbox{\mathversion{bold}$v$}_{s}|^{2}\mbox{\mathversion{bold}$v$}_{s}\cdot 
(\mbox{\mathversion{bold}$v$}\times \mbox{\mathversion{bold}$v$}_{ss})\big]^{1}_{s=0}
+
2[\mbox{\mathversion{bold}$v$}_{ss}\cdot \mbox{\mathversion{bold}$v$}_{ts}]^{1}_{s=0}
\end{align*}
From \( \mbox{\mathversion{bold}$v$}_{t}|_{s=0,1}=
\mbox{\mathversion{bold}$0$}\), we have
\begin{align*}
\mbox{\mathversion{bold}$v$}\times \mbox{\mathversion{bold}$v$}_{ss}|_{s=0,1}
=
\mbox{\mathversion{bold}$0$},
\end{align*}
which shows that the first boundary term is zero. Direct calculation also yields
\begin{align*}
\mbox{\mathversion{bold}$v$}_{ss}\cdot \mbox{\mathversion{bold}$v$}_{ts}|_{s=0,1}
=
\mbox{\mathversion{bold}$v$}_{ss}\cdot (
\mbox{\mathversion{bold}$v$}\times \mbox{\mathversion{bold}$v$}_{sss}
+
\mbox{\mathversion{bold}$v$}_{s}\times \mbox{\mathversion{bold}$v$}_{ss}
)|_{s=0,1}
&=
\mbox{\mathversion{bold}$v$}_{ss}\cdot 
(\mbox{\mathversion{bold}$v$}\times \mbox{\mathversion{bold}$v$}_{sss})\\[3mm]
&=
-\mbox{\mathversion{bold}$v$}_{sss}\cdot (
\mbox{\mathversion{bold}$v$}\times \mbox{\mathversion{bold}$v$}_{ss})|_{s=0,1}\\[3mm]
&= \mbox{\mathversion{bold}$0$},
\end{align*}
which proves
\begin{align*}
\frac{{\rm d}}{{\rm d}t}\bigg\{
\| \mbox{\mathversion{bold}$v$}_{ss}\|^{2}-\frac{5}{4}
\big\| |\mbox{\mathversion{bold}$v$}_{s}|^{2}\big\|^{2}
\bigg\}
=
0.
\end{align*}
Similarly, we can prove that
\begin{align*}
\frac{{\rm d}}{{\rm d}t}\bigg\{
\|\mbox{\mathversion{bold}$v$}_{sss}\|^{2}
&-\frac{7}{2}\big\| |\mbox{\mathversion{bold}$v$}_{s}|
|\mbox{\mathversion{bold}$v$}_{ss}|\big\|^{2}-
14\| \mbox{\mathversion{bold}$v$}_{s}\cdot \mbox{\mathversion{bold}$v$}_{ss}\|^{2}
+
\frac{21}{8}\big\|
|\mbox{\mathversion{bold}$v$}_{s}|^{3}\big\| ^{2}
\bigg\}\\[3mm]
&=
2[\mbox{\mathversion{bold}$v$}_{sss}\cdot \mbox{\mathversion{bold}$v$}_{tss}]^{1}_{s=0}\\[3mm]
& \qquad +
\bigg[
18(\mbox{\mathversion{bold}$v$}_{s}\cdot \mbox{\mathversion{bold}$v$}_{ss})
\mbox{\mathversion{bold}$v$}_{sss}\cdot (\mbox{\mathversion{bold}$v$}\times 
\mbox{\mathversion{bold}$v$}_{s}) 
+
5|\mbox{\mathversion{bold}$v$}_{s}|^{2}\mbox{\mathversion{bold}$v$}_{sss}\cdot 
(\mbox{\mathversion{bold}$v$}\times \mbox{\mathversion{bold}$v$}_{ss})\\[3mm]
& \qquad -
2(\mbox{\mathversion{bold}$v$}_{s}\cdot \mbox{\mathversion{bold}$v$}_{sss})
(\mbox{\mathversion{bold}$v$}\times \mbox{\mathversion{bold}$v$}_{s})\cdot 
\mbox{\mathversion{bold}$v$}_{ss}
+
|\mbox{\mathversion{bold}$v$}_{ss}|^{2}\mbox{\mathversion{bold}$v$}_{ss}\cdot 
(\mbox{\mathversion{bold}$v$}\times \mbox{\mathversion{bold}$v$}_{s})\\[3mm]
& \qquad -
\frac{27}{4}|\mbox{\mathversion{bold}$v$}_{s}|^{2}
(\mbox{\mathversion{bold}$v$}\times \mbox{\mathversion{bold}$v$}_{s})\cdot 
\mbox{\mathversion{bold}$v$}_{ss}
\bigg]^{1}_{s=0}\\[3mm]
&=
0.
\end{align*}
From these estimates and the interpolation inequality, we see that 
\begin{align*}
\sup_{0\leq t\leq T}\| \mbox{\mathversion{bold}$v$}_{s}(t)\|_{2}
\leq
C\|\mbox{\mathversion{bold}$v$}_{0s}\|^{2}_{2}
(1+\| \mbox{\mathversion{bold}$v$}_{0s}\|^{2}_{2})^{4},
\end{align*}
where \( C>0\) is monotone increasing with respect to \( T>0\).
We continue by induction. Suppose that
\( \mbox{\mathversion{bold}$v$}\) satisfies
\begin{align*}
\sup_{0\leq t\leq T}
\|
\mbox{\mathversion{bold}$v$}_{s}(t)
\|_{k-1}
\leq C_{\ast},
\end{align*}
where \( C_{\ast}>0 \) depends on \( \| \mbox{\mathversion{bold}$v$}_{0s}\|_{k-1}\) 
and is monotone increasing in \( T>0\) for some \( k\geq 3 \).
We calculate
\begin{align*}
\frac{1}{2}\frac{{\rm d}}{{\rm d}t}
\| \partial^{k+1}_{s}\mbox{\mathversion{bold}$v$}\|^{2}
&=
(\partial^{k+1}_{s}\mbox{\mathversion{bold}$v$},
\partial^{k+1}_{s}\mbox{\mathversion{bold}$v$}_{t})
=
(\partial^{k+1}_{s}\mbox{\mathversion{bold}$v$},\mbox{\mathversion{bold}$v$}\times 
\partial^{k+3}_{s}\mbox{\mathversion{bold}$v$})
+
k( \partial^{k+1}_{s}\mbox{\mathversion{bold}$v$}, \mbox{\mathversion{bold}$v$}_{s}\times
\partial^{k+2}_{s}\mbox{\mathversion{bold}$v$})\\[3mm]
& \qquad +
(\partial^{k+1}_{s}\mbox{\mathversion{bold}$v$}, \mbox{\mathversion{bold}$h$}_{k+1}),
\end{align*}
where \( \mbox{\mathversion{bold}$h$}_{k+1}\) are terms that satisfy
\begin{align*}
\|
\mbox{\mathversion{bold}$h$}_{k+1}
\|
\leq
C\| \mbox{\mathversion{bold}$v$}_{s}\|_{k},
\end{align*}
with \( C>0\) depending on \( C_{\ast}\).
We further calculate
\begin{align*}
\frac{1}{2}\frac{{\rm d}}{{\rm d}t}
\| \partial^{k+1}_{s}\mbox{\mathversion{bold}$v$}\|^{2}
&=
-k(\partial^{k+1}_{s}\mbox{\mathversion{bold}$v$},\mbox{\mathversion{bold}$v$}_{s}\times 
\partial^{k+2}_{s}\mbox{\mathversion{bold}$v$})
+
[\partial^{k+1}_{s}\mbox{\mathversion{bold}$v$}\cdot (\mbox{\mathversion{bold}$v$}\times 
\partial^{k+2}_{s}\mbox{\mathversion{bold}$v$})]^{1}_{s=0}
+
(\partial^{k+1}_{s}\mbox{\mathversion{bold}$v$},\mbox{\mathversion{bold}$h$}_{k+1})\\[3mm]
&=
-k(\partial^{k+1}_{s}\mbox{\mathversion{bold}$v$},
(\mbox{\mathversion{bold}$v$}\cdot \partial^{k+2}_{s}\mbox{\mathversion{bold}$v$})
\mbox{\mathversion{bold}$v$}\times \mbox{\mathversion{bold}$v$}_{s})
+
k(\partial^{k+1}_{s}\mbox{\mathversion{bold}$v$},
[(\mbox{\mathversion{bold}$v$}\times \mbox{\mathversion{bold}$v$}_{s})\cdot
\partial^{k+2}_{s}\mbox{\mathversion{bold}$v$}]\mbox{\mathversion{bold}$v$})\\[3mm]
& \qquad +
[\partial^{k+1}_{s}\mbox{\mathversion{bold}$v$}\cdot 
(\mbox{\mathversion{bold}$v$}\times \partial^{k+2}_{s}\mbox{\mathversion{bold}$v$})]^{1}_{s=0}
+
(\partial^{k+1}_{s}\mbox{\mathversion{bold}$v$},\mbox{\mathversion{bold}$h$}_{k+1}),
\end{align*}
where (\ref{decom}) was used in the last equality.
Furthermore, from (\ref{one}) we have
\begin{align*}
\frac{1}{2}\frac{{\rm d}}{{\rm d}t}
\| \partial^{k+1}_{s}\mbox{\mathversion{bold}$v$}\|^{2}
&=
\frac{k}{2}\sum^{k+1}_{j=1}
\left(
\begin{array}{c}
k+2\\
j
\end{array}
\right)
(\partial^{k+1}_{s}\mbox{\mathversion{bold}$v$},
( \partial^{j}_{s}\mbox{\mathversion{bold}$v$}\cdot 
\partial^{k+2-j}_{s}\mbox{\mathversion{bold}$v$})\mbox{\mathversion{bold}$v$}\times 
\mbox{\mathversion{bold}$v$}_{s})\\[3mm]
& \qquad -
\frac{k}{2}\sum^{k}_{j=1}
\left(
\begin{array}{c}
k+1\\
j
\end{array}
\right)
(\partial^{k+1-j}_{s}\mbox{\mathversion{bold}$v$},
[(\mbox{\mathversion{bold}$v$}\times \mbox{\mathversion{bold}$v$}_{s})\cdot 
\partial^{k+2}_{s}\mbox{\mathversion{bold}$v$}]
\partial^{j}_{s}\mbox{\mathversion{bold}$v$})\\[3mm]
& \qquad +
[\partial^{k+1}_{s}\mbox{\mathversion{bold}$v$}\cdot 
(\mbox{\mathversion{bold}$v$}\times \partial^{k+2}_{s}\mbox{\mathversion{bold}$v$})]^{1}_{s=0}
+
(\partial^{k+1}_{s}\mbox{\mathversion{bold}$v$},\mbox{\mathversion{bold}$h$}_{k+1})\\[3mm]
& \leq 
C\|\mbox{\mathversion{bold}$v$}_{s}\|^{2}_{k}
-
\frac{k}{2}\sum^{k}_{j=1}
\left(
\begin{array}{c}
k+1\\
j
\end{array}
\right)
(\partial^{k+1-j}_{s}\mbox{\mathversion{bold}$v$},
[(\mbox{\mathversion{bold}$v$}\times \mbox{\mathversion{bold}$v$}_{s})\cdot 
\partial^{k+2}_{s}\mbox{\mathversion{bold}$v$}]
\partial^{j}_{s}\mbox{\mathversion{bold}$v$})\\[3mm]
& \qquad +
[\partial^{k+1}_{s}\mbox{\mathversion{bold}$v$}\cdot 
(\mbox{\mathversion{bold}$v$}\times \partial^{k+2}_{s}\mbox{\mathversion{bold}$v$})]^{1}_{s=0},
\end{align*}
where \( C>0\) depends on \( C_{\ast}\).
Taking the limit \( \varepsilon \to +0 \) in Lemma \ref{vc2m}, we see that 
\( \mbox{\mathversion{bold}$v$}\) satisfies
\begin{align}
\label{1}
\mbox{\mathversion{bold}$v$}\times \partial^{2m}_{s}\mbox{\mathversion{bold}$v$}|_{s=0,1}
&=
\mbox{\mathversion{bold}$0$},\\[3mm]
\partial^{i}_{s}\mbox{\mathversion{bold}$v$}\cdot 
\partial^{n}_{s}\mbox{\mathversion{bold}$v$}|_{s=0,1}
&=
0 \quad (i+n = 2m+1),
\label{2}
\end{align}
for \( m \geq 1\). Hence, the boundary term 
\( [\partial^{k+1}_{s}\mbox{\mathversion{bold}$v$}\cdot 
(\mbox{\mathversion{bold}$v$}\times \partial^{k+2}_{s}\mbox{\mathversion{bold}$v$})]^{1}_{s=0}\)
is zero from (\ref{1})
because either \( k+1 \) or \( k+2 \) is even.
This shows that
\begin{align*}
\frac{1}{2}\frac{{\rm d}}{{\rm d}t}
\| \partial^{k+1}_{s}\mbox{\mathversion{bold}$v$}\|^{2}
& \leq
C\|\mbox{\mathversion{bold}$v$}_{s}\|^{2}_{k}
-
\frac{k}{2}\sum^{k}_{j=1}
\left(
\begin{array}{c}
k+1\\
j
\end{array}
\right)
(\partial^{k+1-j}_{s}\mbox{\mathversion{bold}$v$},
[(\mbox{\mathversion{bold}$v$}\times \mbox{\mathversion{bold}$v$}_{s})\cdot 
\partial^{k+2}_{s}\mbox{\mathversion{bold}$v$}]
\partial^{j}_{s}\mbox{\mathversion{bold}$v$}).
\end{align*}
We continue as
\begin{align*}
\frac{1}{2}\frac{{\rm d}}{{\rm d}t}
\| \partial^{k+1}_{s}\mbox{\mathversion{bold}$v$}\|^{2}
& \leq
C\|\mbox{\mathversion{bold}$v$}_{s}\|^{2}_{k}
+
\frac{k}{2}\sum^{k}_{j=1}
\left(
\begin{array}{c}
k+1\\
j
\end{array}
\right)
\bigg\{
(\partial^{k+j}_{s}\mbox{\mathversion{bold}$v$},
[(\mbox{\mathversion{bold}$v$}\times \mbox{\mathversion{bold}$v$}_{s})\cdot 
\partial^{k+1}_{s}\mbox{\mathversion{bold}$v$}]
\partial^{j}_{s}\mbox{\mathversion{bold}$v$})\\[3mm]
& \qquad +
(\partial^{k+1-j}_{s}\mbox{\mathversion{bold}$v$},
[(\mbox{\mathversion{bold}$v$}\times \mbox{\mathversion{bold}$v$}_{ss})\cdot
\partial^{k+1}_{s}\mbox{\mathversion{bold}$v$}]
\partial^{j}_{s}\mbox{\mathversion{bold}$v$})
+
(\partial^{k+1-j}_{s}\mbox{\mathversion{bold}$v$},
[(\mbox{\mathversion{bold}$v$}\times \mbox{\mathversion{bold}$v$}_{s})\cdot 
\partial^{k+1}_{s}\mbox{\mathversion{bold}$v$}]
\partial^{j+1}_{s}\mbox{\mathversion{bold}$v$})\\[3mm]
& \qquad 
-
[(\partial^{k+1-j}_{s}\mbox{\mathversion{bold}$v$}\cdot 
\partial^{j}_{s}\mbox{\mathversion{bold}$v$})
\big(
(\mbox{\mathversion{bold}$v$}\times \mbox{\mathversion{bold}$v$}_{s})\cdot 
\partial^{k+1}_{s}\mbox{\mathversion{bold}$v$}\big) ]^{1}_{s=0}\bigg\}\\[3mm]
& \leq 
C\|\mbox{\mathversion{bold}$v$}_{s}\|^{2}_{k}
-
\frac{k}{2}\sum^{k}_{j=1}
\left(
\begin{array}{c}
k+1\\
j
\end{array}
\right)
[(\partial^{k+1-j}_{s}\mbox{\mathversion{bold}$v$}\cdot 
\partial^{j}_{s}\mbox{\mathversion{bold}$v$})
\big(
(\mbox{\mathversion{bold}$v$}\times \mbox{\mathversion{bold}$v$}_{s})\cdot 
\partial^{k+1}_{s}\mbox{\mathversion{bold}$v$}\big) ]^{1}_{s=0}.
\end{align*}
When \( k=2m-1 \), we have
\begin{align*}
&\frac{1}{2}\frac{{\rm d}}{{\rm d}t}
\| \partial^{2m}_{s}\mbox{\mathversion{bold}$v$}\|^{2}\\[3mm]
& \quad \leq
C\|\mbox{\mathversion{bold}$v$}_{s}\|^{2}_{2m-1}
-
\frac{(2m-1)}{2}\sum^{2m-2}_{j=1}
\left(
\begin{array}{c}
2m\\
j
\end{array}
\right)
[(\partial^{2m-j}_{s}\mbox{\mathversion{bold}$v$}\cdot 
\partial^{j}_{s}\mbox{\mathversion{bold}$v$})
\big(
(\mbox{\mathversion{bold}$v$}\times \mbox{\mathversion{bold}$v$}_{s})\cdot 
\partial^{2m}_{s}\mbox{\mathversion{bold}$v$}\big) ]^{1}_{s=0}.
\end{align*}
From (\ref{1}), we see that 
\begin{align*}
(\mbox{\mathversion{bold}$v$}\times \mbox{\mathversion{bold}$v$}_{s})\cdot 
\partial^{2m}_{s}\mbox{\mathversion{bold}$v$}|_{s=0,1}
=
-(\mbox{\mathversion{bold}$v$}\times \partial^{2m}_{s}\mbox{\mathversion{bold}$v$})
\cdot \mbox{\mathversion{bold}$v$}_{s}|_{s=0,1}
=
0,
\end{align*}
and we have
\begin{align*}
\frac{1}{2}\frac{{\rm d}}{{\rm d}t}
\| \partial^{2m}_{s}\mbox{\mathversion{bold}$v$}\|^{2}
\leq 
C\| \mbox{\mathversion{bold}$v$}_{s}\|^{2}_{2m-1}.
\end{align*}
When \( k=2m\), we have
\begin{align*}
&\frac{1}{2}\frac{{\rm d}}{{\rm d}t}\| 
\partial^{2m+1}_{s}\mbox{\mathversion{bold}$v$}\|^{2}\\[3mm]
& \quad \leq
C\|\mbox{\mathversion{bold}$v$}_{s}\|^{2}
-
m\sum^{2m}_{j=1}
\left(
\begin{array}{c}
2m+1\\
j
\end{array}
\right)
[(\partial^{2m+1-j}_{s}\mbox{\mathversion{bold}$v$}\cdot 
\partial^{j}_{s}\mbox{\mathversion{bold}$v$})
\big(
(\mbox{\mathversion{bold}$v$}\times \mbox{\mathversion{bold}$v$}_{s})\cdot 
\partial^{2m+1}_{s}\mbox{\mathversion{bold}$v$}\big) ]^{1}_{s=0}.
\end{align*}
We see from (\ref{2}) that 
\begin{align*}
\partial^{2m+1-j}_{s}\mbox{\mathversion{bold}$v$}\cdot 
\partial^{j}_{s}\mbox{\mathversion{bold}$v$}|_{s=0,1}
=
0,
\end{align*}
for all \( j\) with \( 1\leq j \leq 2m \), yielding
\begin{align*}
\frac{1}{2}\frac{{\rm d}}{{\rm d}t}
\| \partial^{2m+1}_{s}\mbox{\mathversion{bold}$v$}\|^{2}
& \leq
C\| \mbox{\mathversion{bold}$v$}_{s}\|^{2}_{2m}.
\end{align*}
In either cases, we have
\begin{align*}
\frac{1}{2}\frac{{\rm d}}{{\rm d}t}
\| \partial^{k+1}_{s}\mbox{\mathversion{bold}$v$}\|^{2}
& \leq
C\| \mbox{\mathversion{bold}$v$}_{s}\|^{2}_{k},
\end{align*}
which, along with the assumption of induction and 
Gronwall's inequality, allows us to conclude that 
\begin{align*}
\sup_{0\leq t\leq T}
\| \mbox{\mathversion{bold}$v$}_{s}\|_{k} \leq C_{\ast},
\end{align*}
where \( C_{\ast}>0 \) depends on \( 
\| \mbox{\mathversion{bold}$v$}_{0s}\|_{k}\) and is 
monotone increasing in \( T>0\). The \( t\) derivatives can be 
estimated from the equation in (\ref{vslant}) and thus, we have derived 
an a priori estimate in the function space stated in 
Theorem \ref{TH}. 

\medskip

Combining Proposition \ref{id1} and \ref{id2}, the time-local existence 
of a smooth solution to (\ref{vslant}), and the a priori estimate just obtained,
we can conclude that by a standard approximation and continuation argument, for an
arbitrary \( T>0\) we have a 
solution \( \mbox{\mathversion{bold}$v$}\) of (\ref{vslant}) satisfying
\begin{align*}
\mbox{\mathversion{bold}$v$}\in \bigcap^{[\frac{l+3}{2}]}_{j=0}
C^{j}\big( [0,T];H^{l+2-2j}(I)\big), \qquad 
\mbox{\mathversion{bold}$v$}\in \bigcap ^{[\frac{l+3}{2}]}_{j=0}
W^{j,\infty}\big(0,T;H^{l+3-2j}(I)\big),
\end{align*}
with initial datum \( \mbox{\mathversion{bold}$v$}_{0}\in 
H^{l+3}(I)\). The uniqueness of the solution is a consequence of a standard energy estimate of 
the differenece of two solutions with the same initial datum.
Finally, since problem (\ref{vslant}) can be solved reverse in time, 
the continuity with respect to \( t\) can be recovered by the same 
argument given in Kato \cite{18}, and this proves Theorem \ref{TH}.

\vspace*{1cm}
\noindent
Masashi Aiki\\
Department of Mathematics\\
Faculty of Science and Technology, Tokyo University of Science\\
2641 Yamazaki, Noda, Chiba 278-8510, Japan\\
E-mail: aiki\verb|_|masashi\verb|@|ma.noda.tus.ac.jp

\end{document}